\documentclass[a4paper,reqno]{amsart}
\usepackage{epsfig}
\usepackage[dvips]{color}
\usepackage{amsmath,amsthm,
amssymb,bbm,mathrsfs,stmaryrd}
\usepackage[all,2cell]{xy}
\xyoption{v2}
\UseAllTwocells
\theoremstyle{plain}
\newtheorem{thm}{Theorem}[section]

\newtheorem{lem}[thm]{Lemma}
\newtheorem{cor}[thm]{Corollary}
\newtheorem{prop}[thm]{Proposition}

\theoremstyle{definition}
\newtheorem{obs}[thm]{Observation}
\newtheorem{defn}{Definition}[section]
\theoremstyle{remark}
\newtheorem{exmp}{Example}[section]

\newcommand{\M}{\ensuremath{\mathscr{M}}}

\newcommand{\Cat}{\ensuremath{\mathbf{Cat}}}
\newcommand{\n}{\ensuremath{\mathrm{n}}}
\newcommand{\e}{\ensuremath{\mathrm{e}}}
\newcommand{\A}{\ensuremath{\mathscr{A}}}
\newcommand{\Ac}{\ensuremath{A^{\circ}}}
\newcommand{\EM}{\ensuremath{\mathsf{EM}}}
\newcommand{\V}{\ensuremath{\mathscr{V}}}
\newcommand{\K}{\ensuremath{\mathscr{K}}}
\newcommand{\B}{\ensuremath{\mathscr{B}}}

\newcommand{\C}{\ensuremath{\mathscr{C}}}

\newcommand{\X}{\ensuremath{\mathscr{X}}}
\newcommand{\SorCh}{Section}

\begin{document}

\title[Hopf modules for autonomous pseudomonoids]
{Hopf modules for 
autonomous pseudomonoids and the monoidal centre}
\author{Ignacio L. Lopez Franco}
\thanks{The author acknowledges the support of an Internal Graduate 
Studentship form Trinity College, Cambridge}
\address{
Department of Pure Mathematics and Mathematical Statistics\\
Centre for Mathematical Sciences\\Wilberforce Road\\Cambridge\\CB3 0WB\\UK}
\email{I.Lopez-Franco@dpmms.cam.ac.uk}

\begin{abstract}
In this work we develop some aspects of the theory of Hopf algebras to the context of
autonomous map pseudomonoids. We concentrate in the Hopf modules and
the Centre or Drinfel'd double. If $A$ is a map pseudomonoid in a monoidal bicategory \M,
the analogue
of the category of Hopf modules for $A$ is an Eilenberg-Moore
construction for a
certain monad in $\mathbf{Hom}(\M^{\mathrm{op}},\mathbf{Cat})$. We study the
existence of the internalisation of this notion, called the Hopf
module construction,  by extending the
completion under Eilenberg-Moore objects of a 2-category to a
endo-homomorphism of tricategories on $\mathbf{Bicat}$. 

Our main result is
the equivalence between the existence of a left dualization for $A$
({\em i.e.}, $A$ is left autonomous) and the validity of an analogue
of the structure theorem of Hopf modules. In this case the Hopf module
construction for $A$ always exists. 

We use these results to study the
lax centre of a left autonomous map pseudomonoid. We show that the lax
centre is the Eilenberg-Moore construction for a certain monad on
$A$ (one existing if the other does). If $A$ is also right autonomous,
then the lax centre equals the centre. We look at the examples of the
bicategories of \V-modules and of comodules in \V, and obtain the Drinfel'd 
double of a coquasi-Hopf algebra $H$ as the centre of $H$. 

\end{abstract}

\maketitle
\section{Introduction}

This work addresses the problem of extending the basic results of the
the theory of Hopf algebras to the context of autonomous
pseudomonoids. 
We will focus mainly on two constructions: Hopf modules and the
Drinfel'd double of a Hopf algebra. 


Left autonomous pseudomonoids, introduced in \cite{dualizations},
 generalise not only Hopf and
(co)quasi-Hopf algebras but also (pro)monoidal enriched categories.
In fact,  this is the conceptual
reason underlying the well-known fact that the category of
finite-dimensional (co)representations of a (co)quasi-Hopf algebra is
left autonomous. 

Our starting point is the so-called {\em theorem of Hopf modules}\/
for (co)quasi-Hopf algebras
\cite{math.QA/9904164,Schauenburg:TwoCharacterizations}, 
that extends the classical result for ordinary Hopf algebras
\cite{Larson-Sweedler:AnAssocOrtogonal}.  A coquasibialgebra $H$,
although not associative in $\mathbf{Vect}$, is an associative
algebra in the category $\mathbf{Comod}(H,H)$ of $H$-bicomodules and
thus we can consider the category of left $H$-modules in
$\mathbf{Comod}(H,H)$. This is by definition the category of $H$-Hopf
modules. There is a monoidal functor from the category of right
$H$-comodules to the category of $H$-Hopf modules sending $M$ to the
tensor product bicomodule 
$H\otimes M$, where $M$ is considered as trivial $H$-comodule on the left. 
It is shown in \cite{Schauenburg:TwoCharacterizations} that
when $H$ is a coquasi-Hopf algebra this functor is an equivalence, and
in a dual fashion, that a finite-dimensional quasibialgebra is
quasi-Hopf if and only if the module version of this functor is an
equivalence. 

We prove that an analogous result holds if we replace
coquasibialgebras by map pseudomonoids ({\em i.e.,}  pseudomonoids
whose multiplication and unit have a right adjoint), Hopf modules by the
Eilenberg-Moore construction for certain monad and coquasi-Hopf
algebras by left autonomous map pseudomonoids. Moreover, in our general
setup no finiteness condition is necessary. We take this as an
indication that the concept of dualization is  more natural than the
one of antipode.

When the monoidal bicategory involved is right closed, and in
particular when it is right autonomous, our Hopf module construction
can be internalised. Naturally, this internalisation need not exist,
being an Eilenberg-Moore construction for some monad on the
endo-hom object $[A,A]$ of the map pseudomonoid $A$. However, it does
exist when the map pseudomonoid is left autonomous, and its object
part is equivalent to $A$.

The centre of a monoidal category was defined in \cite{braidedtencat},
and more recently generalised to the {\em centre construction}\/ for
pseudomonoids \cite{Street:monoidalcentre}. 
A classical result reads that, for a Hopf algebra $H$, the category of
two-sided Hopf modules is monoidally 
equivalent to the centre of the category of
$H$-(co)modules; this has been extended to the case of quasi-Hopf algebras in
\cite{Schauenburg:HopfMods}. Both versions use the category of
Yetter-Drinfel'd modules as an intermediate stage in proving the
equivalence. 

In our general context, the lax centre
of a left autonomous map pseudomonoid is the Eilenberg-Moore construction
for a certain monad on the pseudomonoid, and if this is also right
autonomous, then the lax centre is in fact the centre. In this way we
reduce the problem of the existence of the lax centre of such
pseudomonoids to the problem of the existence of a particular
Eilenberg-Moore object. 
No analogue of
the Yetter-Drinfel'd modules appear in our construction. If we think
of a Hopf algebra as an autonomous pseudomonoid in the appropriate monoidal
bicategory, its centre is equivalent to the {\em Drinfel'd double}\/
of the Hopf algebra.  This shows that the centre construction
generalises the classical quantum double construction for Hopf
algebras. When we apply our results to the bicategory of \V-modules,
we are able to show that any left autonomous map pseudomonoid has a
lax centre. In particular, any left autonomous monoidal \V-category
has a lax centre in $\V\text-\mathbf{Mod}$.

We shall describe the content of each section. 

In Section \ref{thmhopf.s} we introduce 
the Hopf modules for a map pseudomonoid $A$ in a monoidal bicategory $\M$
as the Eilenberg-Moore construction for a certain monad in
$[\M^{\mathrm{op}},\mathbf{Cat}]$,  and  explain what we 
mean by the theorem of Hopf modules. 

Section \ref{opmon.s}  surveys some well-known facts about lax
actions and opmonoidal morphisms.

When the monad in the
definition of Hopf modules  is
representable by a monad $t:[A,A]\to [A,A]$ in \M, we call an 
Eilenberg-Moore construction for it a {\em Hopf module construction}\/ for
$A$. This is introduced in Section \ref{inthopf.s} along with the
proof that $t$ is a opmonoidal monad.
A Hopf module  construction, of course, need not exist in
general. 

In Section \ref{existence.s} we use completions under
Eilenberg-Moore objects to study the existence of Hopf module
constructions. To that end, we extend these completions to a
$\mathbf{Gray}$-functor on $\mathbf{Gray}$, and then to a homomorphism
of tricategories on $\mathbf{Bicat}$.

In Section \ref{maintheorem.s} we prove our main result: a
map pseudomonoid $A$ is left autonomous if and only if the theorem of Hopf
modules holds for $A$. Also, we use the results of the preceding
section to show that a map pseudomonoid is left autonomous if and
only if it has  
a Hopf module construction of a particular form, relating the 
problem of the existence of a dualization with a
completeness problem.

In Section \ref{s:Frob} we look at the relationship between autonomous
and Frobenius pseudomonoids.  

Section \ref{centre.s} deals with the relation between the Hopf module
construction and the lax centre construction. We show that for a left
autonomous map pseudomonoid $A$ there exists a opmonoidal monad on $A$
for which an Eilenberg-Moore construction is exactly the lax centre
of $A$. Moreover, if $A$ is also right autonomous then the lax centre
coincides with   the centre of $A$, and one exists if and only if the
other does.

The last two sections are devoted to the applications. 

In Section \ref{VMod.sec} we treat the example of the bicategory
$\V\text{-}\mathbf{Mod}$ of
\V-modules, for a complete and cocomplete closed symmetric monoidal
category \V. As shown in \cite{dualizations}, a left autonomous
pseudomonoid in $\V\text{-}\mathbf{Mod}$ whose multiplication, unit
and dualization are representable by \V-functors is exactly a left
autonomous \V-category in the usual sense that every object has a left
dual. The bicategory $\V\text{-}\mathbf{Mod}$ has 
Eilenberg-Moore objects for monads
and therefore any map pseudomonoid in it has a Hopf module construction
and any left autonomous map pseudomonoid has a lax centre, 
of which an explicit description is provided. We also exhibit for a
promonoidal \V-category \A, a
canonical equivalence of \V-categories $[Z_\ell(\A),\V]\simeq
Z_\ell([\A,\V])$; here the lax centre on the left hand side is the one
of \A\ as a pseudomonoid in $\V\text{-}\mathbf{Mod}$, while on the
right hand side we have the lax centre in $\V\text{-}\mathbf{Cat}$. 

The example of the bicategory of comodules over a braided monoidal
category is studied in the Section \ref{comodules.sec}. Along
with explicit descriptions of the general constructions of the
previous sections, we show how our work generalises results in
\cite{Schauenburg:TwoCharacterizations} on the theorem of Hopf
modules.
For example, a
coquasibialgebra $C$ has a dualization (which in this case is a left
$C^{\mathrm{cop}}$ and right $C$-bicomodule) if and only if the theorem
of Hopf modules holds for $C$. 
Thus, admitting dualizations instead of mere antipodes
$C^{\mathrm{cop}}\to C$, we are able to drop the finiteness conditions
on C required in \cite{Schauenburg:TwoCharacterizations}. 
We also show that the Hopf module construction always exists for
finite dimensional coquasibialgebras, and that our notion of the theorem
of Hopf modules, which implies 
the one of \cite{Schauenburg:TwoCharacterizations}, is in fact
equivalent to it when $C$ is finite dimensional. Finally, by means of
the results in Section \ref{centre.s}, we show that if  $C$ is a finite
coquasi-Hopf algebra, then the centre of $C$ 
exists and is equivalent to the Drinfel'd double of $C$. 

While preparing a definitive version of this work the author recieved
the preprint \cite{Street:DoublesMonCatsPreprint}, where Pastro and
Street study what they call doubles of monoidal enriched
categories. They are lead to structures in the bicategory of
\V-modules very similar to some of the ones we describe in Section
\ref{VMod.sec}. 

The author would like to thank Martin Hyland for his support
throughout the preparation of this work and reading several versions
of it. Also thanks to Ross Street for 
useful correspondence and spotting some mistakes in a previous version. 

\section{Preliminaries on pseudomonoids}\label{s:prelimin}
In this section we introduce the analogue for a map pseudomonoid of
the category of Hopf modules. 

We call {\em maps}\/ the 1-cells in a bicategory with right adjoint. 

Recall that a Gray monoid \cite{monbicat} 
is a monoid in the monoidal category $\mathbf{Gray}$. For a detailed
treatment of the category $\mathbf{Gray}$ see \cite{tricategories}.
As a category, $\mathbf{Gray}$ is just the
category of 2-categories and 2-functors. However, the
monoidal structure we are interested in is not the one given by the
cartesian product. 

A {\em cubical functor}\/ in two variables is a
pseudofunctor between 2-categories $F:\mathscr K\times\mathscr
L\to\mathscr N$ with the following property: when we fix one of the
variables we get a 2-functor, and  for any 1-cell $(f,g)$ in
$\mathscr K\times\mathscr L$ the constraint $F(1,g)F(f,1)\cong F(f,g)$
is equal to the identity. For any pair of 2-categories $\mathscr K$ and
$\mathscr L$, there is a 2-category $\mathscr K\boxempty\mathscr L$
equipped with a cubical functor $\mathscr K\times\mathscr L\to\mathscr
K\boxempty\mathscr L$ inducing a bijection between cubical functors
$\mathscr K\times\mathscr L\to\mathscr N$ and 2-functors $\mathscr
K\boxempty\mathscr L\to\mathscr N$. This defines a monoidal
structure on $\mathbf{Gray}$, which moreover is symmetric with the
symmetry induced by the usual one for the cartesian product. 
A Gray monoid is the same as a one-object $\mathbf{Gray}$-category in
the sense of enriched category theory, and therefore it can be thought
of as a one-object tricategory, that is, a monoidal bicategory (see
\cite{tricategories}). By
the coherence theorem in \cite{tricategories}, any monoidal bicategory
is monoidally biequivalent (that is, triequivalent as a tricategory)
to a Gray monoid.  This allows us to develop the general theory by using
Gray monoids instead of general monoidal bicategories. 

Our main examples of monoidal bicategories will be the bicategory of
comodules $\mathbf{Comod}(\V)$ in a monoidal category \V (see Section
\ref{comodules.sec}) and the 
bicategory of \V-modules $\V\text{-}\mathbf{Mod}$ (see Section
\ref{VMod.sec}).

Let \M\ be a Gray monoid
and fix a map pseudomonoid $(A, j, p)$
in  \M, that is, a pseudomonoid whose unit $j:I\to A$ and
multiplication $p:A\otimes A\to A$ are maps. Recall from
\cite{monbicat} that a pseudomonoid, in addition to the unit and
multiplication, is equipped with isomorphisms $\phi:p(p\otimes
A)\Rightarrow p(A\otimes p)$, $p(j\otimes A)\Rightarrow 1_A$ and
$p(A\otimes j)\Rightarrow 1_A$ satisfying three axioms which ensure,
as shown in \cite{Lack:Acoherentapproach}, that any 2-cell formed by
pasting of tensor products of these isomorphisms, 1-cells and
pseudonaturality constraints of the Gray monoid is
uniquely determined by its domain and codomain 1-cells.


If $(A,j,p)$ is a map pseudomonoid, then $(A,j^*,p^*)$ is a {\em
pseudocomonoid}, that is, a pseudomonoid in the opposite Gray monoid.
By definition the unit isomorphism $(A\otimes j^*)p^*\cong 1_A$ of the
pseudocomonoid $(A,j^*,p^*)$ is the mate
of the constraint $p(A\otimes j)\cong 1_A$, and thus the
following equality holds.
\begin{equation}\label{eq:jyp}
\xymatrixcolsep{.7cm}
\diagramcompileto{jyp}
A\ar@{=}[rr]\ar[dr]_{p^*}&{}\dtwocell<\omit>{'\cong}&A\ar[dr]^{1\otimes j}&\\
&A^2\ar[ur]|-{1\otimes j^*}\ar@{=}[rr]\rrtwocell<\omit>{<-2.5>}&&
A^2
\enddiagram
=
\diagramcompileto{jyp2}
&A^2\ar@{=}[rr]\rrtwocell<\omit>{<2.5>}\ar[dr]^{p}&&
A^2\\
A\ar@{=}[rr]\ar[ur]^-{1\otimes j}&{}\utwocell<\omit>{'\cong}&A\ar[ur]_{p^*}
\enddiagram
\end{equation}
We mention this because it will be useful in Section \ref{maintheorem.s}.

In order to give a concise and conceptual definition of the Hopf
modules in the next section, we will need to use the Kleisli bicategory
of a pseudocomonad.

One can define a pseudomonad on the 2-category \K\ as a pseudomonoid
in the  monoidal 2-category $\mathbf{Hom}(\K,\K)$ of pseudofunctors,
pseudonatural transformations and modifications. A pseudocomonad is a
pseudocomonoid in the same monoidal 2-category. As before, if $T$ is a
pseudomonad with unit $\eta:1\Rightarrow T$ and multiplication
$\mu:T^2\Rightarrow T$ which are maps, then $T$ together with $\eta^*$
and $\mu^*$ have a canonical structure of a pseudocomonad on \K. 

A lax $T$-algebra is an arrow $a:TA\to A$ in \K\ equipped with a
2-cell $a(Ta)\Rightarrow a\mu_A:T^2A\to A$ satisfying the axioms in
\cite[p. 39]{Marmolejo:Doctrines} and \cite{Lack:Acoherentapproach},
but without the requirement of the invertibility of the 2-cell. 

Let $G$ be a pseudocomonad on the 2-category $\K$, and denote its
comultiplication and counit by $\delta$ and $\epsilon$,
respectively. The Kleisli bicategory $\mathsf{Kl}(T)$ of \K\ has the
same objects as \K, and hom-categories
$\mathsf{Kl}(T)(X,Y)=\K(GX,Y)$. We denote the 1-cells of
$\mathsf{Kl}(T)$ by $f:X \nrightarrow Y$. The composition of this $f$
with $g:Y\nrightarrow Z$ is given by $g (Gf)\delta_X:GX\to Z$, while
the identity of the object $X$ is $\varepsilon_X:GX\to X$. 

The following is a slight generalisation of part of
\cite[Prop. 4.6]{Hermida:FromCoherent}.

\begin{lem}\label{l:oplax=monad}
  Let $T:\K\to\K$ be a pseudomonad whose unit $\eta$ and multiplication $\mu$
  are maps. Then, there exists a bijection between the following
  structures on an arrow $a:TA\to A$ in \K: structures of an lax
  $T$- algebra and structures of a monad  in $\mathsf{Kl}(T)$. 
\end{lem}

A structure of a monad in $\mathsf{Kl}(T)$ on $a:A\nrightarrow A$ is
given by a pair of 2-cells $a (Ta)\mu_A^*\Rightarrow a$ and
$\eta_A^*\Rightarrow a$ in \K. The bijection above is given by 
$${
\diagramcompileto{oplaxmonad1}
T^2A\ar[r]^-{Ta}\ar[d]_{\mu_A}\drtwocell<\omit>{}&
TA\ar[d]^a\\
TA\ar[r]_-a&
A
\enddiagram}
\mapsto
{\diagramcompileto{oplaxmonad2}
&T^2A\ar[r]^-{Ta}\ar[d]_{\mu_A}\drtwocell<\omit>{}&
TA\ar[d]^a\\
TA\ar[ur]^{\mu_A^*} \ar@{=}[r]\rtwocell<\omit>{<-2>} &TA\ar[r]_-a&
A
\enddiagram}
$$

\section{The theorem of Hopf modules}\label{thmhopf.s}
If $(A,j,p)$ is a map pseudomonoid in the Gray monoid \M, the 2-functor
$A\otimes-$ has the structure of a pseudomonad with unit $j\otimes
X:X\to A\otimes X$ and multiplication $p\otimes X:A\otimes A\otimes
X\to A\otimes X$, and also the structure of a pseudocomonad with counit
$j^*\otimes X$ and comultiplication $p^*\otimes X$.
The associativity constraint $p(A\otimes p)\Rightarrow
p(p\otimes A)$ endows $p:A\otimes A\to A$ with the structure of a
lax $(A\otimes-)$-algebra, and hence by Lemma \ref{l:oplax=monad},
with the structure of a monad $p:A\nrightarrow A$ in the Kleisli
bicategory $\mathsf{Kl}(A\otimes -)$. 

\begin{defn}\label{d:theta}
  We will denote by $\theta$ the monad $\mathsf{Kl}(A\otimes-)(-,p)$
  in $\mathsf{Kl}(A\otimes -)$. Hence,
  $\theta$ is a monad on the 2-functor $\M(A\otimes-,A)$ in the
  2-category 
  $\mathbf{Hom}(\M^{\mathrm{op}},\mathbf{Cat})$ of pseudofunctors,
  pseudonatural transformations and modifications. 
\end{defn}

Explicitly, $\theta_X(f)=p(A\otimes f)(p^*\otimes X)$ and  the
multiplication and unit of the monad, depicted in \eqref{mult.eq} and
\eqref{unit.eq}, are induced by the counits of the
adjunctions $p\dashv p^*$ and $j\dashv j^*$ respectively.
\begin{equation}\label{mult.eq}
{\xymatrixcolsep{1cm}
\diagramcompileto{mult}
A\otimes X\ar[r]^-{p^*\otimes1}\ar[dr]_{p^*\otimes1}&
A^2\otimes X\ar[r]^-{1\otimes p^*\otimes1}&
A^3\otimes X\ar[r]^-{A^2\otimes f}\ar[dr]|-{p\otimes1\otimes1}&
A^3\ar[r]^-{1\otimes p}\ar[dr]|-{p\otimes1}&
A^2\ar[r]^-p&
A\\
&A^2\otimes X\ar[ur]|-{p^*\otimes1\otimes1}\ar@{=}[rr]
\rrtwocell<\omit>{<-2.5>}\ar@{}[u]|-{\phi^*\cong}&&
A^2\ar[r]_-{1\otimes f}\ar@{}[u]|-{c\ \cong}&
A^2\ar[ur]_{p}\ar@{}[u]|-{\cong}
\enddiagram}
\end{equation}
\begin{equation}\label{unit.eq}
{\diagramcompileto{unit}
&&A\otimes X\drrtwocell<\omit>{'\cong\ c}
\ar[r]^-f\ar[dr]|-{j\otimes 1\otimes1}&
A\ar[dr]|-{j\otimes 1}\ar@/^1pc/[drr]^1\drrtwocell<\omit>{'\cong}&\\
A\otimes X\ar@/^1pc/[urr]^1\urrtwocell<\omit>{'{\cong}}\ar[r]_-{p^*\otimes1}&
A^2\otimes X\ar[ur]|-{j^*\otimes 1\otimes 1}\ar@{=}[rr]
\rrtwocell<\omit>{<-2.5>}&&
A^2\otimes X\ar[r]_-{1\otimes f}&
A^2\ar[r]_-p&
A
\enddiagram}
\end{equation}

Our generalisation of the category of Hopf modules is 
the Eilenberg-Moore construction $\upsilon:\M(A\otimes
-,A)^\theta\to \M(A\otimes -,A)$ for the monad $\theta$ in
$\mathbf{Hom}(\M^{\mathrm{op}},\Cat)$. We denote by $\varphi$ the left adjoint of
$\upsilon$. As we explain in \SorCh \ref{comodules.sec}, when $A$ is
the pseudomonoid in $\mathbf{Comod}(\mathbf{Vect})$ induced by a
coquasi-bialgebra, $\mathbf{Comod}(\mathbf{Vect})(A,A)^{\theta_{I}}$
is the category of Hopf modules described in
\cite{Schauenburg:TwoCharacterizations}. In particular, when $A$ is
induced by a bialgebra, we get the usual category of Hopf modules
(simultaneous $A$-bicomodules and left $A$-modules, plus compatibility
between the two structures).

\begin{defn}\label{lambda.def}
We say that {\em the theorem of Hopf modules holds}\/ for a map 
pseudomonoid $A$
if the pseudonatural transformation $\lambda$ given by
$$
\M(-,A)\xrightarrow{\M(j^*\otimes-,A)}\M(A\otimes-,A)\xrightarrow{\varphi}
\M(A\otimes-,A)^\theta 
$$
is an equivalence.
\end{defn}
\begin{obs}\label{upsilonlambda.obs}
The composition
$\upsilon_X\lambda_X=\theta_X\M(j^*\otimes X,A):\M(X,A)\to\M(A\otimes
X,A)$ is, up to isomorphism, the functor given by
$$
(X\xrightarrow{f}A)\longmapsto (A\otimes X\xrightarrow{1\otimes
f}A\otimes A\xrightarrow{p}A).
$$

\end{obs}

Recall that a 1-cell in a bicategory is  {\em fully faithful}\/
if it is a map and the unit of the adjunction is an isomorphism.

\begin{prop}\label{lambdaff.prop}
The pseudonatural transformation $\lambda$ is fully faithful.
\end{prop}
\begin{proof}
It is clear that $\lambda$ has right adjoint $\M(j\otimes
-,A)\upsilon$. The component at $X$ of the 
unit of this adjunction is the composite 
$1\to\M(j\otimes X, A)\M(j^*\otimes X,A)\to \M(j\otimes X,
A)\upsilon_X\varphi_X\M(j^*\otimes X,A)$, where each arrow is induced
by the corresponding unit. Evaluation of this natural transformation
at $f:X\to A$ gives the pasted composite (where the unlabelled 2-cells
denote the obvious counits)
$$
\diagramcompileto{fullyfaithful1}
X\ar[d]_{j\otimes1}\ar@(ur,ul)[rrr]^1\rrtwocell<\omit>{<-1>}&&
A\otimes X\ar[dr]|-{j\otimes1\otimes 1}\ar[r]_-{j^*\otimes 1}&
X\ar[r]^-{f}\ar[dr]|-{j\otimes 1}\dtwocell<\omit>{'\cong}&
A\ar[dr]|-{j\otimes1}\ar@/^1pc/[drr]^1\drrtwocell<\omit>{'\cong}
\dtwocell<\omit>{'\cong}&&\\
A\otimes
X\ar[r]_-{p^*\otimes1}\ar@/^1pc/[urr]|-1\urrtwocell<\omit>{'\cong}&
A^2\otimes
X\ar[ur]|-{j^*\otimes1\otimes1}\ar@{=}[rr]\rrtwocell<\omit>{<-2.5>} &&
A^2\otimes X\ar[r]_-{1\otimes j^*\otimes1}&
A\otimes X\ar[r]_-{1\otimes f\otimes 1}&
A^2\ar[r]_-{p}&
A
\enddiagram
$$
which by \eqref{eq:jyp} is equal to 
$$
\diagramcompileto{fullyfaithful2}
X\ar[rr]^1\rtwocell<\omit>{<2>}\ar[d]_{j\otimes1}&&
X\ar[r]^f\ar[dr]|-{j\otimes1}\dtwocell<\omit>{'\cong\ c}&
A\ar[dr]|-{j\otimes 1}\ar@/^1pc/[drr]^1\dtwocell<\omit>{'\cong\ c}
\drrtwocell<\omit>{'\cong}&&\\
A\otimes X\ar[urr]|-{j^*\otimes1}\ar[r]_-{j\otimes1\otimes1}
\drlowertwocell_1{\omit}\drtwocell<\omit>{'\cong}&
A^2\otimes X\ar@{=}[r]\ar[d]|-{p\otimes1}\rtwocell<\omit>{<2>}&
A^2\otimes X\ar[r]_-{1\otimes j^*\otimes1}&
A\otimes X\ar[r]_-{1\otimes f}&
A^2\ar[r]_-p&
A\\
&A\otimes X\ar[ur]_{p^*\otimes 1}
\enddiagram
$$
and thus to
$$
\diagramcompileto{fullyfaithful3}
X\ar[rrrr]^f\ar[d]_{j\otimes1}&&&&
A\ar[d]_{j\otimes1}\druppertwocell^1{\omit}\drtwocell<\omit>{'\cong}&\\
A\otimes X\ar@(ur,ul)[rrr]^1\ar[r]_{1\otimes
  j\otimes1}\drlowertwocell_1{\omit}\drtwocell<\omit>{'\cong}&
A^2\otimes X\ar[d]|-{p\otimes 1}\ar@{=}[r]
\rtwocell<\omit>{<2>}\rtwocell<\omit>{<-2>}&
A^2\otimes X\ar[r]_-{1\otimes j^*\otimes 1}&
A\otimes X\ar[r]_-{1\otimes f}\ultwocell<\omit>{'\cong}&
A^2\ar[r]_-p&
A\\
&A\otimes X\ar[ur]_{p^*\otimes 1}
\enddiagram
$$
This last pasting composite clearly is an isomorphism because of the
isomorphism $p(A\otimes j)\cong 1$.
\end{proof}

The following observation will be of use in \SorCh\
\ref{maintheorem.s}

\begin{obs}\label{obs1.obs}
Consider the following modification, where $\varepsilon $ denotes the
counit of the adjunction $\lambda\dashv \M(j\otimes -,A)\upsilon$.
$$
\diagramcompileto{counit1}
\M(A\otimes-,A)\ar[r]^-{\varphi}&
\M(A\otimes-,A)^\theta\ar[ddrr]_1\ar[r]^-\upsilon&
\M(A\otimes-,A)\ar[r]^-{\M(j\otimes-,A)}\drtwocell<\omit>{\varepsilon}&
\M(-,A)\ar[d]^{\M(j^*\otimes-,A)}\\
&&&\M(A\otimes-,A)\ar[d]^{\varphi}\\
&&&\M(A\otimes-,A)^\theta\ar[d]^\upsilon\\
&&&\M(A\otimes-,A)
\enddiagram
$$
Observe that   $\M(A\otimes X,A)^{\theta_X}$ is the closure under
$\upsilon_X$-split coequalizers 
of the full subcategory determined by the image of the
functor $\varphi_X$, and these coequalizers are preserved by 
$\varphi_X\M(jj^*\otimes-,A)\upsilon_X$, since they become absolute
coequalizers after applying $\upsilon_X$.
It follows that $\varepsilon_X$ is an isomorphism if and only
if $\varepsilon _X \varphi_X$ is an isomorphism. 
Using the fact that each $\upsilon_X$ is
conservative, we deduce that $\varepsilon $ is an isomorphism if and
only if $\upsilon\varepsilon\phi$ is so. 
\end{obs}

\begin{obs}\label{o:1monoid}
  There is another equivalent way of defining Hopf modules. The
  category $\M(A,A)$ has a {\em convolution}\/ monoidal structure,
  with tensor product $f*g=p(A\otimes g)(f\otimes A)p^*$ and unit
  $jj^*$. This monoidal category acts on the pseudofunctor $\M(A\otimes
  -,A):\M^{\mathrm{op}}\to\mathbf{Cat}$ by sending $h:A\otimes X\to A$
  to $p(A\otimes h)(p^*\otimes X)$, in the sense that this defines a
  monoidal functor from $\M(A,A)$ to
  $\mathbf{Hom}(\M^{\mathrm{op}},\mathbf{Cat})(\M(A\otimes
  -,A),\M(A\otimes -,A))$. Now, $1_A:A\to A$ has a canonical structure
  of a monoid in $\M(A,A)$, with multiplication $pp^*\Rightarrow 1$
  and $jj^*\Rightarrow 1$ the respective counits of the
  adjunctions. Hence, $1_A$ defines via the action described above, a
  monad on $\M(A\otimes -,A)$ in
  $\mathbf{Hom}(\M^{\mathrm{op}},\mathbf{Cat})$. This monad is just
  the monad $\theta$ of Definition \ref{d:theta}.
\end{obs}

\section{Opmonoidal morphisms and oplax actions}\label{opmon.s}

In this section we  spell out the relation between opmonoidal morphisms
and right oplax actions in a right closed Gray monoid.  
Everything in this section is well-known, though we have not found 
the present formulation in the literature.  The case when the
monoidal 2-category is strict 
and has certain completeness conditions is studied in 
\cite{Kelly-Lack:MonFunctors}. 

Let $A$ be a pseudomonoid in \M. Briefly, a {\em right oplax action}\/
of $A$ on an object $B$ is an oplax algebra for the pseudomonad
$-\otimes A$ on $\M$. This amounts to a 1-cell $h:B\otimes A\to B$
together with 2-cells
$$
\diagramcompileto{oplaxact1}
B\otimes A\otimes A\ar[rr]^-{h\otimes 1}\ar[d]_{1\otimes p}&&
B\otimes A\ar[d]^h\\
B\otimes A\rrlowertwocell<0>_{h}{^<-3>h_2}&&
A
\enddiagram
\qquad
\diagramcompileto{oplaxact2}
&B\ar[dl]_{1\otimes j}\ar[dr]^1&\\
B\otimes A\rrlowertwocell<0>_{h}{^<-3>h_0}&&
A
\enddiagram
$$
satisfying axioms dual to those in \cite[p. 39]{Marmolejo:Doctrines}
or \cite{Lack:Acoherentapproach} but without the invertibility
requirement on the 2-cells.  A {\em morphism of right oplax
  actions}\/ on $B$ from $h$ to $k:B\otimes A\to B$ is a 2-cell
$\tau:h\Rightarrow k$  such that 
$$
\diagramcompileto{oplaxact2}
B\otimes
A^2\rrtwocell^{k\otimes1}_{h\otimes1}{^\tau\otimes1\hole\hole}\ar[d]_{1\otimes
  p}&&
B\otimes A\dtwocell_h^k{^\tau}\\
B\otimes A\rrlowertwocell<0>_k{^<-3>h_2}&&
B
\enddiagram
=
\diagramcompileto{oplaxact2}
B\otimes A^2\ar[rr]^-{k\otimes 1}\ar[d]_{1\otimes p}&&
B\otimes A\ar[d]^k\\
B\otimes A\rrtwocell_h^k{^\tau}\rrtwocell<\omit>{^<-5>k_2}&&
B
\enddiagram
$$
$$
\diagramcompileto{oplaxact3}
&B\ar[dl]_{1\otimes j}\ar[dr]^1&\\
B\otimes A\rrlowertwocell<\omit>_{h}{^<-4.5>k_0}\rrtwocell_h^k{^\tau}&&
A
\enddiagram
=
\diagramcompileto{oplaxact4}
&B\ar[dl]_{1\otimes j}\ar[dr]^1&\\
B\otimes A\rrlowertwocell<0>_{h}{^<-3>h_0}&&
A
\enddiagram
$$

Right oplax actions of $A$
on $B$ and their morphisms form a category 
$\mathbf{Opact}_A(B)$ which  comes equipped with a canonical forgetful
functor to $\M(B\otimes A,B)$.

For each Gray monoid \M\ we have a 2-category $\mathbf{Mon}(\M)$
whose objects, 1-cells and 2-cells are respectively 
pseudomonoids in \M, lax monoidal
morphisms and monoidal 2-cells. See \cite{McCrudden:Maschkean} and
references therein. 
Define
$\mathbf{Opmon}(\M)=\mathbf{Mon}(\M^{\mathrm{co}})^{\mathrm{co}}$. The
objects of $\mathbf{Opmon}(\M)$ may be identified with the pseudomonoids, 
the 1-cells, called {\em opmonoidal morphisms}, are
1-cells $f:A\to B$ of \M\ equipped with 2-cells
$$
\diagramcompileto{opmon1}
A\otimes A\ar[rr]^-{(B\otimes f)(f\otimes A)}\ar[d]_{p}&&
B\otimes B\ar[d]^p\\
A\rrlowertwocell<0>_f{^<-3>f_2}&&
B
\enddiagram
\qquad
\diagramcompileto{opmon2}
&I\ar[dl]_j\ar[dr]^j&\\
A\rrlowertwocell<0>_f{^<-3>f_0}&&
B
\enddiagram
$$
satisfying the obvious equations, and the 2-cells $f\Rightarrow g$ 
are the 2-cells of
\M\ satisfying compatibility conditions with $f_2$, $g_2$ and $f_0$,
$g_0$.

Now suppose that $\M$ is a right closed  Gray monoid in the sense
of \cite{monbicat}, that is, 
there is a pseudofunctor $[-,-]:\M^{\mathrm{op}}\times\M\to\M$ and a
  pseudonatural equivalence
\begin{equation}\label{rightclosed}
\M(X\otimes  Y,Z)\simeq\M(X,[Y,Z]).
\end{equation}
Equivalently, 
for each pair of objects $Y$, $Z$ of \M\ there is another one denoted
by $[Y,Z]$ and an evaluation 1-cell
$\mathrm{ev}_{Y,Z}:Y\otimes[Y,Z]\to Z$ inducing \eqref{rightclosed}.
For any object $X$ of \M, the internal hom $[X,X]$ has a canonical
structure of a pseudomonoid; namely, there are composition and
identity 1-cells $\mathrm{comp}:[X,X]\otimes[X,X]\to[X,X] $ and
$\mathrm{id}:I\to[X,X]$ corresponding respectively to 
$$
X\otimes
[X,X]\otimes[X,X]\xrightarrow{\mathrm{ev}\otimes1}X\otimes[X,X]
\xrightarrow{\mathrm{ev}}X \quad\text{and}\quad
X\xrightarrow{1_X}X.
$$

\begin{prop}\label{opacteqopmon.p}
For any pseudomonoid $A$ and any object $B$,  
the closedness equivalence $\M(B\otimes A,B)\simeq\M(A,[B,B])$
lifts to an equivalence
$$
\mathbf{Opact}_A(B)\simeq\mathbf{Opmon}(\M)(A,[B,B]).
$$ 
Moreover,
under this equivalence pseudoactions correspond to pseudomonoidal
morphisms. 
\end{prop}

\begin{prop}\label{opmonmorphisms.prop}
\begin{enumerate}
\item
For any map $f:X\to Y$ the 1-cell $[f^*,f]$ from $[X,X]$ to $[Y,Y]$ has a
canonical structure of an opmonoidal morphism. If $\tau:f\Rightarrow
g$ is an invertible 2-cell then
$[(\tau^{-1})^*,\tau]:[f^*,f]\Rightarrow[g^*,g]$ is an invertible
monoidal 2-cell. 
\item
For any pair of objects $X,Y$ of \M, the 1-cell
$i_{X}^Y:[X,X]\to[Y\otimes X,Y\otimes X]$
corresponding to 
$Y\otimes \mathrm{ev}:Y\otimes X\otimes [X,X]\to Y\otimes X$
has a canonical structure of a strong monoidal morphism. Moreover,
there are canonical monoidal isomorphisms 
$(i^W_{Y\otimes X})(i^Y_{X})\cong  i^{W\otimes Y}_{X}$.
\item
For any map $f:X\to Z$ and any object $Y$ there exists a canonical
monoidal isomorphism
\begin{equation}\label{fandiiso.eq}
\diagramcompileto{fandiiso}
[X,X]\ar[r]^-{i^Y_X}\ar[d]_{[f^*,f]}\drtwocell<\omit>{'\cong}&
[Y\otimes X,Y\otimes X]\ar[d]^{[1\otimes f^*,1\otimes f]}\\
[Z,Z]\ar[r]_-{i^Y_Z}&
[Y\otimes Z,Y\otimes Z]
\enddiagram
\end{equation} 
\item
Given a map $f:Y\to Z$ and an object $X$, the counit of $f\dashv f^*$
induces a monoidal 2-cell
\begin{equation}\label{fiimon.eq}
\xymatrixcolsep{1.5cm}
\diagramcompileto{fiimon}
[X,X]\ar[r]^-{i^Y_X}\drlowertwocell<0>_{i^Z_X}{<-2>}&
[Y\otimes X,Y\otimes X]\ar[d]^{[f^*\otimes1,f\otimes1]}\\
&[Z\otimes X,Z\otimes X]
\enddiagram
\end{equation}
\end{enumerate}
\end{prop}
\begin{proof} (1)
It is not hard to show that the 2-cells \eqref{fcolaxact1.e} and
\eqref{fcolaxact2.e} equip 
$$
Y\otimes
[X,X]\xrightarrow{f^*\otimes1}X\otimes[X,X]
\xrightarrow{\mathrm{ev}}X\xrightarrow{f} Y
$$ with a structure of 
right oplax action of $[X,X]$ on $Y$, 
\begin{figure}\label{fcolaxact.f}
\begin{equation}\label{fcolaxact1.e}
\begin{split}
{\diagramcompileto{fcolaxact1}
Y\otimes[X,X]^2\ar[r]^-{f^*\otimes1\otimes1}\ar[ddd]_{1\otimes\mathrm{comp}}&
X\otimes[X,X]^2\ar[r]^-{\mathrm{ev}\otimes1}\ar[ddd]|-{1\otimes\mathrm{comp}}&
X\otimes[X,X]\ar[r]^-{f\otimes1}\ar@{=}[dr]\drtwocell<\omit>{^<-2>}&
Y\otimes [X,X]\ar[d]^{f^*\otimes1}\\
{}\drtwocell<\omit>{'\cong}&
{}\drrtwocell<\omit>{'\cong}&
&X\otimes[X,X]\ar[d]^{\mathrm{ev}}\\
&&& X\ar[d]^f\\
Y\otimes[X,X]\ar[r]_-{f^*\otimes1}&
X\otimes[X,X]\ar[r]_-{\mathrm{ev}}&
X\ar@{=}[ur]\ar[r]_-f&
Y
\enddiagram}
\end{split}
\end{equation}
\begin{equation}\label{fcolaxact2.e}
\begin{split}
{\diagramcompileto{fcolaxact2}
&Y\ar[d]_{f^*}\ar[ddl]_{1\otimes\mathrm{id}}\ar@{=}[ddrr]&&\\
&X\ar[d]_{1\otimes\mathrm{id}}\ar@{=}[dr]\drtwocell<\omit>{^<-3>}&\\
Y\otimes[X,X]\ar[r]_-{f^*\otimes1}\urtwocell<\omit>{'\cong}&
X\otimes[X,X]\ar[r]_-{\mathrm{ev}}&
X\ar[r]_-f&
Y
\enddiagram}
\end{split}
\end{equation}
\caption{}
\end{figure}
and that 
$$
\xymatrixcolsep{1.7cm}
\diagramcompileto{tauactionmorph}
Y\otimes[X,X]\rrtwocell^{f^*\otimes1}_{g^*\otimes
1}{\hole\hole\hole\hole\hole\hole(\tau^{-1})^*\otimes1}&&
X\otimes [X,X]\ar[r]^-{\mathrm{ev}}&
X\rtwocell^f_g{\tau}&
Y
\enddiagram
$$
is a morphism of right oplax actions on $Y$.

(2)
The evaluation $\mathrm{ev}:X\otimes [X,X]\to X$ has a canonical
structure of right oplax action (in fact, pseudoaction) and it is
obvious that any 2-functor $Y\otimes -$ preserves right oplax
actions. This shows that $i^Y_{X}$ has 
a canonical opmonoidal structure. The existence of the isomorphism
$(i^W_{Y\otimes X})(i^Y_{X})\cong  i^{W\otimes Y}_{X}$
follows from the fact
that both 1-cells correspond to the right pseudoaction $W\otimes
Y\otimes\mathrm{ev}:W\otimes Y\otimes X\otimes [X,X]\to W\otimes Y\otimes X$.

(3)
The two legs of the rectangle \eqref{fandiiso.eq} correspond, up to
isomorphism, to the 1-cell 
$$
Y\otimes Z\otimes[X,X]\xrightarrow{1\otimes f^*\otimes 1}Y\otimes
X\otimes [X,X]\xrightarrow{1\otimes\mathrm{ev}}Y\otimes
X\xrightarrow{1\otimes f} Y\otimes Z
$$
and therefore there exists an isomorphism as claimed. Moreover, this
isomorphism is monoidal by Proposition \ref{opacteqopmon.p}.

(4)
The 2-cell \eqref{fiimon.eq} corresponds under the closedness
equivalence to 
$$
\diagramcompileto{fiimon1}
Z\otimes X\otimes[X,X]\ar[r]^-{1\otimes\mathrm{ev}}&
Z\otimes X\ar[r]^-{f^*\otimes1}\rrlowertwocell_1{<2>\hole\hole\varepsilon\otimes1}&
Y\otimes X\ar[r]^-{f\otimes 1}&
Z\otimes X.
\enddiagram
$$
This 2-cell is readily shown to be a morphism of right $[X,X]$-actions
on $Z\otimes X$.
\end{proof}

\section{The object of Hopf modules}\label{inthopf.s}


In this section we shall assume that $A$ is a map pseudomonoid in a 
Gray monoid \M\ such that the 2-functor $A\otimes -$ has right
biadjoint $[A,-]$. This is true, for instance, when \M\ is closed; see
Section \ref{opmon.s}. 
Under these assumptions the monad $\theta$ on $\M(A\otimes-,A)$ is
representable by a monad $t:[A,A]\to[A,A]$; that is, there is
an isomorphism
$$
\diagramcompileto{t1}
\M(A\otimes X,A)\drrtwocell<\omit>{'\cong}\ar[rr]^-{\theta_X}\ar[d]_\simeq&&
\M(A\otimes X,A)\ar[d]^\simeq\\
\M(X,[A,A])\ar[rr]_-{\M(X,t)}&&
\M(X,[A,A])
\enddiagram
$$
pseudonatural in $X$.
More explicitly, $t$ is the 1-cell 
\begin{equation}\label{eq:t1}
[A,A]\xrightarrow{i^A_A}[A\otimes A,A\otimes A]\xrightarrow{[p^*,p]}[A,A]
\end{equation}
where $i^A_A$ was defined in Proposition \ref{opmonmorphisms.prop}. 
The multiplication and unit of $t$ are respectively
$$\xymatrixcolsep{1.6cm}
\diagramcompileto{multt}
[A,A]\ar[r]^-{i^A_A}\ar[dr]_{i^{A^2}_A}\ddrlowertwocell_{i_A^A}{}&
[A^2,A^2]\ar@{}[dl]|-(.2){\cong}
\ar[r]^-{[p^*,p]}\ar[d]_{i^A_{A^2}}\drtwocell<\omit>{'\cong}&
[A,A]\ar[d]^{i^A_A}\\
&[A^3,A^3]\ar[r]^-{[1\otimes p^*,1\otimes
p]}\ar[d]|-{[p^*\otimes1,p\otimes 1]}\drtwocell<\omit>{'\cong}&
[A^2,A^2]\ar[d]^{[p^*,p]}\\
&[A^2,A^2]\ar[r]_-{[p^*,p]}&
[A,A]
\enddiagram
$$
$$
\diagramcompileto{unitt}
&[A,A]\ar[dr]|-{[j^*\otimes1,j\otimes 1]}\ar@/^/[drrr]^1&&&\\
[A,A]\ar[ur]^{1=i^I_A}\rrlowertwocell<0>_{i^A_A}{<-3>}&&
[A^2,A^2]\ar@{}[ur]|-(.2){\cong}\ar[rr]_-{[p^*,p]}&&
[A,A]
\enddiagram
$$
where the unlabelled 2-cells are ones defined in Proposition \ref{opmonmorphisms.prop}.4.
Recall that an {\em opmonoidal monad}\/ is a monad in
$\mathbf{Opmon}(\M)$ (see Section \ref{opmon.s}).

\begin{prop}
The monad $t:[A,A]\to[A,A]$ is opmonoidal. 
\end{prop}
\begin{proof}
It is consequence of Proposition \ref{opmonmorphisms.prop} and the
description of the multiplication and unit of $t$ above. 
\end{proof}

Recall that a (bicategorical) Eilenberg-Moore construction for a monad
$s:B\to B$ in 
a bicategory $\mathscr B$ is a birepresentation of the pseudofunctor
$\mathscr B(-,B)^{\mathscr B(-,s)}:\mathscr
B^{\mathrm{op}}\to\mathbf{Cat}$, or equivalently, the unit $u:B^s\to
B$ of that birepresentation.
Opmonoidal monads $s:B\to B$
have the property that if they have an Eilenberg-Moore
construction $u:B^s\to B$ in \M, 
then this construction lifts to $\mathbf{Opmon}(\M)$; in other words, the
forgetful 2-functor $\mathbf{Opmon}(\M)\to\M$ creates Eilenberg-Moore
objects. 
Moreover, $u:B^s\to B$ is strong monoidal
and an arrow $g:C\to B^s$ is opmonoidal (strong monoidal) if and only
if $ug$ is so. The case of $\B=\mathbf{Cat}$ can be found in
\cite{McCrudden:opmonoidalmonads}, while the general case is in
\cite[Lemma 3.2]{Day-Street:quantumcat}.

\begin{defn}
Suppose that the monad $t$ has an Eilenberg-Moore construction
$u:[A,A]^t\to[A,A]$, with $f\dashv u$. So, $[A,A]^t$ has a unique
(up to isomorphism) structure of a pseudomonoid such that $u$ is
strong monoidal.
The Eilenberg-Moore construction 
$u:[A,A]^t\to [A,A]$ is called {\em a
Hopf module construction on $A$.}
\end{defn}

For a justification for the name see \SorCh\ \ref{comodules.sec} below.
The Hopf module construction, of course,
need not  exist in general, and this problem is addressed in the
subsequent sections.

\begin{obs}\label{o:l}
When $A$ has a Hopf module construction the pseudonatural
transformation $\lambda$ in Definition \ref{lambda.def} is
representable by 
\begin{equation}\label{lamrep.eq}
\ell:A\xrightarrow{[j^*,1]}[A, A]\xrightarrow{f}[A,A]^t.
\end{equation}
There exist isomorphisms as depicted below, where $w$ is the 1-cell
corresponding to $1_{A^2}$ under the closedness 
equivalence $\M(A,[A,A^2])\simeq\M(A^2,A^2)$. 
\begin{equation}
  \label{eq:land[1,p]w}
\xymatrixcolsep{1.2cm}
\diagramcompileto{[1,p]i}
A\ar[r]^-{[j^*,1]}\ar[d]_{w}\drtwocell<\omit>{'\cong}&
[A,A]\ar[r]^-f\ar[d]_{i_{A}^A}\drtwocell<\omit>{'\cong}&
[A,A]^t\ar[d]^u\\
[A,A^2]\ar[r]^-{[1\otimes j^*,1]}
\rrlowertwocell_{[1,p]}{\omit}
\ar@{}@<-10pt>[rr]|-\cong
&
[A^2,A^2]\ar[r]^-{[p^*,p]}&
[A,A]
\enddiagram  
\end{equation}
The isomorphism on the right hand side is the isomorphism of
$t$-algebras $uf\cong t$ induced by the universal property of $u$. We
consider $[A,p]w$ as equipped with the unique $t$-algebra structure
such that \eqref{eq:land[1,p]w} is a morphism of
$t$-algebras. Explicitly, this $t$-algebra structure is $t[A,p]w\cong
t[p^*,p]i^A_A[j^*,A]=tt[j^*,A]\to t[j^*,A]\cong[A,p]w$, where the
non-isomorphic arrow induced by the multiplication of $t$.
\end{obs}

\begin{prop}\label{lffstrongmon.prop}
Suppose that $A$ has a Hopf module construction. 
The 1-cell $\ell$ in \eqref{lamrep.eq} is fully faithful and strong
monoidal. Moreover, $\ell$ is an equivalence if and only if the
theorem of Hopf modules holds for $A$ (see Definition
\ref{lambda.def}). 
\end{prop}
\begin{proof}
The first and last assertions follow trivially from Proposition
\ref{lambdaff.prop} and Definition \ref{lambda.def},  so we only have to prove
that $\ell$ is strong monoidal, or equivalently, that $u\ell\cong
t[j^*,A]$ is strong monoidal. This 1-cell is isomorphic to $[A,p]w$ as
in Observation \ref{o:l}.
The 1-cell $[A,p]w:A\to[A,A]$ corresponds up to isomorphism under
$\M(A,[A,A])\simeq\M(A\otimes A,A)$  to
$p:A\otimes A\to A$, which is obviously a right pseudoaction of $A$ on
$A$, and hence $[A,p]w$ is strong monoidal by Proposition
\ref{opacteqopmon.p}. This endows $\ell$ with the structure of a
strong monoidal morphism, by transport of structure. 
\end{proof}


\begin{cor}\label{[1,p]i.cor}
The theorem of Hopf modules holds for $A$ if and only if the 
1-cell 
\begin{equation}\label{[1,p]i.eq}
A\xrightarrow{w}[A,A^2]\xrightarrow{[A,p]}[A,A]
\end{equation}
provides a Hopf module construction for $A$.
\end{cor}
\begin{proof}
The pseudonatural transformation $\lambda$ in Definition
\ref{lambda.def} is an equivalence if and only if the composition
$$
\upsilon\lambda:\M(-,A)\to\M(A\otimes-,A)^\theta\to\M(A\otimes-,A)
$$
is an Eilenberg-Moore construction for the monad $\theta$ in
$[\M^{\mathrm{op}},\mathbf{Cat}]$. But $\upsilon\lambda$ is
represented by the 1-cell \eqref{[1,p]i.eq} and $\theta$ is
represented by $t$, and the result follows. 
\end{proof}

\begin{cor}\label{EMKleisli.cor}
\begin{enumerate}
\item
Suppose that the monad $t$ has an Eilenberg-Moore construction
$f\dashv u:[A,A]^t\to[A,A]$. If the theorem of Hopf modules holds
for $A$ then $f$ is a Kleisli construction for $t$.
\item
Suppose that the monad $t$ has a Kleisli construction
$k:[A,A]\to[A,A]_t$. If the theorem of Hopf modules holds for $A$ 
then $k^*$ is an Eilenberg-Moore construction for $t$.
\end{enumerate}
\end{cor}
\begin{proof}
Let $\mathscr C\subset\M(A\otimes X,A)^{\theta_X}$ be the full image
of the free $\theta_X$-algebra functor $\varphi_X:\M(A\otimes
X,A)\to\M(A\otimes X,A)^{\theta_X}$. When thought of as with codomain
$\mathscr C$, $\varphi_X$ provides a Kleisli construction for
$\theta_X$. The theorem of Hopf modules holds if and only if
$\lambda_X=\varphi_X\M(j^*\otimes X,A)$ is an essentially surjective
on objects, since it is always fully faithful by Proposition
\ref{lambdaff.prop}. Hence, the theorem of Hopf modules holds if and
only if the inclusion of $\mathscr C$ into $\M(A\otimes
X,A)^{\theta_X}$ is an equivalence, which is equivalent to saying that
$\varphi_X$ is a 
(bicategorical) 
Kleisli construction for $\theta$. This proves (1) since $t$ and $f$
represent $\theta$ and $\varphi$ respectively. To show (2), since
$\varphi_X:\M(A\otimes X,A)\to \mathscr C$ is a Kleisli construction
for $\theta_X$, the 1-cell $k^*$ is an Eilenberg-Moore construction
for $t$ if and only if the right adjoint of $\varphi_X$,
$\mathscr C\hookrightarrow\M(A\otimes
X,A)^{\theta_X}\to\M(A\otimes X,A)$, is an Eilenberg-Moore construction
for $\theta_X$ and this happens only if the inclusion $\mathscr
C\hookrightarrow\M(A\otimes X,A)^{\theta_X}$ is an equivalence.
\end{proof}

\section{On the existence of hopf modules}\label{existence.s}

In this section we study the existence of the Hopf
module construction for an arbitrary map pseudomonoid. 
Since this construction is an Eilenberg-Moore
construction for a certain monad, it is natural to embed \M\ into a
2-category where this exists, and the obvious choice is the completion
of \M\ under ($\mathbf{Cat}$-enriched) Eilenberg-Moore objects. This
is a 2-category $\EM(\M)$ with a fully faithful universal 2-functor $E:\M\to
\EM(\M)$. However, in order to speak of the Hopf module construction
for a map pseudomonoid $B$ in $\EM(\M)$ we need $\EM(\M)$ to be
a monoidal 2-category and the pseudofunctor $B\otimes -$ to have right
biadjoint.

We prove that when \M\ is a Gray monoid there exists a model of its
completion under Eilenberg-Moore objects which is also a Gray monoid
and such that the 2-functor $E:\M\to \EM(\M)$ is strict monoidal; this
model is the 2-category explicitly described in
\cite{Lack-Street:ftm2}. In fact, we prove this by extending the
assignment $\M\mapsto\EM(\M)$ to a monoidal functor on the monoidal
category $\mathbf{Gray}$, which turns out to be a
$\mathbf{Gray}$-functor.  In order to show that if $A\otimes -:\M\to
\M$ has right biadjoint then the same is true for $E(A)$ in $\EM(\M)$
we have to move from $\mathbf{Gray}$, where the 1-cells are
2-functors, to $\mathbf{Bicat}$, where 1-cells are pseudofunctors. For
this we extend $\EM$ to a homomorphism of tricategories on
$\mathbf{Bicat}$. 

So far we have only considered bicategorical Eilenberg-Moore
constructions. However, in this section we will use the completion of a
2-category under  $\mathbf{Cat}$-enriched
Eilenberg-Moore objects.
Recall that a $\mathbf{Cat}$-enriched Eilenberg-Moore construction on a monad
$s:Y\to Y$ in a 2-category $\mathscr K$ is a representation of the
2-functor $\mathscr K(-,Y)^{\mathscr K(-,t)}:\mathscr
K^{\mathrm{op}}\to\mathbf{Cat}$. Any 2-categorical
Eilenberg-Moore construction is also a bicategorical one because
2-natural isomorphisms are pseudonatural equivalences. 

From  \cite{Lack-Street:ftm2} we know that $\mathsf{EM}(\mathscr K)$,
the completion under Eilenberg-Moore objects of the 2-category
$\mathscr K$, may be described as the 2-category with objects 
the monads in $\mathscr K$, 1-cells from $(X,r)$ to
$(Y,s)$ monad morphisms, {\em i.e.}, a 1-cells $f:X\to Y$ equipped with a
2-cell $\psi:sf\Rightarrow ft$ satisfying
$$
\diagramcompileto{fmonmor}
X\rruppertwocell^t{^<-2>\mu}\ar[r]^t\ar[d]_f\drtwocell<\omit>{^\psi}&
X\ar[d]^f\ar[r]^t\drtwocell<\omit>{^\psi}&
X\ar[d]^f\\
Y\ar[r]_s&
Y\ar[r]_s&
Y
\enddiagram
=
\xymatrixrowsep{.6cm}
\xymatrixcolsep{.7cm}
\diagramcompileto{fnonmor2}
X\ar[rr]^t\ar[d]_f\drrtwocell<\omit>{^\psi}&&
X\ar[d]^f\\
Y\rrtwocell<\omit>{^<2>\mu}\ar[rr]^s\ar[dr]_s&&
Y\\
&Y\ar[ur]_s
\enddiagram
$$
and
$$
\diagramcompileto{fnonmor3}
X\ar[r]^t\ar[d]_f\drtwocell<\omit>{^\psi}&
X\ar[d]^f\\
Y\ar[r]^s\rlowertwocell_1{^\eta}&
Y
\enddiagram
=
\diagramcompileto{fnonmor4}
X\rtwocell^t_1{^\eta}&X\ar[r]^f&Y
\enddiagram
$$
and  2-cells $(f,\psi)\Rightarrow(g,\chi)$ 
2-cells $\rho:sf\Rightarrow gt$ in $\mathscr K$ such that 
$$
\diagramcompileto{EM2cell1}
X\rruppertwocell^t{^<-2>\mu}\ar[r]^t\ar[d]_f\drtwocell<\omit>{^\rho}&
X\ar[d]^g\ar[r]^t\drtwocell<\omit>{^\chi}&
X\ar[d]^g\\
Y\ar[r]_s&
Y\ar[r]_s&
Y
\enddiagram
=
\xymatrixrowsep{.6cm}
\xymatrixcolsep{.7cm}
\diagramcompileto{EM2cell2}
X\ar[rr]^t\ar[d]_f\drrtwocell<\omit>{^\rho}&&
X\ar[d]^g\\
Y\rrtwocell<\omit>{^<2>\mu}\ar[rr]^s\ar[dr]_s&&
Y\\
&Y\ar[ur]_s
\enddiagram
$$
$$
\diagramcompileto{EM2cell3}
X\rruppertwocell^t{^<-2>\mu}\ar[r]^t\ar[d]_f\drtwocell<\omit>{^\psi}&
X\ar[d]^f\ar[r]^t\drtwocell<\omit>{^\rho}&
X\ar[d]^g\\
Y\ar[r]_s&
Y\ar[r]_s&
Y
\enddiagram
=
\xymatrixrowsep{.6cm}
\xymatrixcolsep{.7cm}
\diagramcompileto{EM2cell4}
X\ar[rr]^t\ar[d]_f\drrtwocell<\omit>{^\rho}&&
X\ar[d]^g\\
Y\rrtwocell<\omit>{^<2>\mu}\ar[rr]^s\ar[dr]_s&&
Y\\
&Y\ar[ur]_s
\enddiagram
$$
This is called the {\em unreduced form}\/ 
of the 2-cells in \cite{Lack-Street:ftm2}.

The completion comes equipped with a fully faithful 2-functor
$E:\mathscr K\to \EM(\mathscr K)$ given on objects by $X\mapsto
(X,1_X)$. This 2-functor has a universal property: for any 2-category
with Eilenberg-Moore objects $\mathscr L$, $E$ induces an isomorphism
of categories $[\EM(\mathscr K),\mathscr L]_{\mathrm{EM}}\to[\mathscr
K,\mathscr L]$, where 
$[\EM(\mathscr K),\mathscr L]_{\mathrm{EM}}\subset[\EM(\mathscr
K),\mathscr L]$ is the full sub 2-category of Eilenberg-Moore
object-preserving 2-functors. Moreover, any object of $\EM(\mathscr
K)$ is the Eilenberg-Moore construction on some monad in the image of
$E$. 

Denote by $\mathbf{Hom}$ the category whose objects are 2-categories
and whose arrows are pseudofunctors. This category is monoidal under
the cartesian product. 
\begin{prop}\label{EMstrongmonoidal.prop}
Completion under Eilenberg-Moore objects defines a strong monoidal
functor $\mathsf{EM}:\mathbf{Hom}\to\mathbf{Hom}$. 
\end{prop}
\begin{proof}
We use the explicit description of the Eilenberg-Moore completion
given in \cite{Lack-Street:ftm2}.
Define $\mathsf{EM}$ on a pseudofunctor $F:\mathscr K\to\mathscr L$ as
sending an object $(X,r)$ to the monad $(FX,Fr)$ in $\mathscr L$, a
1-cell $(f,\psi)$ to $(Ff,F\psi)$ and a 2-cell $\rho$ to $F\rho$. 
The comparison 2-cell $(\EM F(g,\chi))(\EM F(f,\psi))\to\EM
F((g,\chi)(f,\psi))$
is defined to be
$$
(Fr)(Fg)(Ff)\xrightarrow{\cong}F(rgf)\xrightarrow{F((g\psi)\cdot(\chi
f))} F(gft)\xrightarrow{\cong}F(gf)(Ft) 
$$
or what is the same thing
\begin{multline}\label{EMconstraint.eq}
(Fr)(Fg)(Ff)\xrightarrow{\cong}F(rg)(Ff)\xrightarrow{(F\chi)(Ff)}F(sg)(Ff)
\xrightarrow{\cong}(Fg)F(sf)\longrightarrow\\
\xrightarrow{(Fg)(F\psi)}
(Fg)F(ft)\xrightarrow{\cong} F(gf)(Ft) 
\end{multline}
where the unlabelled isomorphisms are (the unique possible) 
compositions of the structural
constraints of the pseudofunctor $F$.
The axioms of a 2-cell in $\EM(\mathscr L)$ follow from the fact that
$(Fg)(Ff)$ and $F(gf)$ are monad morphisms. Similarly, the identity
constraint of $1_{\EM F(X)}\to (\EM F)(1_X) $ is defined as 
$$
\big((Ft)1_{FX}\xrightarrow{(Ft)F_0}
(Ft)(F1_x)\xrightarrow{\cong}(F1_x)(Ft)\big)=
\big(1_{FX}(Ft)\xrightarrow{F_0(Ft)}(F1_X)(Ft)\big)
$$
where $F_0$ is the identity constraint of $F$. 

It is clear that this defines a functor $\mathsf{EM}$. It is also clear
that it is strong monoidal, with constraints the evident isomorphisms 
$\mathsf{EM}(\mathscr K)\times
\mathsf{EM}(\mathscr L)\cong \mathsf{EM}(\mathscr K\times\mathscr L)$ and
$E_1:1\cong \mathsf{EM}(1)$. 
\end{proof}

\begin{obs}\label{o:EMpresbieq}
  If $F:\K\to \mathscr{L}$ is a biequivalence between 2-categories, then $\EM F$
  is a biequivalence too. This is straightforward from the definition
  of $\EM$ on pseudofunctors 
  in the proof of Proposition \ref{EMstrongmonoidal.prop} above. 
\end{obs}

Recall from Section \ref{s:prelimin} the notion of cubical functor.

\begin{cor}\label{EMfcubical.cor}
The pseudofunctor below is a cubical functor whenever $F:\mathscr
K\times\mathscr L\to\mathscr J$ is one. 
$$
\EM(\mathscr K)\times \EM(\mathscr
L)\xrightarrow{\cong}\EM(\K\times\mathscr{L}) \xrightarrow{\EM F}\EM(\mathscr J)
$$
\end{cor} 
\begin{proof}
Consider 1-cells in $\EM(\mathscr K)\times\EM(\mathscr L)$
$$
((X',t'),(X,t))\xrightarrow{((f',\psi'),(f,\psi))}((Y',s'),(Y,s))
\xrightarrow{((g',\chi'),(g,\chi))} ((Z',r'),(Z,r)).
$$
If $(X,t)=(Y,s)$ and $(f,\psi)$ is the identity 1-cell of $(X,t)$,
that is $(f,\psi)=(1_X,1_t)$, then the constraint defined in
\eqref{EMconstraint.eq} above is
\begin{multline*}
F(r',r)F(g',g)F(f',1)=F(r',r)F(g'f',g)\xrightarrow{\cong}
F(r'g'f',rg)\\
\xrightarrow{F((g'\psi')\cdot(\chi' f'),\chi)}
F(g'f't',gt)\xrightarrow{\cong}F(g'f',g)F(t',t')
\end{multline*}
which is exactly the identity 2-cell of the 1-cell $\EM
F((g',\chi')(f',\psi'),(g\chi))$ in the 2-category $\EM(\mathscr J)$. The rest of the
proof is similar.
\end{proof}

Recall from Section \ref{s:prelimin} the Gray tensor product of
2-categories. If $\K,\mathscr L$ are 2-categories, its Gray tensor
product $\K\square\mathscr L$ is a 2-category classifying cubical
functors out of $\K\times\mathscr L$. 

\begin{cor}\label{c:EMmonoidal}
Completion under Eilenberg-Moore objects induces a monoidal functor
$\EM$ from  $\mathbf{Gray}$ to itself. Furthermore, the 2-functors
$E_{\K}:\K\to \EM(\K)$ are the components of a monoidal natural
transformation. 
\end{cor}
\begin{proof}
Define the structural arrow $\EM(\mathscr K)\square\EM(\mathscr
L)\to\EM(\mathscr K\square\mathscr L)$ as corresponding to 
$\EM(\K)\times\EM(\mathscr L)\cong\EM(\mathscr K\times\mathscr
L)\to\EM(\mathscr K\square\mathscr L)$, 
which is a cubical functor by Corollary \ref{EMfcubical.cor},
and the arrow $1\to\EM(1)$ as the universal
$E_1$. Here the symbol $\square$ denotes the Gray tensor product.
The axioms of lax monoidal functor follow from the fact that $\EM$ is
monoidal with respect to the cartesian product. 

The naturality of the arrows $E_{\K}$ follows form the universal
property of the completions under Eilenberg-Moore objects. We only
have to prove that the resulting natural transformation is
monoidal. Consider the diagram
$$
\xymatrixcolsep{.7cm}
{\diagramcompileto{EMmonoidal}
\EM(\K)\square\EM(\mathscr L) \ar@/^/[drrr]&&&\\
&\EM(\K)\times\EM(\mathscr L)\ar[ul]\ar[r]^-\cong&
\EM(\K\times\mathscr L)\ar[r]&
\EM(\K\square\mathscr L)\\
&\K\times\mathscr L\ar[u]^{E_{\K}\times E_{\mathscr
    L}}\ar[ur]_{E_{\K\times\mathscr L}}\ar[dl]&&\\
\K\square\mathscr L\ar[uuu]^{E_{\K}\square E_{\mathscr
    L}}\ar@/_/[uurrr]_{E_{\K\square\mathscr L}}
\enddiagram}
$$
One of the two axioms we have to check is the commutativity of the
exterior diagram. This commutativity can be proven by observing that 
each one of the four internal
diagrams commute and then applying the universal property of
$\K\times\mathscr L\to\K\square\mathscr L$. The other axiom, involving
$E_1:1\to\EM(1)$ is trivial, since $E_1$ itself is the unit
constraint. 
\end{proof}

\begin{cor}\label{c:EMMgraymonoid}
$\mathsf{EM}(\M)$ is a Gray monoid whenever \M\ is a Gray monoid. 
Moreover, the
  2-functor $E_{\M}:\M\to \mathsf{EM}(\M)$ is strict monoidal, so that $\M$
  can be identified with a full monoidal sub 2-category of
  $\mathsf{EM}(\M)$. 
\end{cor}
\begin{proof}
We know that $\mathsf{EM}$ is a monoidal functor, and as such
it preserves monoids. Moreover, $E_{\M}$ is strict monoidal, that is,
a morphism of monoids in $\mathbf{Gray}$, since $E$ is a monoidal
natural transformation (see Corollary \ref{c:EMmonoidal}). 
\end{proof}

The tensor product in $\mathsf{EM}(\M)$ is induced by the one of \M; for
instance, the tensor product of $(X,r)$ with $(Y,s)$, denoted by
$(X,r)\circledcirc(Y,s)$, 
is $(X\otimes
Y,r\otimes s)$.

In order to show that $\EM$ is in fact a $\mathbf{Gray}$-functor we
state the following easy result.
\begin{lem}
Let $\mathscr V$ be a symmetric monoidal closed category and
$F:\mathscr V\to\mathscr V$ be a lax monoidal functor. Then, any 
monoidal natural transformation $\eta:1_{\mathscr V}\Rightarrow F$
induces on $F$ a structure of a $\mathscr V$-functor.
\end{lem}
\begin{proof}
Define $F$ on enriched homs as
$$
F_{X,Y}:[X,Y]\xrightarrow{\eta_{[X,Y]}}F([X,Y])\xrightarrow{\vartheta_{X,Y}}[FX,FY] 
$$
where $\vartheta_{X,Y}$ is the arrow corresponding to 
$
F[X,Y]\otimes FX\longrightarrow F([X,Y]\otimes X)\xrightarrow{F\mathrm{ev}}FY.
$
\end{proof}

\begin{cor}
$\EM:\mathbf{Gray}\to\mathbf{Gray}$ has a canonical structure of
  $\mathbf{Gray}$-functor. 
\end{cor}
\begin{proof}
Let  $\mathscr V$ in the lemma above  be $\mathbf{Gray}$ and  $\eta$ 
be  the transformation defined by the inclusions $E_{\mathscr
  K}:\mathscr K\to\EM(\mathscr K)$, which is easily shown to be a
monoidal transformation. Now apply the lemma.
\end{proof}

Let $\mathbf{Ps}(\K,\mathscr L)$ denote the 2-category of pseudofunctors
from $\K$ to $\mathscr L$, pseudonatural transformations between
them and modifications between these.

\begin{obs}\label{EMonhoms.obs}
In the case of $\EM$, the transformation $\vartheta_{\mathscr K,\mathscr
L}$ 
is defined by
the commutativity of the following diagram
$$
\xymatrixcolsep{1.5cm}
\diagramcompileto{EM3}
{\mathbf{Ps}}(\EM(\mathscr K),\EM(\mathscr L))\otimes\EM(\mathscr
K)\ar[r]^-{\mathrm{ev}}&
\EM(\mathscr L)\\
\EM{\mathbf{Ps}}(\mathscr K,\mathscr L)\otimes\EM(\mathscr
K)\ar[r]\ar[u]^{\vartheta_{\mathscr K,\mathscr L}\otimes1}&
\EM({\mathbf{Ps}}(\mathscr K,\mathscr L)\otimes\mathscr
L)\ar[u]_{\EM\mathrm{ev}}\\
\enddiagram
$$
that is, 
\begin{align*}
(\vartheta_{\mathscr K,\mathscr L}(F,\tau))(X,t)&
=(\EM\mathrm{ev})((F,\tau),(X,t))\\
&=(\mathrm{ev}(F,X),\mathrm{ev}(\tau,t))\\
&=(FX,(Ft)\tau_X),
\end{align*}
and then $\EM$ is defined on homs by the 2-functor 
$$
\vartheta_{\mathscr K,\mathscr
L}E_{\mathscr K,\mathscr L}
:\mathbf{Ps}(\mathscr K,\mathscr L)\to\mathbf{Ps}(\EM(\mathscr
K),\EM(\mathscr L))
$$ 
whose value on a 2-functor $F$ is the 2-functor
sending a monad $(X,t)$ to $(FX,Ft)$.  Then we see that our
$\mathbf{Gray}$-functor has as underlying ordinary functor just the
restriction to $\mathbf{Gray}$ of the functor in Proposition
\ref{EMstrongmonoidal.prop}. 
\end{obs}

Denote by 
$\mathbf{Bicat}$ the tricategory of bicategories,
pseudofunctors, pseudonatural transformations and modifications as
defined in \cite[5.6]{tricategories}. (There is another canonical choice for
a tricategory structure on $\mathbf{Bicat}$, as explained in the
that paper.) 
We shall describe an extension of 
the $\mathbf{Gray}$-functor $\EM$ to a homomorphism of tricategories
$\widetilde{\EM}:\mathbf{Bicat}\to\mathbf{Bicat}$. 
In order to do this we will  use the construction of a homomorphism of 
tricategories
$\mathbf{Bicat}\to\mathbf{Gray}$ given in \cite{tricategories}, of
which we recall some aspects. 
For each bicategory $\mathscr B $ there is a 2-category
$\mathbf{st}\mathscr B$ and a pseudofunctor $\xi_{\mathscr B}:\mathscr
B\to\mathbf{st}\mathscr B$ inducing for each 2-category $\mathscr K$
an {\em isomorphism}\/ of 2-categories $\mathbf{Bicat}(\mathscr
B,\mathscr K)\cong\mathbf{Ps}(\mathbf{st}\mathscr B,\mathscr
K)$. Moreover, $\xi_{\mathscr B}$ is a biequivalence of bicategories. 
As usual, we get a pseudofunctor
$$
\mathbf{st}_{\mathscr A,\mathscr B}:\mathbf{Bicat}(\mathscr A,\mathscr
B)\to\mathbf{Ps}(\mathbf{st}\mathscr A,\mathbf{st}\mathscr B)
$$
which turns out to be an biequivalence. Finally, the object part of
the homomorphism of tricategories $\mathbf{Bicat}\to \mathbf{Gray}$ is given
by $\mathscr B\mapsto\mathbf{st}\mathscr B$ while on hom-bicategories
it is given by the biequivalence $\mathbf{st}_{\mathscr A,\mathscr
B}$. 

Define a homomorphism of tricategories $\widetilde{\EM}$ by
$$
\diagramcompileto{EM4}
{\mathbf{Bicat}}\ar[r]^\sim\ar@{..>}[d]_{\widetilde{\EM}}&
{\mathbf{Gray}}\ar[d]^{\EM}\\
{\mathbf{Bicat}}&
{\hole\mathbf{Gray}\hole}\ar@{_(->}[l]
\enddiagram
$$
It is given on objects by $\mathscr B\mapsto\EM(\mathbf{st}\mathscr B)$ and
on homs by 
$$
\mathbf{Bicat}(\mathscr A,\mathscr
B)\xrightarrow{\mathbf{st}}\mathbf{Ps}(\mathbf{st}\mathscr
A,\mathbf{st}\mathscr
B)\xrightarrow{E}\EM\mathbf{Ps}(\mathbf{st}\mathscr
A,\mathbf{st}\mathscr
B)\xrightarrow{\vartheta}\mathbf{Ps}(\EM\mathbf{st}\mathscr
A,\EM\mathbf{st}\mathscr B),
$$
which by the Observation \ref{EMonhoms.obs} sends a pseudofunctor 
$F:\mathscr A\to\mathscr B$ to the 2-functor $\EM(\mathbf{st} F)$
defined in Proposition \ref{EMstrongmonoidal.prop}.

\begin{prop}
Every biadjunction between pseudofunctors $F\dashv_{\mathrm
b} G:\mathscr L\to\mathscr K$, where  $\mathscr K$ and $\mathscr L$
are 2-categories, induces a biadjunction $\EM
F\dashv_{\mathrm b}\EM G$.
\end{prop}
\begin{proof}
Since $\widetilde{\EM}$ is a homomorphism of tricategories on
$\mathbf{Bicat}$, 
$\widetilde{\EM} F=\mathsf{EM}(\mathbf{st}F)$ is left
biadjoint to $\widetilde{\EM} G=\mathsf{EM}(\mathbf{st}G)$.
The 2-functor $\mathbf{st}F$ is defined as the unique 2-functor such
that $(\mathbf{st}F)\xi_{\mathscr K}=\xi_{\mathscr L}F$, and similarly
for $G$. 
It follows, by
functoriality of $\EM$ with respect to pseudofunctors (Proposition
\ref{EMstrongmonoidal.prop}), that 
$$
\EM(\mathbf{st}F)\EM\xi_{\mathscr K}=\EM\xi_{\mathscr L}\EM
F\quad\text{and}\quad \EM(\mathbf{st}G)\EM\xi_{\mathscr L}=\EM\xi_{\K}\EM G.
$$
Since each component of $\xi$ is a
biequivalence and these are preserved by $\EM$ (see Observation
\ref{o:EMpresbieq}), we have 
$$
\EM F\simeq (\EM\xi_{\mathscr
  L})^*\EM(\mathbf{st}F)\EM\xi_{\mathscr K} \dashv_{\mathrm b} 
(\EM\xi_{\mathscr  K})^*\EM(\mathbf{st}G)\EM\xi_{\mathscr L}\simeq \EM
G
$$
\end{proof}

\begin{cor}\label{c:MclosedEMclosed}
If $X$ is an object in a Gray monoid \M\ such that $X\otimes -$ has
right biadjoint $[X,-]$, then $(EX\circledcirc-):\EM(\M)\to\EM(\M)$
has right biadjoint 
$\langle EX,-\rangle$ given by  $\langle EX,(Y,s)\rangle=([X,Y],[X,s])$.
\end{cor}
\begin{proof}
The 2-functor $(EX\circledcirc-)$ is just $\EM(X\otimes -)$, and
then by the proposition above it has right biadjoint
$\EM([X,-])$. This is given by the stated formula as a consequence of the
description of the effect of $\EM$ on pseudofunctors in the proof of
Proposition \ref{EMstrongmonoidal.prop}.  
\end{proof}

\begin{thm}
For any closed Gray monoid \M\ there exists another Gray monoid $\mathscr N$
and a fully faithful strict monoidal 2-functor $\M\to\mathscr N$ such
that any map pseudomonoid in \M\ has a Hopf module construction in
$\mathscr N$. Moreover, $\mathscr N$ can be taken to be $\EM(\M)$. 
\end{thm}
\begin{proof}
  The proof is only a matter of putting Corollaries
  \ref{c:EMMgraymonoid} and \ref{c:MclosedEMclosed} together with the
  definition of the object of Hopf modules.
\end{proof}


\begin{prop}\label{HopfforA.prop}
Let $A$ be a map pseudomonoid in a Gray monoid \M\ such that $A\otimes
-$ has right biadjoint. Suppose that the theorem of Hopf modules holds for 
$E(A)\in\mathrm{ob}\EM(\M)$; then it also holds for $A$. Moreover, in
this case $A$ has a Hopf module construction provided by 
\begin{equation}\label{[1,p]i.eq2}
A\xrightarrow{w}[A,A\otimes A]\xrightarrow{[A,p]}[A,A]
\end{equation}
as in Corollary \ref{[1,p]i.cor}.
\end{prop}
\begin{proof}
Consider the image of the monad $t$ under the 2-functor
$E:\M\to\EM(\M)$. Denote by $\hat\theta$ the monad
$\EM(\M)(-,Et)$ on $\EM(\M)(-,E[A,A])$ and
$\hat\varphi\dashv\hat\upsilon$ the adjunction arising from its
Eilenberg-Moore construction in
$\mathbf{Hom}(\EM(\M)^{\mathrm{op}},\mathbf{Cat})$. Observe that by
the fully faithfulness of $E$, the monad $\hat\theta_{EXMP}$ can be
identified with the monad $\theta_X$ of Definition \ref{d:theta}, and 
the adjunction $\hat\varphi_{EXMP}\dashv\hat \upsilon_{EXMP}$ with the
adjunction $\varphi_X\dashv\upsilon_X$ corresponding to $\theta$. 

If the theorem of Hopf modules 
holds for $E(A)$ then in particular for each object $X$ of \M\ the
functor 
\begin{multline}\label{hatlambdaEX.eq}
\EM(\M)(E(X),E(A))\xrightarrow{\EM(\M)(1,E([j^*,A]))}
\EM(\M)(E(X),E[A,A])\longrightarrow\\
\xrightarrow{\hat\varphi_{E(X)}}\EM(\M)(E(X),E[A,A])^{\hat\theta_{E(X)}}
\end{multline}
is an equivalence (Definition \ref{lambda.def}). But by the fully
faithfulness of the 2-functor $E$ 
this is, up to composing with suitable isomorphisms, just the functor
$\lambda_X$ in Definition \ref{lambda.def} and then the theorem of
Hopf modules holds for $A$.


The last assertion follows directly from Corollary \ref{[1,p]i.cor}.
\end{proof}

\section{Left autonomous pseudomonoids}\label{maintheorem.s}

In this section we specialise to a special kind of pseudomonoid,
central to our work, namely the autonomous pseudomonoids. We begin by
recalling the necessary background. 

A {\em bidual pair}\/ in a Gray monoid \M\ is a pseudoadjunction (see
for example \cite{Lack:Acoherentapproach}) in the one-object Gray-category
\M. Explicitly, it consists of a pair of 1-cells $\e:X\otimes Y\to I$
and $\n:I\to Y\otimes X$ together with 
invertible 2-cells
\begin{equation}
\begin{split}
{\diagramcompileto{biduals1}
Y\rruppertwocell<0>^1{<3>\eta}\ar[dr]_{\n\otimes 1}&&Y\\
&Y\otimes X\otimes Y\ar[ur]_{1\otimes\e}&
\enddiagram}
{\diagramcompileto{biduals2}
&X\otimes Y\otimes X\ar[dr]^{\e\otimes 1}&\\
X\ar[ur]^{1\otimes \n}\rrlowertwocell<0>_1{<-3>\varepsilon}&&X
\enddiagram}
\end{split}
\end{equation}
such that the
following composites are identities.
\begin{equation}\label{cond1.lax.e}
\begin{split}
\xymatrixrowsep{.9cm}
\xymatrixcolsep{.9cm}
{\diagram
&&f\otimes u\ar[rd]^{\e}&\\
f\otimes u \drrlowertwocell_1{<-.1>\hole\varepsilon\otimes1}
\urruppertwocell^1{<.1>\hole 1\otimes\eta}\ar[r]^-{1\otimes\n\otimes
1}&
f\otimes u\otimes f\otimes u
\ar[ur]|-{1\otimes1\otimes\e}\ar[dr]|-{\e\otimes
1\otimes 1}\rruppertwocell<\omit>{\hole c^{-1}_{\e,\e}}&&1\\
&&f\otimes u\ar[ur]_\e&
\enddiagram}
\end{split}
\end{equation}
\begin{equation}
\begin{split}
\xymatrixrowsep{.9cm}
\xymatrixcolsep{.9cm}
{\diagram
&u\otimes f\ar[dr]|-{\n\otimes1\otimes1}
\drruppertwocell^1{<.1>\hole\eta\otimes1}&&\\
1\ar[ur]^\n\ar[dr]_\n\rruppertwocell<\omit>{\hole c^{-1}_{\n,\n}}&&u\otimes
f\otimes u\otimes f\ar[r]^-{1\otimes \e\otimes 1}&u\otimes f\\
&u\otimes f\ar[ur]|-{1\otimes 1\otimes
\n}\urrlowertwocell_1{<-.1>\hole1\otimes \varepsilon}
\enddiagram}
\end{split}
\end{equation}
The object $X$ is called a {\em right bidual}\/ of $Y$, denoted by
$Y^\circ$, and $Y$ is called a {\em  left bidual}\/ of $X$, denoted by
$X^\vee$. 

A Gray monoid in which every object has a right (left) bidual
is called right (left) autonomous.

If $X$ has a right bidual $X^\circ$, then the 2-functor $X\otimes-$ has a right
biadjoint $X^\circ\otimes-$, and $-\otimes X$ has a left biadjoint
$-\otimes X^\circ$, and dually for left biadjoints. In particular,
any right (left) autonomous Gray monoid is right (left) closed with
internal hom $[X,Y]=X^\circ\otimes Y$ ($[X,Y]=Y\otimes X^\vee$). If
both $X$ and $Y$ have bidual and $f:X\to Y$, the bidual of $f$ is the
1-cell $f^\circ=(X^\circ\otimes\e)(X^\circ \otimes f\otimes
Y^\circ)(\n\otimes Y^\circ):Y^\circ\to X^\circ$. Similarly with
2-cells. If $\mathscr N$ is the full sub-2-category of \M\ whose
objects are the objects with right bidual, we have a monoidal
pseudofunctor $(-)^\circ:(\mathscr N^{\mathrm{op}})^{\mathrm{rev}}\to
\M$, where the superscript rev indicates the reverse monoidal
structure. The structural constraints are given by the canonical
equivalences $I\simeq I^\circ$ and $Y^\circ\otimes
X^\circ\simeq(X\otimes Y)^\circ$. 

Recall from \cite{dualizations} that a {\em left dualization}\/ for a
pseudomonoid $(A,j,p)$ in \M\ is a 1-cell $d:A^\circ\to A$ equipped
with two 2-cells $\alpha:p(d\otimes A)\n\Rightarrow j$ and
$\beta:j\e\Rightarrow p(A\otimes d)$ satisfying two axioms. 
Left dualization structures on
$d:A^\circ\to A$ are in bijection with adjunctions $p(d\otimes
A)\dashv (A^\circ\otimes p)(\n\otimes A)$ and with adjunctions
\begin{equation}\label{p*.eq}
p\dashv (p\otimes A)(A\otimes d\otimes A)(A\otimes\n).
\end{equation} 
For example, given $\alpha $ and
$\beta$ the counit of the corresponding adjunction \eqref{p*.eq} is
$$
\diagramcompileto{counit.alpha}
A\otimes\Ac \otimes A\ar[r]^-{1\otimes d\otimes
  1}\rtwocell<\omit>{<3>\hole\hole1\otimes\alpha}&
A^3\ar[r]^-{p\otimes 1}\ar[d]|-{1\otimes p}\drtwocell<\omit>{'\cong}&
A^2\ar[d]^p\\
A\ar[u]^{1\otimes\n}\ar[r]^-{1\otimes j}
\rrlowertwocell_1{\cong}&
A^2\ar[r]^-p&
A
\enddiagram
$$
(To be precise, in \cite{dualizations} the authors define left
dualization in a right  autonomous Gray monoid, {\em i.e.}, a Gray
monoid where any object has a right bidual, but the only really
necessary condition is that the pseudomonoid itself have a right bidual).

A pseudomonoid equipped with a left dualization is called {\em left
  autonomous.} 
If a left dualization exists, then it is isomorphic to
  $(A\otimes\e)(p^*\otimes\Ac)(j\otimes \Ac)$ \cite[Proposition
  1.2]{dualizations}. Examples of this structure are the left autonomous
  (pro)monoidal $\mathscr V$-categories (for a good monoidal category
  $\mathscr V$) and (co)quasi-Hopf algebras. See \cite{dualizations}
  or \SorCh\ \ref{VMod.sec} and \ref{comodules.sec}.

Given a left autonomous pseudomonoid $A$ define the following 2-cell.
\begin{equation}\label{gamma.eq}
\gamma:=
\diagramcompileto{gamma}
&A^3\ar[r]^-{p\otimes1}\ar[dr]|-{1\otimes p}
\drrtwocell<\omit>{'\cong\ \phi}&
A^2\ar[dr]^p\ar@{=}[rr]\rrtwocell<\omit>{<3>\eta}&&
A^2\\
A^2\ar[ur]^{1\otimes
  p^*}\ar@{=}[rr]\rrtwocell<\omit>{<-3>1\otimes\varepsilon\hole\hole\hole\hole\hole\hole} &&
A^2\ar[r]_-p&
A\ar[ur]_{p^*}
\enddiagram
\end{equation}
In the lemma below we show that this 2-cell $\gamma$ is
invertible, and in fact this property will turn out to be equivalent to
the existence of a left dualization. 

\begin{lem}\label{lemmalapdm.l}
For a left autonomous pseudomonoid $A$ the following equality holds.
\begin{equation}\label{gamma1.eq}
\gamma
=
\xymatrixcolsep{1.5cm}
\diagramcompileto{gamma=V2}
&&&\\
A^2\ar `u[r]`[rrr]^{1\otimes p^*}[rrr]
\ar[d]_{p}\ar[r]^-{A^2\otimes \n}\drtwocell<\omit>{'\cong}&
A^2\otimes \Ac\otimes A\ar[d]|-{p\otimes 1\otimes 1}\ar[r]^-{A^2\otimes
d\otimes 1}\drtwocell<\omit>{'\cong}&
A^4\ar[r]^-{1\otimes p\otimes 1}\ar[d]|-{p\otimes
A^2}\drtwocell<\omit>{'\cong}&
A^3\ar[d]^{p\otimes 1}\\
A\ar@{-<} `d  [r] `[rrr]_{p^*}  [rrr]\ar[r]_-{1\otimes \n}&
A\otimes \Ac\otimes A\ar[r]_-{1\otimes d\otimes 1}&
A^3\ar[r]_-{p\otimes 1}&
A^2\\
&&&
\enddiagram
\end{equation}
In particular, $\gamma $ is invertible.
\end{lem}
\begin{proof}
The 2-cell on the right hand of \eqref{gamma1.eq} pasted with the
counit of the adjunction \eqref{p*.eq} gives the following 2-cell
$$
\xymatrixcolsep{1.5cm}
\diagramcompileto{gammaV2proof1}
A^2\ar[r]^-{A^2\otimes \n}\ar[d]_{p}\drtwocell<\omit>{'\cong}&
A^2\otimes \Ac\otimes A\ar[d]|-{p\otimes 1\otimes 1}\ar[r]^-{A^2\otimes
d\otimes1}\drtwocell<\omit>{'\cong}&
A^4\ar[r]^-{1\otimes p\otimes
1}\ar[d]|-{p\otimes1\otimes1}\drtwocell<\omit>{'\cong}& 
A^3\ar[d]^{p\otimes 1}\\
A \ar@(d,d)[drrr]_1 \ar@{}@<-20pt>[drrr]|-{\cong}\ar[r]^-{1\otimes \n}\ar@/_/[drr]_{1\otimes j}&
A\otimes \Ac\otimes A\ar[r]_-{1\otimes d\otimes 1}&
A^3\ar[r]_-{p\otimes1}\ar[d]|-{1\otimes p}\drtwocell<\omit>{'\cong}&
A^2\ar[d]^{p}\\
&&A^2\ar[r]_-p\ultwocell<\omit>{^\hole\hole 1\otimes \alpha}&
A
\enddiagram
$$
which itself is equal to 
$$
\xymatrixcolsep{1.5cm}
\diagramcompileto{gammaV2proof2}
A^2\drrtwocell<\omit>{<-1>\hole\hole A^2\otimes \alpha}\ar[dd]_p\ar[r]^-{A^2\otimes \n}\ar@/_/[drr]_{A^2\otimes j}&
A^2\otimes \Ac\otimes A\ar[r]^-{A^2\otimes d\otimes 1}&
A^4\ar[d]_{A^2\otimes p}\ar[dr]|-{p\otimes A^2}\ar[r]^-{1\otimes
p\otimes1}&
A^3\ar[dr]^{p\otimes1}\dtwocell<\omit>{'\cong}&\\
&&A^3\ar[d]_{p\otimes1}\rtwocell<\omit>{'\cong}&
A^3\ar[r]_-{p\otimes 1}\ar[dl]|-{1\otimes p}&
A^2\ar[dl]^p\\
A\urrtwocell<\omit>{'\cong}
\ar@{-<} `d [r]  `[rrr]_1 [rrr]\ar@{}@<-10pt>[rrr]|-{\cong} \ar[rr]_-{1\otimes j}&&
A^2\ar[r]_-p\urrtwocell<\omit>{'\cong}&
A\\
&&&&
\enddiagram
$$
$$
=
\xymatrixcolsep{1.5cm}
\diagramcompileto{gammaV2proof3}
A^2\drrtwocell<\omit>{<-1>\hole\hole A^2\otimes \alpha}
\ar[dd]_{p\otimes 1}\ar[r]^-{A^2\otimes \n}\ar@/_/[drr]_{A^2\otimes j}&
A^2\otimes\Ac\otimes A\ar[r]^-{A^\otimes
d\otimes1}&
A^4\ar[d]|-{A^2\otimes p}\ar[r]^-{1\otimes p\otimes
1}\drtwocell<\omit>{'\cong}&
A^3\ar[d]|-{1\otimes p}\ar[dr]^{p\otimes 1}&\\
&&A^3\ar[d]_{p\otimes 1}\drtwocell<\omit>{'\cong}\ar[r]_-{1\otimes p}&
A^2\ar[d]|-{p}\rtwocell<\omit>{'\cong}&
A^2\ar[dl]^{p}\\
A\urrtwocell<\omit>{'\cong}\ar[rr]_-{1\otimes j} \ar@{-<} `d [r]  `[rrr]_1
[rrr]\ar@{}@<-10pt>[rrr]|-{\cong} &&
A^2\ar[r]_-p&
A
\enddiagram
$$
$$
=
{\xymatrixcolsep{1.5cm}
\diagramcompileto{gammaV2proof4}
A^2 \ar@(d,d)[drrr]_1\ar@{}@<-20pt>[drrr]|-\cong\ar[r]^-{A^2\otimes \n}\ar@/_/[drr]_{A^2\otimes
j}\drrtwocell<\omit>{<-1>\hole\hole A^2\otimes \alpha} &
A^2\otimes \Ac\otimes A\ar[r]^-{A^2\otimes d\otimes 1}&
A^4\ar[r]^-{1\otimes p\otimes
1}\drtwocell<\omit>{'\cong}\ar[d]|-{A^2\otimes p}&
A^3\ar[d]|-{1\otimes p}\ar[dr]^{p\otimes 1}&\\
&&A^3\ar[r]_-{1\otimes p}&
A^2\ar[d]_{p}\rtwocell<\omit>{'\cong}&
A^2\ar[dl]^p\\
&&&A&
\enddiagram}
$$
$$
=
{\xymatrixcolsep{1.5cm}
\diagramcompileto{gammaV2proof5}
A^2\ar[r]^-{1\otimes p^*}\drtwocell<0>_1{<-2>\hole\hole1\otimes\varepsilon}&
A^3\ar[d]|-{1\otimes p}\ar[dr]^{p\otimes1}&\\
&A^2\ar[d]_p\rtwocell<\omit>{'\cong}&
A^2\ar[dl]^p\\
&A
\enddiagram}
$$
The result follows.
\end{proof}

Define the 2-cell $\omega$ as
\begin{equation}\label{omega.eq}
\diagramcompileto{omega}
A\ar[r]^-{1\otimes j}\drlowertwocell_1{\omit}\drtwocell<\omit>{'\cong}&
A^2\ar[d]^p\ar[r]^-{1\otimes p^*}\rtwocell<\omit>{<4>\gamma}&
A^3\ar[d]^{p\otimes1}\\
&A\ar[r]_{p^*}&
A^2
\enddiagram
\end{equation}

Now we state the basic result of this work. 

\begin{thm}\label{maintheorem.th}
Let $(A, j, p)$ be a map pseudomonoid in a  Gray
monoid \M\ and suppose that $A$ has a right bidual.
Then, the following assertions are equivalent.
\begin{enumerate}
\item
$A$ is left autonomous.
\item\label{maintheoremgamma.it}
The 2-cell $\gamma$ in \eqref{gamma.eq} is invertible.
\item 
The 2-cell $\omega$ in \eqref{omega.eq} is invertible.
\item 
The theorem of Hopf
modules holds for $A$.
\item
The functor
$\lambda_{\Ac}:\M(\Ac,A)\xrightarrow{\M(j^*\otimes1,1)}
\M(A\otimes\Ac,A)\xrightarrow{\varphi_{\Ac}}
\M(A\otimes\Ac,A)^{\theta_{\Ac}}$ is an equivalence.
\end{enumerate}
\end{thm}
\begin{proof}
(1) implies (2) by Lemma \ref{lemmalapdm.l}.3, and (3) follows
trivially from (2)  as (5) does from (4).
By Observation \ref{obs1.obs}, 
to prove that (3) implies (4) it is enough to show that for each
object $X$ the natural transformation
$\upsilon_X\varepsilon_X\varphi_X$ is an isomorphism.  For
$g\in\M(A\otimes X,A)^{\theta_X}$, the component
$\upsilon_X(\varepsilon_X)_{g}$ is the pasting
$$
\diagramcompileto{thm1}
A\otimes X\ar[r]^-{p^*\otimes1}\ar@/_3pc/[rrrrr]_g&
A^2\otimes X\ar[r]^-{1\otimes
  j^*\otimes1}\rrlowertwocell<-7>_1{<2>}&
A\otimes X\ar[r]^-{1\otimes j\otimes 1}\rtwocell<\omit>{<5>\nu}&
A^2\otimes X\ar[r]^-{1\otimes g}&
A^2\ar[r]^-p&
A
\enddiagram
$$
where $\nu$ is the action of $\theta_X$ on $g$ and the unlabelled
arrow is induced by the counit of $j\dashv j^*$.
This 2-cell pasted with $1_{A\otimes X}\cong (A\otimes j^*\otimes
X)(p^*\otimes X)$ gives, by the equality \eqref{eq:jyp},
\begin{equation}\label{thm1.eq}
\diagramcompileto{thm2}
A\otimes X\ar[r]^-{1\otimes j\otimes 1}\drrlowertwocell_1{\cong}&
A^2\otimes
X\ar@{=}[rr]\rrtwocell<\omit>{<2.5>\hole\hole\eta\otimes1}\ar[dr]_{p\otimes1}&&
A^2\otimes X\ar[r]^-{1\otimes g}\rtwocell<\omit>{<3>\nu}&
A^2\ar[r]^-{p}&
A\\
&&A\otimes X\ar[ur]_{p^*\otimes 1}\ar@/_1pc/[rrru]_g
\enddiagram
\end{equation}
When $g=\varphi_X(h)$ 
for some $h\in\M(A\otimes X,A)$, that is 
$g=\theta_X(h)=p(A\otimes h)(p^*\otimes X)$ and $\nu$ is equal to
$$
\diagramcompileto{thm3}
A\otimes X\ar[r]^-{p^*\otimes 1}\ar[dr]_{p^*\otimes 1}&
A^2\otimes X\ar[r]^-{1\otimes p^*\otimes
  1}\dtwocell<\omit>{'\cong\phi^*\otimes1}&
A^3\otimes X\ar[dr]|-{p\otimes1\otimes1}\ar[r]^-{1\otimes1\otimes
  h}\drrtwocell<\omit>{'\cong\ c}&
A^3\ar[dr]_{p\otimes 1}\ar[r]^-{1\otimes p}&
A^2\ar[r]^-p\dtwocell<\omit>{'\cong\ \phi^{-1}}&
A\\
&A^2\otimes X\ar[ur]|-{p^*\otimes 1\otimes
  1}\rrtwocell<\omit>{<-3>}\ar@{=}[rr]&&
A^2\otimes X\ar[r]_-{1\otimes h}&
A^2\ar[ur]_p 
\enddiagram
$$
then \eqref{thm1.eq} is equal to the pasting of $\phi^{-1}:p(A\otimes
p)\Rightarrow p(p\otimes A)$ with the following 2-cell
$$
\diagramcompileto{thm4}
A^2\otimes X\ar[r]^-{p\otimes 1}\rruppertwocell<8>^1{}&
A\otimes X\ar[r]_-{p^*\otimes1}\ar[dr]_{p^*\otimes 1}&
A^2\otimes X\dtwocell<\omit>{'\cong\ \phi^*\otimes 1}\ar[r]^-{1\otimes
  p^*\otimes1} &
A^3\otimes X\ar[r]^-{A^2\otimes h}\ar[dr]|-{p\otimes 1\otimes
  1}\drrtwocell<\omit>{'\cong}&
A^3\ar[dr]|-{p\otimes 1}&\\
A\otimes X\urtwocell<\omit>{<-2>\cong}\ar[u]^{1\otimes j\otimes 1}\ar[ur]&
&
A^2\otimes X\rrtwocell<\omit>{<-3>}\ar[ur]|-{p^*\otimes1\otimes1}\ar@{=}[rr]&&
A^2\otimes X\ar[r]_-{1\otimes h}&
A^2
\enddiagram
$$
which is nothing but $\omega\otimes X$ pasted on the right with an isomorphism,
and so it is itself an isomorphism.

Now we show that (5) implies (1).
Suppose that $\lambda_{\Ac}$ is an equivalence.
Define a 1-cell $b=A\otimes
\Ac\xrightarrow{p^*\otimes1}A^2\otimes\Ac\xrightarrow{1\otimes\e} A$; it
has a structure of $\theta_{\Ac}$-algebra $\theta_{\Ac}(b)\Rightarrow b$
given by
$$
\diagramcompileto{thm5}
A\otimes \Ac\ar[r]^-{p^*\otimes1}\ar[dr]_{p^*\otimes1}&
A^2\otimes \Ac\dtwocell<\omit>{'\cong}\ar[r]^-{1\otimes p^*\otimes 1}&
A^3\otimes \Ac\ar[r]^-{A^2\otimes\e}\ar[dr]|-{p\otimes1\otimes1}&
 A^2\ar[r]^-{p}\dtwocell<\omit>{'\cong}&
A\\
&A^2\otimes\Ac\ar[ur]|-{p^*\otimes1\otimes1}\rrtwocell<\omit>{<-3>}\ar@{=}[rr]&&
A^2\otimes \Ac\ar[ur]_{1\otimes \e}
\enddiagram
$$
Denote by $d:\Ac\to A$ a 2-cell corresponding (up to isomorphism) 
to $b$; that is, $b\cong p(A\otimes d)$ in $\M(A\otimes
\Ac,A)^{\theta_{\Ac}}$. Then we have
\begin{align*}
(p\otimes A)(A\otimes d\otimes A)(A\otimes \n)&\cong
(A\otimes\e\otimes A)(p^*\otimes\Ac\otimes A)(A\otimes\n)\\
&\cong (A\otimes\e\otimes A)(A\otimes A\otimes \n)p^*\\
&\cong p^*
\end{align*}
showing that $d$ is a left dualization for $A$. 
\end{proof}
If $A$ has a right bidual the 2-functor $A\otimes -$ has right
biadjoint given by $[A,-]=A^\circ\otimes-$ (see the discussion on
biduals at the beginning of the section).
In this
case, the monad $t$ of \eqref{eq:t1} can be expressed as
$$
t:A^\circ\otimes A\xrightarrow{1\otimes \n\otimes 1}A^\circ\otimes
A^\circ\otimes A\otimes
A\xrightarrow{(p^*)^\circ\otimes1\otimes1}A^\circ\otimes A\otimes
  A\xrightarrow{1\otimes p}A^\circ\otimes A
$$
or
\begin{multline}
  \label{eq:twhenbiduals}
  A^\circ\otimes A\xrightarrow{\n\otimes1\otimes1}A^\circ\otimes
  A\otimes A^\circ\otimes A\xrightarrow{1\otimes p^*\otimes1\otimes1}A^\circ\otimes
  A\otimes A\otimes A^\circ\otimes A\to\\\xrightarrow{1\otimes
    1\otimes\e\otimes1}A^\circ\otimes A\otimes A\xrightarrow{1\otimes
    p}A^\circ\otimes A
\end{multline}
(we omitted the canonical equivalence $A^\circ\otimes
A^\circ\simeq(A\otimes A)^\circ$), and the 1-cell $\ell$ in
\eqref{lamrep.eq} as
$$
A\xrightarrow{(j^*)^\circ\otimes 1}A^\circ\otimes A\xrightarrow{\hole
f\hole}(A^\circ\otimes A)^t.
$$
The 1-cell \eqref{[1,p]i.eq} can
be expressed as $(A^\circ\otimes p)(\n\otimes A):A\to A^\circ\otimes
A\otimes A\to A\otimes A$. 
Recall that this 1-cell has a canonical $t$-algebra structure,
described in Observation \ref{o:l}. 

\begin{thm}\label{laut=hmc.thm}
For any
map pseudomonoid $A$ with right bidual the following are equivalent.
\begin{enumerate}
\item
$A$ is left autonomous.
\item
$A$ has a Hopf module construction
provided by 
\begin{equation}\label{Hopfaut.eq}
A\xrightarrow{\n\otimes 1}\Ac\otimes A\otimes A\xrightarrow{1\otimes
p}\Ac\otimes A.
\end{equation}
\end{enumerate}
Moreover, in this case the dualization is given by $A^\circ
\xrightarrow{1\otimes j}A^\circ \otimes A\xrightarrow{f} A$, where $f$
is left adjoint to \eqref{Hopfaut.eq} and thus, by Corollary
\ref{EMKleisli.cor},  a Kleisli 
construction for the monad $t$.
\end{thm}
\begin{proof}
By Corollary \ref{[1,p]i.cor}, \eqref{Hopfaut.eq} is a Hopf module construction
for $A$ if and only if the theorem of Hopf modules holds for $A$, and
this is equivalent to the existence of a left dualization by Theorem
\ref{maintheorem.th}. 
\end{proof}

As we already mentioned at the beginning of the section, 
to provide a 1-cell $d:\Ac\to A$ with a structure of a dualization is
to provide an adjunction $p(d\otimes A)\dashv (\Ac\otimes p)(\n\otimes A)$.

\begin{prop}\label{p:HMforleftaut}
For a left autonomous map pseudomonoid $A$ 
the adjunction $p(d\otimes
A)\dashv (A^\circ \otimes p)(\n\otimes A)$ induces the monad $t$. 
Moreover, this adjunction is monadic.  
\end{prop}
\begin{proof}
By Corollary \ref{[1,p]i.cor} we know that $(A^\circ \otimes p)(\n\otimes
A):A\to A^\circ \otimes A$ provides an Eilenberg-Moore construction
for $t$.  
\end{proof}

By definition \cite{dualizations}, a {\em right dualization}\/
  $d':A^\vee\to A$
  for a pseudomonoid $A$ in \M\ is a left dualization for $A$ in
  $\M^{\mathrm{rev}}$, \M\ with the reverse tensor product (or the
  opposite tricategory, when we think of \M\ as a one-object tricategory).
In particular, $A^\vee$ is a left bidual for $A$. A pseudomonoid
  equipped with a right dualization is called {\em right autonomous}\/
  and a left and right autonomous pseudomonoid is simply called {\em
  autonomous}. A left autonomous {\em map}\/ pseudomonoid with
  dualization $d$ is autonomous if and only if $d$ is an
  equivalence \cite[Propositions 1.4 and 1.5]{dualizations}.

Recall from \cite{street.ftm} that, given  monads $s$ on $X$ and $s'$
on $X'$, a morphism of monads is a pair $(f,\phi)$ where
$f:X\to X'$ is a 1-cell and $\phi:s'f\Rightarrow fs$ is a 2-cell
compatible with the multiplications and units; these compatibility
conditions can be found in Section \ref{existence.s}. 
With the obvious definition for the 2-cells, we have a 2-category of
monads in a given 2-category $\mathscr K$,
denoted by $\mathrm{Mnd}(\mathscr K)$;
a morphism of monads $(f,\phi)$ 
is an equivalence in $\mathrm{Mnd}(\mathscr K)$ precisely
when $f$ is an equivalence and $\phi$ is invertible.

\begin{cor}
Suppose that $A$ is an autonomous map pseudomonoid.
Then there exists an equivalence of monads
$$
\diagramcompileto{tp^*p}
\Ac\otimes A\ar[d]_{d\otimes1}\ar[rr]^-t\drrtwocell<\omit>{'\cong}&&
\Ac\otimes A\ar[d]^{d\otimes 1}\\
A\otimes A\ar[rr]_{p^*p}&&
A\otimes A
\enddiagram
$$
and, moreover, $p^*:A\to A\otimes A$ is monadic.
\end{cor}
\begin{proof}
The first assertion is clear since $d$ is an equivalence and $t$ is
induced by $p(d\otimes A)\dashv (d^*\otimes A)p^*$; see Proposition
\ref{laut=hmc.thm}. By the same theorem, $(d^*\otimes A)p^*$ is
monadic, and then so is $p^*$  since $d$ is an equivalence. 
\end{proof}

\begin{prop}\label{dmonoidal.prop}
For any left autonomous map pseudomonoid $A$ the left dualization
$d:\Ac\to A$ has the structure of a strong monoidal morphism from
$(\Ac,(j^*)^\circ,(p^*)^\circ)$ to $(A,j,p)$.
\end{prop}
\begin{proof}
It is enough to show that 
\begin{equation}\label{dmonoidal.eq}
\Ac\xrightarrow{d} A \xrightarrow{\n\otimes
  1}\Ac\otimes A\otimes A\xrightarrow{1\otimes p}\Ac\otimes A
\end{equation} 
is strong  monoidal, since $(\Ac\otimes p)(\n\otimes A)$ is an
Eilenberg-Moore object in the 2-category $\mathbf{Opmon}(\M)$.
In the proof of  Theorem \ref{maintheorem.th} we saw that $p(A\otimes
d)\cong (A\otimes \e)(p^*\otimes \Ac)$, so we have to show that
$(\Ac\otimes A\otimes \e)(\Ac \otimes p^*\otimes \Ac)(\n\otimes\Ac)$ is a
strong monoidal 
morphism, or equivalently, by Proposition \ref{opacteqopmon.p}, that
  $(A\otimes \e)(p^*\otimes 
\Ac):A\otimes\Ac\to A$ is a right pseudoaction of $\Ac$ on $A$ ({\em i.e.}, a
$(-\otimes \Ac)$-pseudoalgebra structure on $A$). This itself turns to
be equivalent to say that $p^*:A\to A\otimes A$ is a right
pseudocoaction of $A$ on $A$ ({\em i.e.}, a $(-\otimes
A)$-pseudocoalgebra structure on $A$), which is obviously true.
\end{proof}

We finish this section with some comments on autonomous monoidal lax functors.
The notion of {\em right autonomous monoidal lax functor}\/ was
introduced in \cite{dualizations}, and it consists of a monoidal lax
functor equipped with structure necessary to ensure that it preserves,
in some lax sense, right biduals. Another way of looking at this
concept is as a monoidal lax functor with extra structure 
such that when we take the domain Gray monoid
as the unit Gray monoid ({\em i.e.}, the Gray monoid whose only cells are
identities), then we get a {\em left}\/ autonomous pseudomonoid.  
More explicitly, a right autonomous monoidal lax functor is a monoidal
lax functor $F$ equipped with a pseudonatural transformation
$\kappa_X:(FX)^\circ\to F(X^\circ)$ and modifications
$$
{\xymatrixrowsep{.5cm}
\xymatrixcolsep{.3cm}
\diagramcompileto{autmonlaxfun1}
F(X)\otimes
F(X)^\circ\ar[rr]^-{1\otimes\kappa_X}\ar[d]_{\e}\drrtwocell<\omit>{^\xi}
&&
F(X)\otimes F(X^\circ)\ar[d]^{\chi_{X,X^{\circ}}}\\
I\ar[rd]_{\iota} &&
F(X\otimes X^\circ)\ar[dl]^{F\e}\\
&F(I)
\enddiagram}
{\xymatrixrowsep{.5cm}
\xymatrixcolsep{.3cm}
\diagramcompileto{autmonlaxfun2}
F(X)^\circ\otimes F(X)\ar[rr]^{\kappa_X\otimes 1}&& 
F(X^\circ)\otimes F(X)\ar[d]^{\chi_{X^\circ, X}}\\
I\ar[u]^{\n}\ar[dr]_\iota\rrtwocell<\omit>{\zeta}&&
F(X^\circ\otimes X)\\
&I\ar[ur]_{F\n}
\enddiagram}
$$
satisfying two axioms.

What is proved in \cite{dualizations} is that if $F:\M\to\mathscr N$
is a monoidal special lax functor and $A$ is a left autonomous
pseudomonoid in \M\ with left dualization $d$, then $F(A)$ is left
autonomous with left dualization $F(d)\kappa_A:F(A)^\circ\to
F(A)$. The term {\em special}\/ means that $F$ is normal (in the sense
that the constraint $1_{FX}\to FX$ is an isomorphism for all $X$) and 
the constraints
$(Fg)(Ff)\Rightarrow F(gf)$ are isomorphisms 
whenever $f$ is a map. Special lax functors have the property of
preserving adjunctions. 

If we restrict ourselves to map pseudomonoids, as 
application of Theorem \ref{maintheorem.th}, we can deduce the
following result.

\begin{prop}\label{Fmonspecial.p}
Let $F:\M\to \mathscr N$ be a monoidal special lax functor between
right autonomous Gray monoids and $A$ be a left autonomous map
pseudomonoid in \M. Assume $F$ has the following two
properties:
the monoidal constraints $\iota:I\to FI$ and $\chi_{A,A}:F(A)\otimes
F(A)\to F(A\otimes A)$ are maps, and the 2-cell below is invertible.
$$
\diagramcompileto{specialfun}
&F(A)^3\ar[r]^-{Fp\chi\otimes
  1}\ar[dr]|-{1\otimes Fp\chi}&
F(A)^2\ar@{=}[rr]\ar[dr]^{Fp\chi}\dtwocell<\omit>{'\cong}
\rrtwocell<\omit>{<3>\eta}&& 
F(A)^2\\
F(A)^2 \ar[ur]^{1\otimes(Fp\chi)^*}\ar@{=}[rr]
\rrtwocell<\omit>{<-3>\hole\hole 1\otimes\varepsilon} &&
F(A)^2\ar[r]_-{Fp\chi}&
F(A)\ar[ur]_{(Fp\chi)^*}
\enddiagram
$$
Then, the map pseudomonoid $F(A)$ is left autonomous with left
dualization
\begin{equation}\label{Fmonspecial.e}
F(A)^\circ\xrightarrow{(Fj)\iota\otimes1}F(I)\otimes
F(A)^\circ\xrightarrow{\chi^* (Fp^*)\otimes 1}
F(A)^2\otimes F(A)^\circ\xrightarrow{1\otimes \e}F(A).
\end{equation}
\end{prop}
\begin{proof}
Recall that $F(A)$ has multiplication $F(p)\chi:F(A)\otimes F(A)\to
F(A)$ and unit $F(j)\iota:I\to F(A)$, so that it is a map
pseudomonoid. 
Using the conditions above plus the fact that \eqref{gamma.eq} is
invertible, it can be shown that the corresponding 2-cell
\eqref{gamma.eq} for  $F(A)$ is invertible, and hence $F(A)$ is
left autonomous. The formula for the left dualization is just the
general expression of any left dualization in terms of the product,
unit and evaluation. 
\end{proof}

If $F$ is strong monoidal (sometimes called weak monoidal) in the
sense that $\iota$ and $\chi$ are equivalences, then $F$ preserves
biduals; more explicitly,  there exists
$\kappa:F(A)^\circ\to F(A^\circ)$, unique up to isomorphism, such that 
\begin{equation}\label{kappadef.e}
(F(A)\otimes F(A)^\circ\xrightarrow{1\otimes\kappa}F(A)\otimes
F(A^\circ)\xrightarrow{\chi_{A,A^\circ}}F(A\otimes
A^\circ)\xrightarrow{F\e}FI)
\cong \iota\e, 
\end{equation}
and $\kappa$ is {\em a fortiori}\/ an equivalence. 

\begin{prop}\label{p:monsplaxpreslaut}
Suppose $F:\M\to \mathscr N$ is a  strong monoidal 
special lax functor between Gray monoids and $A$ is a left autonomous
map pseudomonoid in \M\ with left dualization $d$. Then $FA$ is a left
autonomous map pseudomonoid too, with left dualization
$(Fd)\kappa:(FA)^\circ\to F(A^\circ)\to FA$.
\end{prop}
\begin{proof}
The fact that \eqref{gamma.eq} is invertible and that 
$\chi:F(A)\otimes F(A)\to F(A\otimes A)$ is an
equivalence ensures that the hypotheses of Proposition
\ref{Fmonspecial.p} are satisfied, and hence $F(A)$ is left autonomous. The 
formula for the dualization follows from \eqref{Fmonspecial.e} using
\eqref{kappadef.e} and the fact that $\chi$ is an equivalence. 
\end{proof}

Note that although losing some generality, we gain in simplicity by
restricting to the case of left autonomous {\em map}\/ pseudomonoids,
in what our proofs are not based on big diagrams but on the theory of
Hopf modules. 

The following is a simple corollary about left autonomous map
pseudomonoids in a braided Gray monoid. For background about braidings
in a Gray monoid see Section \ref{s:thelaxcentre} and references
therein. 

\begin{cor}\label{dulizforAotimesB.c}
If $A$ and $B$ are left autonomous map pseudomonoids, with left
dualizations $d_A$ and $d_B$ respectively, in a braided Gray monoid
\M, then $A\otimes B$ is a left autonomous map pseudomonoid too, with
left dualization
$$
B^\circ\otimes A^\circ\xrightarrow{c_{B^\circ,A^\circ}}A^\circ\otimes
B^\circ\xrightarrow{(d_A\otimes 1)(1\otimes d_B)}A\otimes B.
$$
\end{cor}
\begin{proof}
As pointed out in \cite{monbicat}, the tensor product
$\otimes:\M\times\M\to M$ is a strong monoidal pseudofunctor with
$\chi_{(X,Y),(Z,W)}=1\otimes c_{Y,Z}\otimes1:X\otimes Y\otimes Z\otimes W\to
X\otimes Z\otimes Y\otimes W$, and $\iota$ equal to the identity. Then,
it is easy to show that, if we take $B^\circ\otimes A^\circ$ as right
bidual for $A\otimes B$, the corresponding 1-cell $\kappa$ is just
$c_{B^\circ,A^\circ}$.
\end{proof}

\section{Frobenius and autonomous map pseudomonoids}\label{s:Frob}

In this section we study the relationship between autonomous
pseudomonoids, the condition \ref{maintheoremgamma.it} in Theorem
\ref{maintheorem.th} and Frobenius pseudomonoids. 
In \cite{Day-Street:quantumcat} it is shown that any autonomous
pseudomonoid is Frobenius, and 
we showed in Theorem \ref{maintheoremgamma.it} that autonomy is
equivalent to the invertibility of the 2-cell $\gamma$ in
\eqref{gamma.eq} and its dual, {\em   i.e.}, the
corresponding 2-cell $\gamma'$ in $\M^{\mathrm{rev}}$. We show 
a converse in absence of biduals, namely: if $\gamma$ and $\gamma'$
are invertible, then $A$ is Frobenius, and as such it has right and
left bidual,  and moreover $A$ is autonomous.

A {\em Frobenius structure}\/ for a pseudomonoid $A$ is a 1-cell
$\varepsilon:A\to I$ such that $\varepsilon p:A\otimes A\to I$ is the
evaluation of a bidual pair; 

\begin{lem}\label{frobenius1.l}
Let $A$ be a pseudomonoid whose multiplication $p$ is a map, 
and call 
$\gamma$ and $\gamma'$, respectively, the following  2-cells.
$$
\diagramcompileto{gammaagain}
&A^3\ar[r]^-{p\otimes1}\ar[dr]|-{1\otimes p}
\drrtwocell<\omit>{'\cong\ \phi}&
A^2\ar[dr]^p\ar@{=}[rr]\rrtwocell<\omit>{<3>\eta}&&
A^2\\
A^2\ar[ur]^{1\otimes
  p^*}\ar@{=}[rr]\rrtwocell<\omit>{<-3>1\otimes\varepsilon\hole\hole\hole\hole\hole\hole} &&
A^2\ar[r]_-p&
A\ar[ur]_{p^*}
\enddiagram
\!\!\!
\diagramcompileto{gammaagain2}
&A^3\ar[r]^-{1\otimes p}\ar[dr]|-{p\otimes 1}
\drrtwocell<\omit>{'\cong\ \phi^{-1}}&
A^2\ar[dr]^p\ar@{=}[rr]\rrtwocell<\omit>{<3>\eta}&&
A^2\\
A^2\ar[ur]^{p^*\otimes
  1}\ar@{=}[rr]\rrtwocell<\omit>{<-3>\varepsilon\otimes 1 \hole\hole\hole\hole\hole\hole} &&
A^2\ar[r]_-p&
A\ar[ur]_{p^*}
\enddiagram
$$
Then the following equalities hold
$$
\diagramcompileto{froblemma1}
A^4\ar[rr]^-{1\otimes p\otimes1}\rrtwocell<\omit>{<3>1\otimes\gamma\hole\hole\hole\hole\hole\hole\hole}&&
A^3&\\
A^3\ar[u]^{A^2\otimes p}\rruppertwocell<0>^{1\otimes p}{<3>\gamma'}&&
A^2\ar[u]|-{1\otimes p^*}\ar@{}[r]|-\cong&
A^2\ar[ul]_{p^*\otimes1}\\
A^2\ar[u]^{p^*\otimes 1}\ar[rr]_-p&&
A\ar[ur]_{p^*}\ar[u]^{p^*}
\enddiagram
=
\diagramcompileto{froblemma2}
&A^4\rruppertwocell<0>^{1\otimes p\otimes
1}{<3>\hole\hole\hole\gamma'\otimes1}&&
A^3\\
A^3\ar[ur]^{A^2\otimes p^*}\ar@{}[r]|-{\cong}&
A^3\ar[u]_{p^*\otimes A^2}\rruppertwocell<0>^{p\otimes 1}{<3>\gamma}&&
A^2\ar[u]_{p^*\otimes 1}\\
&A^2\ar[rr]_-p\ar[ul]^{p^*\otimes 1}\ar[u]_{1\otimes p^*}&&
A\ar[u]_{p^*}
\enddiagram
$$
$$
\xymatrixcolsep{.65cm}
\diagramcompileto{froblemma3}
A^4\rruppertwocell<0>^{A^2\otimes
p}{<3>1\otimes\gamma'\hole\hole\hole\hole\hole\hole} &&
A^3\rruppertwocell<0>^{p\otimes1}{<3>\gamma}&&
A^2\\
A^3\ar[u]^{1\otimes p^*\otimes 1}\ar[rr]^-{1\otimes
p}\ar[drr]_{p\otimes1}&&
A^2\ar[u]|-{1\otimes p^*}\ar[rr]^-p\ar@{}[d]|-{\cong}&&
A\ar[u]_{p^*}\\
&&A^2\ar[urr]_p
\enddiagram
=
\diagramcompileto{froblemma4}
&&A^3\ar[drr]^{p\otimes1}\ar@{}[d]|-\cong&&\\
A^4\ar[urr]^{A^2\otimes p}\rruppertwocell<0>^{p\otimes
A^2}{<3>\gamma\otimes1\hole\hole\hole\hole\hole\hole}&& 
A^3\rruppertwocell<0>^{1\otimes p}{<3>\gamma'}&&
A^2\\
A^3\ar[u]^{1\otimes p^*\otimes 1}\ar[rr]_-{p\otimes 1}&&
A^2\ar[rr]_-p\ar[u]|-{p^*\otimes 1}&&
A\ar[u]_{p^*}
\enddiagram
$$
\end{lem}
\begin{proof}
The proof is a standard calculation involving mates and the axioms of
a pseudomonoid.
\end{proof}

\begin{prop}
Suppose $A$ is a map pseudomonoid and that the 2-cells $\gamma $ and
$\gamma'$ in Lemma \ref{frobenius1.l} are invertible. Then 
$j^*p:A\otimes A\to I$ and $p^*j:I\to A\otimes A$ have the structure
of a bidual pair. In particular, $A$ is a Frobenius pseudomonoid and 
given a choice of right and left
biduals, $A$ is autonomous.
\end{prop}
\begin{proof}
  The 2-cells 
$$
(j^*\otimes A)(p\otimes A)(A\otimes p^*)(A\otimes
  j)\xrightarrow{(j^*\otimes A)\gamma(A\otimes j)}
 (j^*\otimes   A)p^*p(A\otimes j)\cong 1_A
$$
$$
(A\otimes j^*)(A\otimes p)(p^*\otimes A)(j\otimes
A)\xrightarrow{(A\otimes j^*)\gamma'(j\otimes A)}
(A\otimes j^*)p^*p(j\otimes A)\cong 1_A 
$$
endow $j^*p$ and $p^*j$ with the structure of a bidual pair. The
axioms of a bidual pair follow form Lemma \ref{frobenius1.l}.
\end{proof}

\begin{obs}\label{autfrob.obs}
In the hypothesis of the proposition above, different choices of a
bidual for $A$ give rise to different dualizations. For example, when we take
the bidual pair $j^*p, p^*j$, so that $A$ is right and left bidual of
itself, the resulting left and right dualizations are just the
identity $1_A$. 
Slightly more generally, given any equivalence $f:B\to A$, $B$ has a canonical
structure of right bidual of $A$ such that the corresponding left
dualization is (isomorphic) to $f$. To see this just consider the
evaluation $j^*p(A\otimes f):A\otimes B\to I$ and the coevaluation
$(f^*\otimes A)p^*j:I\to B\otimes A$.
\end{obs}

\section{Hopf modules and the centre construction}\label{centre.s}

The most classical notion of the centre of an algebraic structure is the
centre of a monoid. If $M$ is a monoid, its centre is the set of
elements of $M$ with the {\em property}\/ of commuting with every
element of $M$. Anyone would agree if we slightly change our point of
view and said that the centre of $M$ is the set whose elements are
pairs $(x, (x\cdot-)=(-\cdot x))$: elements of $x\in M$ equipped with the
extra {\em structure}\/ of an equality between the multiplication with $x$ on the
left and on the right. The centre of a monoidal category, defined in
\cite{braidedtencat}, follows the spirit of the latter: from the
algebraic structure of a monoidal category \C\ one forms a new algebraic
structure $Z\C$, called the centre of \C. What we actually have is a functor $Z\C\to\C$, and $Z\C$ has a monoidal structure such
that this functor is strong monoidal. Moreover, $Z\C$ has a canonical
braiding. The objects of $Z\C$ are pairs $(x,\gamma_x)$ where
$\gamma_x:(-\otimes x)\Rightarrow(x\otimes-)$ is an invertible natural
transformation. In this context one can also consider the {\em lax
  centre}\/ of \C, simply by dropping the requirement of the
invertibility of $\gamma_x$. See Example \ref{ex:centremonVcat}. 
The functor $Z\C\to\C$ is the universal one satisfying certain
commutation properties, as we shall see later.

Another centre-like object classically considered is the Drinfel'd
double of a finite-dimensional Hopf algebra, or, more recently, of a (co)quasi-Hopf
algebra. See \cite{Majid:QDQHA,Schauenburg:HopfMods}.
Here the concept is not the one of the
object classifying  maps with certain commutation properties, but it is
a representational one.
Roughly speaking, the Drinfel'd double of a finite 
dimensional Hopf algebra $H$ is a Hopf algebra $D(H)$ such that 
the category of representations of $D(H)$ is monoidally equivalent 
to the centre of the category of representations of $H$. 

In this section we study centres and lax centres of autonomous
pseudomonoids by means of the theory of Hopf modules developed in the
previous sections. When applied to the bicategory of
comodules, this approach proves the existence of the centre of a
finite dimensional coquasi-Hopf algebra (considered as a pseudomonoid)
and, moreover, this centre is equivalent to the 
Drinfel'd double (see Section \ref{comodules.sec}). When applied to the
bicategory of \V-modules,  
we see that left autonomous promonoidal \V-categories always have lax
centres (see Section \ref{VMod.sec}).

\subsection{The lax centre}\label{s:thelaxcentre}

We shall work in a {\em braided Gray monoid}, in the sense of
\cite{monbicat}. A braided Gray monoid is a Gray monoid equipped with
pseudonatural equivalences $c_{X,Y}:X\otimes Y\to Y\otimes X$ and
invertible 2-cells 
$$(X\otimes c_{W\otimes Y,Z})(c_{W\otimes X}Y\otimes
Z)\cong (c_{W,X\otimes Z}\otimes Y)(W\otimes X\otimes c_{Y,Z}).
$$
satisfying three axioms. 
These axioms imply that the tensor pseudofunctor
$\otimes:\M\times\M\to\M$ is (strong) monoidal.


The centre of a pseudomonoid was defined in
\cite{Street:monoidalcentre}. Here we  will be interested in the lax
version of the centre, called the {\em lax centre}\/ of a
pseudomonoid. The definition is exactly the same as that of the centre
but for the fact that we drop the requirement of the invertibility of
certain 2-cells. 

\begin{defn}\label{d:CPl}
Given a pseudomonoid in a
braided Gray monoid \M\
define for each object $X$  a category ${CP}_\ell(X,A)$. The
objects, called {\em lax centre pieces},  are pairs $(f,\gamma)$ where
$f:X\to A$ is a 1-cell and $\gamma$ is a 2-cell
\begin{equation}\label{gammacentre1.eq}
\diagramcompileto{cpl}
A\otimes X\ar[d]_{1\otimes f}\ar@{<-}[rr]^{c_{X,A}}&
{}\ddtwocell<\omit>{\gamma}&
X\otimes A\ar[d]^{f\otimes 1}\\
A\otimes A\ar[dr]_p&&
A\otimes A\ar[dl]^p\\
&A
\enddiagram
\end{equation}
satisfying axioms \eqref{asscentre.eq} and \eqref{unitcentre.eq} in
Figure \ref{fig:centrepiece}. The arrows $(f,\gamma)\to(f',\gamma')$
are the 2-cells $f\Rightarrow g$ which are compatible with $\gamma$
and $\gamma'$ in the obvious sense.

This is the object part of a
  pseudofunctor $CP_\ell(-,A):\M^{\mathrm{op}}\to \mathbf{Cat}$, that
  is defined on 1-cells and 2-cells just by precomposition.  
When
  $CP_\ell $ is birepresentable we call a birepresentation $z_\ell:Z_\ell
  A\to A$ a {\em
  lax centre}\/ of the pseudomonoid $A$.  

A {\em centre piece}\/ is a lax centre piece $(f,\gamma)$ such that
$\gamma$ is invertible. The full subcategories $CP(X,A)\subset
CP_\ell(X,A)$ with objects the centre pieces define a pseudofunctor
$CP(-,A):\M^{\mathrm{op}}\to\mathbf{Cat}$, and we call a
birepresentation of it a {\em centre}\/ of $A$, denoted by $z:ZA\to A$. 
\end{defn}
\begin{figure}\label{axiomslaxcentre.fig}
\begin{equation*}
\diagramcompileto{cplax1a}
A\otimes A\otimes X\ddrtwocell<\omit>{'\cong}\ar[d]_{1\otimes 1\otimes f}
\ar[dr]^{p\otimes 1}\ar@{<-}[rrrr]^-{c_{X,A\otimes
    A}}&&
{}\dtwocell<\omit>{'\cong}&&
X\otimes A\otimes A\ar[dl]_{1\otimes p}\ar[d]^{f\otimes 1\otimes 1}
\ddltwocell<\omit>{'\cong} \\
A\otimes A\otimes A\ar[dr]_{p\otimes1}\ar[d]_{1\otimes p}&
A\otimes X\ar@{<-}[rr]^-{c_{X,A}}\ar[d]|-{1\otimes f}&
{}\ddtwocell<\omit>{\gamma}&
X\otimes A\ar[d]|-{f\otimes 1}&
A\otimes A\otimes A\ar[d]^{p\otimes1}\ar[dl]^{1\otimes p}\\
A\otimes A\ar[drr]_{p}\rtwocell<\omit>{'\cong}&
A\otimes A\ar[dr]^-{p}&&
A\otimes A\ar[dl]_{p}\rtwocell<\omit>{'\cong}&
A\otimes A\ar[dll]^p\\
&&A&&
\enddiagram
\end{equation*}
\begin{equation}\label{asscentre.eq}
\parallel
\end{equation}
\begin{equation*}
\xymatrixcolsep{.8cm}
\diagramcompileto{cplax1b}
{A\otimes A\otimes X}\ar@{<-} `u `[rrrr]^{c_{X,A\otimes A}} [rrrr] 
\ar@{<-}[rr]^-{1\otimes c_{X,A}}\ar[d]_{1\otimes1\otimes
  f}&
{}\ddtwocell<\omit>{1\otimes \gamma}\ar@{}@<10pt>[rr]|-\cong&
A\otimes X\otimes A\ar[d]|-{1\otimes f\otimes
  1}\ar@{<-}[rr]^-{c_{X,A}\otimes1}&
{}\ddtwocell<\omit>{\gamma\otimes 1}&
X\otimes A\otimes A\ar[d]^{f\otimes 1\otimes 1}\\
A\otimes A\otimes A\ar[dr]_{1\otimes p}&&
A\otimes A\otimes A\ar[dl]|-{1\otimes p}\ar[dr]|-{p\otimes 1}&&
A\otimes A\otimes A\ar[dl]^{p\otimes 1}\\
&A\otimes A\ar[dr]_{p}\rrtwocell<\omit>{'\cong}&&
A\otimes A\ar[dl]^p&\\
&&A
\enddiagram
\end{equation*}
\begin{equation}\label{unitcentre.eq}
\diagramcompileto{cplax2}
X\ar@{=}[rrrr]\ar[dr]^{j\otimes1}\ar[dd]_{f}\ddrtwocell<\omit>{'\cong}&
&{}\dtwocell<\omit>{'\cong}&&
X\ar[dl]_{1\otimes j}\ar[dd]^{f}\ddltwocell<\omit>{'\cong}\\
&A\otimes X\ar@{<-}[rr]^-{c_{X,A}}\ar[d]_{1\otimes
  f}&{}\ddtwocell<\omit>{\gamma}& 
X\otimes A\ar[d]^{f\otimes 1}&\\
A\ar[r]^-{j\otimes 1}\ar@/_/[drr]_1\drrtwocell<\omit>{'\cong}&
A\otimes A\ar[dr]_p&&
A\otimes A\ar[dl]^p&
A\ar[l]_-{1\otimes j}\ar@/^/[dll]^1\dlltwocell<\omit>{'\cong}\\
&&A&&
\enddiagram
=1_f
\end{equation}
\caption{Lax centre piece axioms}\label{fig:centrepiece}
\end{figure}

\begin{defn}\label{d:zc}
The inclusion $CP(-,A)\hookrightarrow CP_\ell(-,A)$ induces a 1-cell
$z_c:ZA\to Z_\ell A$, unique up to isomorphism such that $z_\ell
z_c\cong z$ as centre pieces. 
When $z_c$ is an equivalence we will say that the centre
of $A$ coincides with the lax centre. 
\end{defn}

\begin{exmp}\label{ex:centremonVcat}
  The centre of a pseudomonoid in $\mathbf{Cat}$, that is, of a monoidal
  category, is the usual centre defined in \cite{braidedtencat}. In
  fact, lax centres and centres of pseudomonoids in
  $\V\text{-}\mathbf{Cat}$ exist and are given by the constructions in
  \cite{Day:LaxCentre}. If \A\ is a monoidal \V-category, its lax
  centre $Z_\ell\C$ has objects pairs $(x,\gamma)$ where $x$ is an
  object of \C\ and $\gamma:(-\otimes x)\Rightarrow(x\otimes -)$ is a
  \V-natural transformation. The \V-enriched hom $Z_\ell\C((x,\gamma),
  (y,\delta))$ is the equalizer of the pair of arrows
$${\xymatrixcolsep{1.4cm}
\diagramcompileto{ZmonVcathom}
\C(x,y)\ar[r]\ar[d]&
[\C,\C](-\otimes x,-\otimes y)\ar[d]^{[\C,\C](\gamma,1)}\\
[\C,\C](x\otimes -,y\otimes -)\ar[r]_{[\C,\C](1,\delta)}&
[\C,\C](x\otimes-,-\otimes y)\\
\enddiagram}
$$
\end{exmp}

\begin{obs}
By \cite{Street:monoidalcentre}, in a bicategory with finite
products, iso-inserters and cotensoring with the arrow category {\em
any}\/ pseudomonoid has a centre. 
\end{obs}

We would like to exhibit an equivalence $\M(I,Z_\ell A)\simeq
Z_\ell(\M(I,A))$. Our leading example is the one of the bicategory
$\V\text{-}\mathbf{Mod}$ of \V-categories and \V-modules. For details
about this bicategory see \SorCh\ \ref{VMod.sec}. Henceforth, we shall
assume our Gray monoid \M\ satisfies additional properties, which we explain
below. 

 Recall that a 2-cell 
$$
\diagramcompileto{rlifting}
&Y\ar[d]^g\dllowertwocell<0>_{{}^f\! g}{^<-2>\lambda}\\
X\ar[r]_f&Z
\enddiagram
$$
in a bicategory \B\
is said to {\em exhibit $ {}^f\! g$ as the right lifting of $g$ through $f$}
if it induces a bijection  $\B(Y,X)(k,{}^f\! g)\cong\B(Y,Z)(fk,g)$, natural in
$k$. Clearly, right liftings are unique up to compatible isomorphisms.
See \cite{yonedastructures}. 

We shall assume that our braided Gray monoid \M\ is closed (see
Section \ref{opmon.s} and references therein) and has
{\em right liftings}\/ of arrows 
out of $I$
through arrows out of $I$. As explained in \cite{monbicat}, this endows
each
$\M(X,Y)$ with the structure of a \V-category. Here $\V=\M(I,I)$ is a 
symmetric monoidal closed category whose tensor product is given by
composition. The \V-enriched hom $\M(X,Y)(f,g)$ is ${}^{\hat f}\! \hat g$, the
right lifting of $\hat g:I\to [X,Y]$ through $\hat f:I\to [X,Y]$,
where these two arrows correspond to $f$ and $g$ under the closedness
biadjunction. Both $\hat f$ and $\hat g$ are determined up to
isomorphism, and then so is $\M(X,Y)(f,g)$. The
compositions $\M(X,Y)(g,h)\otimes\M(X,Y)( f,g)\to \M(X,Y)(f,h)$ and
units $1_I\to \M(X,Y)(f,f)$, along 
with the \V-category axioms, are easily deduced from the universal
property of the right liftings.
Observe that the underlying category of the \V-category $\M(X,Y)$ is
the hom-category $\M(X,Y)$. For, $\V(1_I,\M(X,Y)(f,g))=\V(1_I,{}^{\hat
f}\! \hat g)\cong\M(I,[X,Y])(\hat f,\hat g)\cong\M(X,Y)(f,g)$. 

One can define {\em composition \V-functors}\/
$\M(Y,Z)\otimes\M(X,Y)\to\M(X,Z)$ on objects just by composition in
\M\ and on \V-enriched homs in the following way. Given $f,h:Y\to Z$
and $g,k:X\to Y$, define an arrow
$\M(I,[Y,Z])(\hat f,\hat h)\otimes\M(I,[X,Y])(\hat g,\hat
k)\to\M(I,[X,Z])(\widehat{fg},\widehat{hk})$ as the 2-cell in \M\
corresponding to the following pasting. 
$$
\xymatrixcolsep{1.5cm}
\diagramcompileto{Venrichedcomposition}
&&I\ar[dl]_{{}^{\hat g}\! \hat k}\dltwocell<\omit>{^<-2>}\ar[d]^{\hat k}\ar[dddr]^{\widehat{hk}}&\\ 
&I\ar[r]_-{\hat g}\ar[d]^{\hat h}\ar@{}[dr]|-\cong\ar[dl]_{{}^{\hat f}\! \hat
  h} \dltwocell<\omit>{^<-2>}&
[X,Y]\ar[d]^{\hat h\otimes 1}\ar@{}[ddr]|-\cong&\\
I\ar[r]^-{\hat f}\ar[drrr]_{\widehat{fg}} &
[Y,Z]\ar[r]_-{1\otimes \hat g}\ar@{}@<5pt>[drr]|-\cong&
[Y,Z]\otimes [X,Y]\ar[dr]|-{\mathrm{comp}}&\\
&&&[X,Z]
\enddiagram
$$
There are also {\em  identity \V-functors}\/  from the trivial
\V-category to $\M(X,X)$. On objects they just pick the identity
1-cells $1_X$ and homs they are given by the arrows $1_I\to
{}^{(\hat{1_X})}\!\hat{1_X}$ corresponding to the identity 2-cells
$\hat{1_X}\Rightarrow \hat{1_X}$. These composition and identity
\V-functors endow \M\ with the structure of a category weakly enriched
in $\V\text-\mathbf{Cat}$, in the sense that the category axioms hold
only up to specified \V-natural isomorphisms ({\em e.g.} when $\V$ is
the category of sets, we get a bicategory with locally small
hom-categories). 


Now we shall further suppose that the category $\V=\M(I,I)$ is
complete. This allows us to consider functor \V-categories. In this
situation, the composition \V-functors induce \V-functors
$\M(X,-)_{Y,Z}:\M(Y,Z)\to[\M(X,Y),\M(X,Z)]$ making 
the
pseudofunctor $\M(X,-):\M\to\V\text-\mathbf{Cat}$ locally a
\V-functor. 



\begin{lem}\label{l:CPVcat}
  Under the hypothesis above, if $A$ is a pseudomonoid in \M,
  $CP_\ell(I,A)$ has a canonical structure of a \V-category such that
  the forgetful  functor $CP_\ell(I,A)\to\M(I,A)$ is the underlying
  functor of a \V-functor. Moreover, 
  $CP(I,A)$ is a full sub-\V-category of $CP_\ell(I,A)$. 
\end{lem}
\begin{proof}
We give only a sketch of a proof; the details are an exercise
in the universal property of right liftings.  
Given two lax centre pieces $(f,\alpha)$ and $(g,\beta)$, define the
  \V-enriched hom $CP_\ell(I,A)((f,\alpha),(g,\beta))$ as the
  equalizer in \V\ of the pair
\begin{equation}\label{eq:CPVhom}
\xymatrixcolsep{1.6cm}
\diagramcompileto{CPell(I,A)}
\M(I,A)(f,g)\ar[r]\ar[d]&
\M(A,A)(p(A\otimes f),p(A\otimes g))\ar[d]^{\M(A,A)(\alpha,1)}\\
\M(A,A)(p(f\otimes A),p(g\otimes A))\ar[r]_{\M(A,A)(1,\beta)}&
\M(A,A)(p(f\otimes A),p(A\otimes g))
\enddiagram
\end{equation}
where the unlabelled arrows are induced by the universal property of
right liftings under postcomposition with the arrows $A\to [A,A]$
corresponding to $p$ and $pc_{A,A}$. With this definition, an arrow
$1_I\to CP_\ell(I,A)((f,\alpha),(g,\beta))$ in $\V=\M(I,I)$
corresponds to an arrow $(f,\alpha)\to(g,\beta)$ in the ordinary
category $CP_\ell(I,A)$. The composition
$CP_\ell(I,A)((g,\beta),(h,\gamma))\otimes
CP_\ell(I,A)((f,\alpha),(g,\beta)) \to
CP_\ell(I,A)((f,\alpha),(h,\gamma))$ is induced by the composition
$\M(I,A)(g,h)\otimes\M(I,A)(f,g)\to\M(I,A)(f,h)$ and the universal property
of the equalizers, and likewise for the identities. 
\end{proof}


\begin{prop}\label{p:MIZCP}
  Assume the lax centre of $A$ exists, with universal centre piece
  $(z_\ell,\gamma)$. 
  Under the hypothesis above,
  $(z_\ell,\gamma)$ induces a \V-enriched equivalence $U$ making the
  following diagram commute. 
$$
\diagramcompileto{MIZPC}
\M(I,Z_\ell A)\ar[dr]_{\M(I,z_\ell)}\ar[rr]^-U&&
CP_\ell(I,A)\ar[dl]\\
&\M(I,A)&
\enddiagram
$$
  Moreover, the same is true if the centre of $A$ exists and we use
  $CP(I,A)$ instead of $CP_\ell(I,A)$. 
\end{prop}
\begin{proof}
  Let $(z_\ell,\gamma)\in CP_\ell(Z_\ell A,A)$ be the universal lax
  centre piece.
  On objects, $U$ is equal to the usual functor, that is,
  it sends $f:I\to Z_\ell A$ to the lax centre piece $(z_\ell f,\gamma
  (f\otimes A))$. Next we describe our \V-functor $U$ on homs. Define
  $\varrho$ by the following equality, where $\pi$ exhibits ${}^h\! k$ as a
  right lifting of $k$ through $h$ and $\varpi$ exhibits ${}^{(z_\ell
  h)}\! (z_\ell k)$ as
  a right lifting of $z_\ell k$ through $z_\ell h$. 
\begin{equation}\label{eq:varrhopi}
\diagramcompileto{vfunct1}
{}&I\dllowertwocell<\omit>{^<-2>\pi}\ar[dl]_{{}^h\! k}\ar[d]^k&{}\\
I\ar[r]_-h&
Z_\ell A\ar[dr]_-{z_\ell}&\\
&&A
\enddiagram
=
\diagramcompileto{vfunct2}
{}&&I\ar[ddll]|-(.2){{}^{(z_\ell h)}\! (z_\ell
  k)}\ddlllowertwocell<-10>_{{}^h\! k}{^\varrho}\ar[d]^k\\ 
{}&&Z_\ell A\ar[d]^{z_\ell}\dltwocell<\omit>{^\varpi}\\
I\ar[r]_-h&Z_\ell A\ar[r]_-{z_\ell}&A
\enddiagram
\end{equation}
This pasted composite is trivially a morphism of lax centre pieces
$U(h({}^h\! k))\to U(k)$, and this means exactly that $\varrho$ factors
through the equalizer 
$$
CP_\ell(I,A)(U(h),U(k))\rightarrowtail {}^{(z_\ell h)}\! (z_\ell
k)=\M(I,A)(z_\ell h,z_\ell k);
$$
in \eqref{eq:CPVhom} defining $CP_\ell(I,A)(U(h),U(k))$ on \V-enriched homs. 
Denote by $\tilde\varrho:{}^h\! k=\M(I,A)(h,k)\to CP_\ell(I,A)(U(h),U(k))$ the
resulting arrow in \V. This is by definition the effect of $U$ on
enriched homs. 

Observe that the underlying ordinary functor of $U$ is the usual
equivalence  given by the universal property of the lax centre. Hence,
$U$ is essentially surjective on objects as a \V-functor. It is
sufficient, then, to show that $U$ is fully faithful, or, in other
words, that 
$\tilde\varrho$ is invertible. To do this, we shall show that
$\varrho$ has the universal property of the equalizer defining
$CP_\ell(I,A)(U(h),U(k))$. 

Suppose $\nu:v\to {}^{(z_\ell h)}\! (z_\ell k)$ is an arrow in \V\ equalizing the pair of arrows
${}^{(z_\ell h)}\! (z_\ell k)\to\M(A,A)(p(z_\ell h\otimes A),p(A\otimes z_\ell
k))$ analogues to \eqref{eq:CPVhom}. 
If one unravels this condition, one gets the following equality.
$$
{\diagramcompileto{vfunct3}
{}&&A\ddlllowertwocell<-10>_{v\otimes  1}{^\hole\hole\nu\otimes1}
\ar[ddll]|-(.7){{}^{(z_\ell h)}\! (z_\ell k)\otimes 1}\ar[d]|-{k\otimes  1} 
\ar[drr]^{1\otimes k} \ar@{}@<-.8pc>[drr]|-(.3){\cong}&&\\
{}&& Z_\ell A\otimes A \ar[rr]^\simeq\ar[d]|-{z_\ell \otimes 1}&
{}\ddtwocell<\omit>{^\gamma}&
A\otimes Z_\ell A\ar[d]^{1\otimes z_\ell }\\
A\ar[r]_-{h\otimes 1}&
 Z_\ell A\otimes A
\ar[r]_-{z_\ell \otimes 1}\urtwocell<\omit>{<-1>\hole\hole\varpi\otimes1}&
A^2\ar[dr]_p&&
A^2 \ar[dl]^p \\
&&&A
\enddiagram}
$$
$$
\parallel
$$
$$
\xymatrixcolsep{1.4cm}
{\diagramcompileto{vfunct4}
{}&&A \ddlllowertwocell<-10>_{1\otimes v}{^\hole\hole1\otimes \nu}
\ar[ddll]|-(.7){ 1\otimes {}^{(z_\ell h)}\! (z_\ell k)}\ar[d]^{1\otimes k} &\\
&& A\otimes Z_\ell A\ar[d]^{1\otimes z_\ell}&\\
A \ar@{}@<.7pc>[dr]|-(.7)\cong \ar[r]^-{1\otimes h}\ar[dr]_{h\otimes 1}&
A\otimes Z_\ell A\ar[r]^-{1\otimes z_\ell}
\urtwocell<\omit>{<-1>\hole\hole1\otimes\varpi}&
A^2\ar[r]^p&
A\\
& Z_\ell A\otimes A\ar[r]_-{z_\ell \otimes1}
\ar[u]^\simeq\urrtwocell<\omit>{^\gamma}& 
A^2\ar[ur]_p&
\enddiagram}
$$
This means that the 2-cell $\varpi(z_\ell h\nu)$ is an arrow in the
ordinary category $CP_\ell(I,A)$ from $U(hv)=(z_\ell hv,\gamma((
hv)\otimes A))$ to $U(k)=(z_\ell k,\gamma(k\otimes A))$, and therefore
there exists a 
unique 2-cell $\tau:hv\Rightarrow k:I\to Z_\ell A$ such that
$z_\ell \tau=\varpi(z_\ell h\nu)$. From the universal property of right
liftings, we deduce the existence of a unique $\tau':v\Rightarrow {}^h\! k$ such
that $\pi(h\tau')=\tau$. In order to show that $\varrho:{}^h\! k\Rightarrow
{}^{(z_\ell h)}\!(z_\ell k)$ has the universal property  of the equalizer as
explained above, we have to show that $\varrho\tau'=\nu$. But the
pasting of 
$\varrho\tau'$ with $\varpi$, $\varpi(z_\ell h(\varrho\tau'))$,  is
equal, by definition of $\varrho$, to 
$z_\ell (\pi(h\tau'))=z_\ell \tau=\varpi(z_\ell h\nu)$. It follows that
$\varrho\tau'=\nu$. 

The case of the centre is completely analogous to the case of the lax
centre. The \V-functor $U$ is defined on objects by sending $f:I\to
ZA$ to the centre piece $(zf,\gamma(f\otimes A))$, where $(z,\gamma)$
is the universal centre piece. The definition of $U$ on \V-enriched
homs is the same as in the case of the lax centre above.
\end{proof}

In order to exhibit the desired equivalence $\M(I,Z_\ell A)\simeq
Z_\ell(\M(I,A))$, 
we shall require of our closed braided Gray
monoid \M\ three further  properties. 

Firstly, we require
the monoidal closed category $\V$ be complete. 
This
allows us to talk about functor \V-categories. 

Secondly,
the pseudofunctor
$\M(I,-):\M\to\V\text{-}\mathbf{Cat}$ must be locally 
faithful. In other words, for every pair of 1-cells $f,g$, the
following must be a monic arrow in \V.
\begin{equation}
  \label{eq:M(IA)faith}
\M(X,Y)(f,g)\to[\M(I,X),\M(I,Y)](\M(I,f),\M(I,g))  
\end{equation}

Finally, for any $f,g:X\to Y$, the image of the arrow
\eqref{eq:M(IA)faith}
under
$\V(I,-):\V\to\mathbf{Set}$ must be surjective. This condition is
saying that every \V-natural transformation
$\M(I,f)\Rightarrow\M(I,g)$ is induced by a 2-cell $f\Rightarrow g$;
this 2-cell is unique by the  condition in the previous paragraph.



All these properties are satisfied by our main example of
$\V\text-\mathbf{Mod}$, as we shall see later. 

\begin{thm}\label{t:MZsimeqZM}
  Under the hypothesis above, if $A$ has a lax centre then 
  there exists a \V-enriched equivalence
  making the following diagram commutes up to a canonical isomorphism. 
$$
\diagramcompileto{MZsimeqZM1}
\M(I,Z_\ell A)\ar[rr]^-\simeq\ar[dr]_{\M(I,z_\ell)}&&
Z_\ell(\M(I,A))\ar[dl]^V\\
&\M(I,A)
\enddiagram
$$
Here the \V-category on the right hand side is a lax centre in
$\V\text{-}\mathbf{Cat}$ and $V$ is the forgetful \V-functor. 
Furthermore, the result remains true if we write centres in place of
lax centres.
\end{thm}
\begin{proof}
By Proposition \ref{p:MIZCP} it is enough to exhibit a \V-enriched
equivalence between $CP_\ell(I,A)$ and $Z_\ell(\M(I,A))$ commuting
with the forgetful functors. 

  Define a \V-functor $\Phi:CP_\ell(I,A)\to Z_\ell(\M(I,A))$ as
  follows. On objects $\Phi(f,\alpha)=(f,\Phi_1(\alpha))$ where
  $$
  \Phi_1(\alpha)_h:h*f\cong p(A\otimes f)h\xrightarrow{\alpha h}
  p(f\otimes A)h\cong f*h.
  $$
Recall that the \V-enriched hom $CP_\ell(I,A)((f,\alpha),(g,\beta))$
is the equalizer of 
\eqref{eq:CPVhom} and $Z_\ell(\M(I,A))(\Phi(f,\alpha),\Phi(g,\beta))$ is
the equalizer of the diagram in Example \ref{ex:centremonVcat}, where
$\C=\M(I,A)$, $x=f$, $y=g$, $\gamma=\Phi_1(\alpha)$ and
$\delta=\Phi_1(\beta)$. We can draw a diagram 
$$
\diagramcompileto{MZsimeqZM2}
CP_\ell(I,A)((f,\alpha),(g,\beta))\rightarrowtail\M(I,A)(f,g)
\ar@<2pt>[r] \ar@<-2pt>[r]\ar@<-2pt>[dr]\ar@<2pt>[dr]&
\M(A,A)(p(f\otimes A),p(A\otimes g))\ar[d]^{\M(I,-)}\\
&[\M(I,A),\M(I,A)](f*-,-*g)
\enddiagram
$$
where  $CP_\ell(I,A)((f,\alpha),(g,\beta))$ is the equalizer of the
pair of arrows in the top row and
$Z_\ell(\M(I,A))(\Phi(f,\alpha),\Phi(g,\beta))$ is the equalizer of
the other diagonal pair of arrows. Moreover, the diagram serially
commutes, as the vertical arrow is induced by the effect  of the
pseudofunctor $\M(I,-):\M\to\V\text{-}\mathbf{Cat}$ on
\V-enriched homs, and hence monic by hypothesis. It follows that there
exists an isomorphism $CP_\ell(I,A)((f,\alpha),(g,\beta))\to
Z_\ell(\M(I,A))(\Phi(f,\alpha),\Phi(g,\beta))$. One can check that
these isomorphisms are part of a \V-functor $\Phi$, which, obviously,
is fully faithful. 

It only rests to prove that $\Phi$ is essentially surjective on
objects. Here is where the hypothesis on
$\M(I,-):\M\to\V\text{-}\mathbf{Cat}$  come into play.  
An object $(f,\gamma)$ of $Z_\ell(\M(I,A))$ gives rise to a
\V-natural transformation 
$$
\gamma'_h:p(A\otimes f)h\cong h*f\xrightarrow{\gamma_h}f*h\cong
p(f\otimes A)h.
$$
By hypothesis, $\gamma'$ is induced by a unique $\alpha:p(A\otimes
f)\Rightarrow p(f\otimes A)$. The equalities \eqref{asscentre.eq} and
\eqref{unitcentre.eq} for the 2-cell $\alpha$ 
follow from the fact that $(f,\gamma)$ is an object in the lax centre
of $\M(I,A)$ and the fact that 
$\M(A^2,A)\to [\M(I,A^2),\M(I,A)]$ is fully faithful.
Now observe that $\Phi(f,\alpha)=(f,\gamma)$. Finally, $\alpha$ is
invertible if and only if $\gamma$ is invertible, so that proof
also applies to centres.
\end{proof}

Recall from \cite{dualizations} that for a right autonomous pseudomonoid $A$, with right dualization
$\bar d:A^\vee\to A$, every map  $f:I\to A$ has a right dual in the
monoidal \V-category $\M(I,A)$. A right dual
of $f$ is given by $\bar d (f^*)^\vee$, where $f^*$ is a right adjoint to
$f$. Then the full subcategory $\mathrm{Map}\M(I,A)$ of $\M(I,A)$ is
right autonomous (in the classical sense that it has right duals). 

\begin{thm}\label{t:laxcentre=centre-maps}
In addition  to the hypothesis above, assume the following: 
\V\ is complete and cocomplete monoidal closed category, \M\ has
{\em all} right liftings, the inclusion \V-functor 
$\mathrm{Map}\M(I,A)\to\M(I,A)$ is dense and
$\M(I,-):\M\to\mathbf{Cat}$ reflects equivalences. If $A$ is
left autonomous, then the
centre of $A$ coincides with the lax centre whenever both exist.
\end{thm}
\begin{proof}
  By Theorem \ref{t:MZsimeqZM}, there exists an isomorphism as
  depicted below. 
$$\xymatrixcolsep{1.6cm}
\diagramcompileto{Z=Zl}
\M(I,Z A)\ar[r]^-\simeq\ar[d]_{\M(I,z_c)}\ar@{}[dr]|-\cong&
Z(\M(I,A))\ar@{_(->}[d]\\
\M(I,Z_\ell A)\ar[r]^-\simeq&
Z_\ell(\M(I,A))
\enddiagram
$$
A straightforward modification of \cite[Prop. 6]{monbicat} (using the
property of the right liftings with respect to composition dual to
\cite[Prop. 1]{yonedastructures}) shows that the monoidal \V-category
$\M(I,A)$ is closed as a \V-category. It follows that the \V-functors
$(f*-)=p(f\otimes A)-:\M(I,A)\to\M(I,A)$ 
given by tensoring with an object $f$ are cocontinuous. 
As $\M(I,A)$ has a dense sub \V-category of which every object has a right
dual, the hypotheses of \cite[Theorem
3.4]{Day:LaxCentre} are satisfied, and we deduce that the inclusion
$Z(\M(I,A))\hookrightarrow Z_\ell(\M(I,A))$ is the identity. 
It follows that $\M(I,z_c)$ is an equivalence and
hence $z_c$ is an equivalence.
\end{proof}

\subsection{Lax centres of autonomous pseudomonoids}

The lax centre of a pseudomonoid was defined as a birepresentation of
the pseudofunctor $CP_\ell(-,A)$. An object of the category
$CP_\ell(X,A)$, {\em i.e.}, a lax centre piece, is a 2-cell
$p(f\otimes A)\Rightarrow p(A\otimes f)c_{X,A}$. We observe that the
same notion of lax centre can be defined by using $c^*$ instead of
$c$. In an entirely analogous
way to Definition \ref{d:CPl}, one defines a category $CP^*_\ell(X,A)$
as follows. It has objects $(f,\gamma)$
where $f:X\to A$ and 
$\gamma:p(f\otimes A)c_{X,A}^*\Rightarrow p(A\otimes f)$, and
arrows $(f,\gamma)\to(g,\delta)$ those 2-cells $f\Rightarrow g$
which are compatible with $\gamma$ and $\delta$. Pasting with the
structural isomorphism $c_{X,A}c_{X,A}^*\cong 1_{X\otimes A}$ induces
pseudonatural equivalences $CP_\ell(X,A)\to CP_\ell^*(X,A)$. This is
the reason why the $c^*$ appears in the following definition.

\begin{defn}\label{d:sigma}
Given a map pseudomonoid $A$ in a braided Gray monoid \M\ 
define a pseudonatural transformation $\sigma:\M(A\otimes
-,A)\Rightarrow\M(A\otimes -,A)$ with components 
$$
\sigma_X(g)=\big(A\otimes X\xrightarrow{p^*\otimes 1}A^2\otimes
X\xrightarrow{1\otimes c^*_{X,A}}A\otimes X\otimes
A\xrightarrow{g\otimes 1}A^2\xrightarrow{p}A \big).
$$
\end{defn}
\begin{lem}
The pseudonatural transformation $\sigma$ has a canonical structure of
a monad. 
\end{lem}
\begin{proof}
Just note that $\sigma$ is isomorphic to the monad $\theta$ of Definition
\ref{d:theta} for the map pseudomonoid $(A,j,pc_{A,A}^*)$.
\end{proof}

Explicitly, the multiplication of $\sigma$ is given by components
\begin{equation}\label{sigmamult.eq}
\def\objectstyle{\scriptstyle}
\def\labelstyle{\scriptscriptstyle}
\xymatrixcolsep{.9cm}
{\diagramcompileto{simamult}
A\otimes X\otimes A\ar[r]^-{p^*\otimes 1\otimes 1}&
A^2\otimes X\otimes A\drtwocell<\omit>{'\cong}
\ar@/^1pc/[dr]^{1\otimes c^{*}_{X,A}\otimes 1}&&&&\\
A^2\otimes X\ar[r]^-{p^*\otimes1\otimes 1}\ar[u]|-{1\otimes
c^{*}_{X,A}}\urtwocell<\omit>{'\cong} &
A^3\otimes X\ar[u]|-{A^2\otimes c^{*}_{X,A}}\ar[r]^-{1\otimes
c^*_{X,A^2}}\ar[dr]|-{1\otimes p\otimes 1 }&
A\otimes X\otimes A^2\ar[dr]|-{1\otimes 1\otimes p}\ar[r]^-{g\otimes A^2}
\dtwocell<\omit>{'\cong}&
A^3\ar[r]^-{p\otimes 1}\ar[dr]|-{1\otimes p}\dtwocell<\omit>{'\cong}&
A^2\ar[dr]|-{p}\dtwocell<\omit>{'\cong}&
\\
A\otimes X\ar[u]|-{p^*\otimes 1}\ar[r]_-{p^*\otimes
1}\urtwocell<\omit>{'\cong}&
A^2\otimes X\ar[u]|-{1\otimes p^*\otimes 1}\ar@{=}[r]\rtwocell<\omit>{<-1.5>}&
A^2\otimes X\ar[r]_{1\otimes c^{*}_{X,A}}&
A\otimes X\otimes A\ar[r]_-{g\otimes 1}&
A^2\ar[r]_-{p}&
A
\enddiagram}
\end{equation}
and the unit by
\begin{equation}\label{sigmaunit.eq}
\xymatrixcolsep{.9cm}
\diagramcompileto{sigmaunit}
&A\otimes X\ar@{=}[r]\drtwocell<\omit>{'\cong}&
A\otimes X\ar[r]^-{g}\drtwocell<\omit>{'\cong}&
A\ar[dr]^{1\otimes j}\ar@(r,ul)[drr]^1&&\\
A\otimes X\ar[r]_-{p^*\otimes 1}\ar@/^/[ur]&
A^2\otimes X\ar[r]_-{1\otimes c^*_{X,A}}\ar[u]|-{1\otimes j^*\otimes
  1}\ar@{}[ul]|-(.2){\cong} &
A\otimes X\otimes A\ar[r]_-{g\otimes 1}\ar[u]|-{1\otimes
  1\otimes j^*}&
A^2\ar[u]|-{1\otimes j^*}\ar@{=}[r]\rtwocell<\omit>{<-1.5>}&
A^2\ar[r]_-p&
A\ar@{}[ull]|-(.3){\cong}
\enddiagram
\end{equation}

When $A\otimes -$ has right biadjoint the monads $\theta$ and $\sigma$
are represented by monads $t$ and $s:[A,A]\to [A,A]$. The monad
$s$ is
\begin{equation}\label{sdefinition.eq}
[A,A]\xrightarrow{i_{A}^A}[A\otimes A,A\otimes A]
\xrightarrow{[c_{A,A},c^*_{A,A}]}[A\otimes A,A\otimes
A]\xrightarrow{[p^*,p]}[A,A],
\end{equation}
which is the monad $t$ for the opposite pseudomonoid of
$A$ with respect to $c^*$, in other words, $(A,j,pc_{A,A}^*)$. 
Alternatively, 
$t$ and  $s$ can be taken respectively as
\begin{equation}\label{tasaction.e}
[A,A]\xrightarrow{\mathrm{id}\otimes1}[A,A]\otimes
[A,A]\longrightarrow[A\otimes A,A\otimes A]\xrightarrow{[p^*,p]}[A,A]
\end{equation}
and
\begin{equation}\label{sasaction.e}
[A,A]\xrightarrow{1\otimes\mathrm{id}}[A,A]\otimes[A,A]
\longrightarrow[A\otimes A,A\otimes A]\xrightarrow{[p^*,p]}[A,A] 
\end{equation}
where $\mathrm{id}:I\to[A,A]$ is the 1-cell corresponding to $1_A$
under the equivalence $\M(A,A)\simeq\M(I,[A,A])$.

\begin{obs}\label{o:[A,A]twostr}
  At this point we should remark that for a map pseudomonoid $A$,
  $[A,A]$ has two pseudomonoid structures. The one we have considered
  so far is the {\em composition}\/ pseudomonoid structure, but we
  also have the {\em convolution}\/ pseudomonoid structure. 

  If $(C,\varepsilon,\delta)$ is a pseudocomonoid in the braided Gray monoid \M\ 
  such that the 2-functor $C\otimes -$
  has a right biadjoint $[C,-]$, this is lax monoidal in the
  standard way. The unit constraint $I\to[C,I]$ corresponds under the
  closedness equivalence to the counit $\varepsilon:C\to I$ and the 1-cells
  $[C,X]\otimes[C,Y]\to[C,X\otimes Y]$ correspond 
$$
C\otimes[C,X]\otimes[C,Y]\xrightarrow{\delta\otimes1\otimes
  1}C^2\otimes [C,X]\otimes [C,Y]\xrightarrow{1\otimes c\otimes 
  1}(C\otimes[C,X])^2\xrightarrow{(\mathrm{ev}\otimes
  1)(1\otimes 1\otimes \mathrm{ev})}X\otimes Y.
$$
  In particular, for a pseudomonoid $A$, $[C,A]$ has a canonical {\em
    convolution}\/ 
  pseudomonoid structure. This structure corresponds to the usual
  convolution tensor product in $\M(C,A)$ given by $f*g=p(A\otimes
  g)(f\otimes C)\delta$ with unit $j\varepsilon$. As we saw in
  Observation \ref{o:1monoid}, for a map pseudomonoid $A$, the
  identity $1_A$ has a canonical structure of a monoid in the
  convolution monoidal category $\M(A,A)$. It follows that the
  corresponding 1-cell $\mathrm{id}:I\to [A,A]$ is a monoid in
  $\M(I,[A,A])$. 


\end{obs}


\begin{obs}\label{o:actions}
Let $B$ be a pseudomonoid in \M\ and consider $\M(I,B)$ and $\M(B,B)$
as monoidal categories with the convolution and the composition tensor
product respectively. We have monoidal functors
$L,R:\M(I,B)\to \M(B,B)$ given by $L(f)=p(f\otimes B)$ and
$R(f)=p(B\otimes f)$. The associativity constraint of $B$ induces
isomorphisms $L(f)R(g)\cong R(g)L(f)$, natural in $f$ and $g$. If $m$
and $n$ are monoids in $\M(I,B)$, then these isomorphisms form an
invertible distributive law between the monads $L(m)$ and $R(n)$. 

The moniodal functors $L,R$ are compatible with monoidal
pseudofunctors: if $F:\M\to \mathscr N$ is a monoidal pseudofunctor,
then there is a monoidal isomorphism of monoidal functors 
\begin{equation*}
\xymatrixrowsep{.5cm}
{\diagramcompileto{LRcompF}
\M(I,B)\ar[r]^{L,R}\ar[d]_{F_{I,B}}\ar@{}[ddr]|-\cong&
\M(B,B)\ar[dd]^{F_{B,B}}\\
{\mathscr{N}}(FI,FB)\ar[d]_{\simeq}&{}\\
{\mathscr{N}}(I,FB)\ar[r]^-{L,R}& 
{\mathscr{N}}(FB,FB)
\enddiagram
}
\end{equation*}
In particular, if $m$ is a monoid in $M(I,B)$, we have an isomorphisms 
$F(L(m))\cong L(Fm)$ and $F(R(m))\cong R(Fm)$ of monoids in $\M(B,B)$. 
\end{obs}
\begin{prop}\label{p:distrlawtands}
There exists an invertible distributive law between the monads $t$
and $s$, and hence between the monads $\theta$ and $\sigma$.
\end{prop}
\begin{proof}
  Apply Observation \ref{o:actions} above to the {\em convoluiton}\/
  pseudomonoid $B=[A,A]$ 
  and the monoid $m=n=\mathrm{id}:I\to[A,A]$, noting that
  $t=L(\mathrm{id})$ and $s=R(\mathrm{id})$. The 1-cell $\mathrm{id}$
  is a monoid with the structure given by Observation \ref{o:[A,A]twostr}. 
\end{proof}

If $t$ has an Eilenberg-Moore construction
$u:[A,A]^t\to[A,A]$ the monad $\hat\sigma$ is represented by some
$\hat s:[A,A]^t\to[A,A]^t$.
\begin{prop}\label{p:sandhatsopmonoidal}
The monads $s$ and $\hat s$ are opmonoidal monads.
\end{prop}
\begin{proof}
As we noted, $s$ is the monad $t$ for the pseudomonoid 
$(A,j,pc_{A,A})$. It can also be regarded as the corresponding monad
$t$ for the pseudomonoid $(A,j,p)$  in $\M^{\mathrm{rev}}$, 
and thus it is opmonoidal in
$\M^{\mathrm{rev}}$, and hence in \M.
The monad $\hat s$ is opmonoidal since
$[A,A]^t$ is an Eilenberg-Moore construction in $\mathbf{Opmon}(\M)$.
\end{proof}

Denote by $\hat\sigma$ the induced monad on $\M(A\otimes-,A)^\theta$
such that 
$$
\diagramcompileto{hatsigma}
\M(A\otimes-,A)^\theta\ar[r]^-{\hat\sigma}\ar[d]_{\upsilon}&
\M(A\otimes -,A)^\theta\ar[d]^\upsilon\\
\M(A\otimes-,A)\ar[r]_-\sigma&
\M(A\otimes -,A)
\enddiagram
$$
commutes. There exists an equivalence $(\M(A\otimes
-,A)^\theta)^{\hat\sigma}\simeq\M(A\otimes-,A)^{\sigma\theta}$. 

Suppose that there exists a pseudonatural transformation
$\tilde\sigma:\M(-,A)\to\M(-,A)$ such that
$\lambda\tilde\sigma\cong\hat\sigma\lambda$; since $\lambda$ is fully
faithful (see Proposition \ref{lambdaff.prop}), this is equivalent to
saying that for each $X$ the monad 
$\hat\sigma_X$ restricts to 
a monad on the replete image of $\lambda_X$ in $\M(A\otimes
X,A)^{\theta_X}$, and in this case
$\tilde\sigma=\lambda^*\hat\sigma\lambda$. Moreover, $\tilde\sigma$
carries the structure of a monad induced by the one of $\hat\sigma$.
Such a monad $\tilde\sigma $ clearly exists if the theorem of Hopf
modules holds for $A$, {\em i.e.}, if $\lambda$ is an equivalence. 

\begin{thm}\label{CPequiv.thm}
There exists an equivalence  in the 2-category
$[\M^{\mathrm{op}},\mathbf{Cat}]$ between $\M(-,A)^{\tilde\sigma}$ 
and $CP_\ell(-,A)$
whenever the monad $\tilde\sigma$ exists. Moreover, this equivalence
commutes with the corresponding forgetful pseudonatural transformations.
\end{thm}
\begin{proof}
We shall consider the restriction of $\hat\sigma_X$ to the replete image
of $\lambda_X$ instead of $\tilde\sigma_X$.
Take $f:X\to A$ and assume that $\lambda_X(f:X\to A)$ has a structure $\nu$ of
$\hat\sigma$-algebra. This means that the action $\nu$ is a 2-cell 
\begin{equation}\label{hatsigmaalg.eq}
\diagramcompileto{nu1}
A\otimes A\otimes X\ar[rr]^-{1\otimes c^*_{X,A}}&
{}\ddtwocell<\omit>{\nu}&
A\otimes X\otimes A\ar[d]^{1\otimes f\otimes 1}\\
A\otimes X\ar[u]^{p^*\otimes1}\ar[d]_{1\otimes f}&&
A\otimes A\otimes A\ar[d]^{p\otimes 1}\\
A\otimes A\ar[dr]_p&&
A\otimes A\ar[dl]^p\\
&A&
\enddiagram
\end{equation}
which is a morphism of $\theta_X$-algebras from
$\hat\sigma_X\lambda_X(f)$ to $\lambda_X(f)$. Furthermore, the pasting
$$
\diagramcompileto{nu2}
&A^2\otimes X \otimes A\ar[rr]^-{1\otimes c^*\otimes1}&
{}\ddtwocell<\omit>{\nu\otimes1}&
A\otimes X\otimes A^2\ar[d]^{1\otimes f\otimes A^2}\\
A^2\otimes X\ar[r]^-{1\otimes c^*}&
A\otimes X\otimes A\ar[d]_{1\otimes f\otimes1}\ar[u]^{p^*\otimes 1\otimes 1}&&
A^4\ar[d]^{p\otimes A^2}\\
A\otimes X\ar[u]^{p^*\otimes 1}\ar[d]_{1\otimes f}&
A^3\ar[dr]_{p\otimes 1}&&
A^3\ar[dl]^{p\otimes1}\\
A^2\ar[dr]_p&&
A^2\ar[dl]^p&\\
&A\uutwocell<\omit>{^\nu}
\enddiagram
$$
should be equal to the composition
$\sigma_X\sigma_X\lambda_X(f)\to\sigma_X\lambda_X(f)\xrightarrow{\nu}
\lambda_X(f)$ of the multiplication of $\sigma_X$ \eqref{sigmamult.eq}
and $\nu$, and the
composition
$\lambda_X(f)\to\sigma_X\lambda_X(f)\xrightarrow{\nu}\lambda_X(f)$ of
the unit of $\sigma$ \eqref{sigmaunit.eq} and $\nu$ is the identity.
The 2-cells \eqref{hatsigmaalg.eq} 
correspond, under pasting with $\phi^{-1}:p(A\otimes
p)\cong p(p\otimes A)$, to 2-cells $p(A\otimes (p(f\otimes
A)c^*_{X,A}))(p^*\otimes X)\Rightarrow p(A\otimes f)$, and then to
2-cells
$p(A\otimes (p(f\otimes A)c^*_{A,X}))\Rightarrow  p(A\otimes f)(p\otimes
A)\cong p(A\otimes p)(A\otimes A\otimes f)$. Since $\lambda_X$ is
fully faithful, and $\hat\sigma$ restricts to its replete image, it
follows that the 2-cells $\nu$ correspond to the  
2-cells $\gamma$ \eqref{gammacentre1.eq}. The axiom of associativity
for the action $\nu$ translates into the axiom \eqref{asscentre.eq} for
$\gamma$ and the axiom of unit for $\nu$ into the axiom
\eqref{unitcentre.eq} for $\gamma$. 
This shows that the composition of the forgetful functor
$V:CP_\ell(X,A)\to\M(X,A)$ with $\lambda_X$ 
factors as $G$ followed by $\hat U$, as depicted below. 
$$
\xymatrixrowsep{.5cm}
\diagramcompileto{CPalg}
CP_\ell(X,A)\ar@{-->}[rr]^-{H_X}\ar[d]_{V_X}\ar[drr]_(.3){G_X}&&
\M(X,A)^{\tilde\sigma_X}\ar[d]^{\tilde\lambda_X}\ar[dll]^(.2){\tilde U_X}\\
\M(X,A)\ar[dr]_{\lambda_X}&&
(\M(A\otimes X,A)^{\theta_X})^{\hat\sigma_X}\ar[dl]^{\hat U_X}\\
&\M(A\otimes X,A)^{\theta_X}
\enddiagram
$$
Moreover, $G_X$ factors through the image of $\tilde\lambda_X$, since
$\tilde U_XG_X$ factors through $\lambda_X$, and in fact $G_X$ is an
equivalence into the image of $\tilde\lambda_X$. 
Here $\tilde\lambda_X$ is the
functor induced on Eilenberg-Moore constructions by $\lambda_X$; in
particular, $\tilde\lambda_X$ is fully faithful since $\lambda_X$ is
fully faithful. Therefore we have an equivalence $H_X$ as in the diagram,
such that $\tilde\lambda_X H_X= G_X$. Hence, $\lambda_X\tilde U_X H_X=\hat
U_X\tilde \lambda_X H_X=\hat U_X G_X=\lambda_X V_X$, and $\tilde U_X
H_X=V_X$ . 
The equivalences $H_X$ are clearly pseudonatural in $X$. 

\end{proof}

\begin{cor}\label{centrecor1}
If the theorem of Hopf modules holds for a map pseudomonoid $A$ then
there exists an equivalence $CP_\ell(-,A)\simeq\M(A\otimes-,A)^{\sigma\theta}$.
\end{cor}
\begin{proof}
$\lambda_X$ is an equivalence and then the monad $\tilde\sigma$ exists and
$$
\M(-,A)^{\tilde\sigma}
\simeq({\M(A\otimes-,A)^{\theta}})^{\hat\sigma}\simeq\M(A\otimes
-,A)^{\sigma\theta}.
$$ 
\end{proof}

\begin{cor}\label{laxcentreEMconstr.cor}
Suppose that the theorem of Hopf modules holds for the map
pseudomonoid $A$ and that it has a Hopf module construction. Then the
lax centre of 
$A$ is the Eilenberg-Moore construction for the opmonoidal monad
$$
\tilde s:=\ell^*\hat s\ell=A\to A
$$
one of them existing if the other does. Moreover, 
$$
\tilde s \cong\big(A\xrightarrow{j\otimes1}A\otimes
A\xrightarrow{p^*\otimes1}A\otimes A\otimes A\xrightarrow{1\otimes
c^*_{A,A}}A\otimes A\otimes A\xrightarrow{p\otimes A}A\otimes
A\xrightarrow{p}A\big).
$$
\end{cor}
\begin{proof}
  The monad $\hat s$ exists and is opmonoidal since $t:[A,A]\to[A,A]$ has an
  Eilenberg-Moore construction in
  $\mathbf{Opmon}(\M)$.  Hence, $\tilde s$ has a canonical opmonoidal
  monad structure induced by the one of $\hat s$. The Theorem
  \ref{CPequiv.thm} implies that 
the lax centre of $A$ exists, that is, $CP(-,A)$ is birepresentable, if
and only if the monad $\tilde s$ has an Eilenberg-Moore construction. 

To obtain an expression for the 1-cell $\tilde s$ recall that, by
definition, $\M(-,\tilde s)$ is isomorphic to
$\lambda^*\hat\sigma\lambda$. It is easy to show that
\begin{align*}
\lambda^*_X\hat\sigma_X\lambda_X(f:X\to A)&
=p(p\otimes A)(A\otimes f\otimes A)(A\otimes c^*_{X,A})(p^*\otimes
X)(j\otimes X)\\
&\cong p(p\otimes A)(A\otimes c^*_{A,A})(p^*\otimes A)(j\otimes A)f;
\end{align*}
see Definitions \ref{lambda.def} and \ref{d:sigma}. It follows that
the expression for $\tilde s$ of the statement holds. 
\end{proof}

\begin{obs}\label{o:corOKAla}
The thesis of Corollary \ref{laxcentreEMconstr.cor} above holds under
the sole hypothesis of that $A$ be left autonomous. This is so because
every left autonomous map pseudomonoid has a Hopf module
construction (Theorem \ref{laut=hmc.thm}).
\end{obs}

\begin{obs}
Suppose that
the theorem of Hopf modules holds for $A$ and that
$A$ has a Hopf module construction ({\em e.g.}, $A$ is left autonomous). 
Then, when the lax centre $Z_\ell(A)$ of $A$ exists, it has a canonical
structure of a pseudomonoid such that the universal $Z_\ell(A)\to A$
is strong monoidal. 
\end{obs}

Now we concentrate in the case of autonomous map pseudomonoids. Let
$A$ be a such pseudomonoid. The internal hom $[A,A]$ is given by
$A^\circ\otimes A$, where $A^\circ$ is a right bidual of $A$. The
1-cell $\mathrm{id}$ is just the coevaluation $\n:I\to A^\circ\otimes
A$. The convolution 
\begin{cor}
Suppose $F:\M\to \mathscr N$ is a pseudofunctor between Gray monoids
with the following properties: $F$ preserves Eilenberg-Moore objects,
is braided and strong monoidal.
Then, $F$ preserves lax centres of left autonomous map pseudomonoids. 
\end{cor}
\begin{proof}
Let $A$ be a left autonomous map pseudomonoid in \M. 
  By Observation \ref{o:corOKAla}, the lax centre of $A$ is the
  Eilenberg-Moore construction for the opmonoidal monad $\tilde s:A\to
  A$, one existing if the other does. On the other hand, $FA$ is also
  a left autonomous map pseudomonoid by Proposition
  \ref{p:monsplaxpreslaut}. Therefore, it is enough to show
  that $F$ preserves the monad $\tilde s$, in the sense that $F\tilde
  s$ is isomophic to the corresponding monad $\tilde s$ for $FA$. 
  
  Since $\tilde s$ is the
  lifting of the monad $s$ on $A^\circ \otimes A$ to the
  Eilenberg-Moore construcion $(A^\circ\otimes p)(\n\otimes A):A\to
  A^\circ\otimes A$ of the monad $t$ (see Theorem \ref{laut=hmc.thm}), 
  it suffices to prove that $F$
  preserves the monads $t$ and $s$. We only work with $t$, the proof
  for the monad $s$ being completely analogous. Now, we know from the proof of
  Proposition \ref{p:distrlawtands} that $t=L(\n)$ and
  $s=R(\n)$, where $L,R:\M(I,A^\circ\otimes A)\to\M(A^\circ\otimes
  A,A^\circ\otimes A)$ are the functors defined in Observation
  \ref{o:actions}. Therefore, $Ft=F(L(\n))\cong
  L(I\xrightarrow{\cong}FI\xrightarrow{F\n_A}F(A^\circ\otimes A))\cong
  L(\n_{FA})$, which is the monad $t$ corresponding to the
  pseudomonoid $FA$.

\end{proof}

\begin{thm}\label{lax=centre.th}
For a (left and right) autonomous map pseudomonoid  the centre
equals the lax centre, either existing if the other does.  
\end{thm}
\begin{proof}
Consider the commutative diagram 
$$
\diagramcompileto{centre21}
(\M(A\otimes X,A)^{\theta_X})^{\hat\sigma_X}\ar[d]\ar[r]&
\M(A\otimes X,A)^{\theta_X}\ar[r]^-{\hat\sigma}\ar[d]&
\M(A\otimes X,A)^{\theta_X}\ar[d]^{\upsilon_X}\\
\M(A\otimes X,A)^{\sigma_X}\ar[r]&
\M(A\otimes X,A)\ar[r]_-\sigma&
\M(A\otimes X,A)
\enddiagram
$$
In Theorem \ref{CPequiv.thm} we proved that any lax centre piece arises as
\begin{equation}\label{centrenu.eq}
\diagramcompileto{centre2nu}
A\otimes A\otimes X\ar@{=}[r]\ar[dr]_{p\otimes 1}\rtwocell<\omit>{<2>}&
A\otimes A\otimes X\ar[rr]^-{1\otimes c_{A,X}}&
{}\ddtwocell<\omit>{\nu}&
A\otimes X\otimes A\ar[d]^{h\otimes 1}\\
&A\otimes X\ar[u]_{p^*\otimes1}\ar[dr]_{h}&&A\otimes A\ar[dl]^p\\
&&A&
\enddiagram
\end{equation}
for some $\hat\sigma_X$-algebra $\nu:\hat\sigma_X(h)\to h$, so we
have to prove that \eqref{centrenu.eq} is invertible. Consider the
canonical split coequalizer
$\hat\sigma_X^2(h)\rightrightarrows\hat\sigma_X(h)\twoheadrightarrow h$
in $\M(A\otimes X,A)^{\theta_X}$, and its image $\nu:\sigma_X(h)\to h$
in $\M(A\otimes X,A)$. The arrow $\nu$ is a morphism of
$\sigma_X$-algebras. This implies that the lower rectangle in the
diagram below commutes. 
$$
\xymatrixcolsep{1.8cm}
\diagramcompileto{centre2h}
p(\sigma_X(h)\otimes A)(A\otimes c_{A,X})\ar[d]|-{p(\sigma_X(h)\otimes
  A)(A\otimes c_{A,X})(\eta\otimes X)}
\ar[r]&p(h\otimes A)(A\otimes c_{A,X})\ar[d]|--{p(h\otimes A)(A\otimes
  c_{A,X})(\eta\otimes X)}\\
p(\sigma_X(h)\otimes A)(A\otimes c_{A,X})(p^*p\otimes X) \ar@{=}[d]&
p(h\otimes A)(A\otimes c_{A,X})(p^*p\otimes X)\ar@{=}[d]\\
\sigma_X^2(h)(p\otimes X)\ar[d]|-{(\mu_X)_h (p\otimes X)}\ar[r]&
\sigma_X(h)(p\otimes X)\ar[d]^{\nu (p\otimes X)}\\
\sigma_X(h)(p\otimes X)\ar[r]&
h(p\otimes X)
\enddiagram
$$
The upper rectangle commutes by naturality of composition. Here $\eta$
denotes the unit of the adjunction $p\dashv p^*$ and $\mu$ the
multiplication of the monad $\sigma$. Observe
that the rows are coequalizers and the right-hand column is just
\eqref{centrenu.eq}. Then, to show that this last arrow is invertible
it suffices to show that the left-hand side column, which is the
pasting of $\eta$ with the multiplication of $\sigma$
\eqref{sigmamult.eq}, is so. 
But this 2-cell is invertible because $A$
is right autonomous and by
the dual of Theorem \ref{maintheorem.th}.\ref{maintheoremgamma.it} the
 2-cell below is invertible. This completes the proof.
$$
\diagramcompileto{centre2prevgamma}
A^2\ar@{=}[rr]\ar[dr]_p\rrtwocell<\omit>{<3>}&&
A^2\dtwocell<\omit>{'\cong}\ar[r]^-{p^*\otimes 1}&
A^3\ar[dr]^{1\otimes p}&\\
&A\ar[ur]_{p^*}\ar[r]_-{p^*}&
A^2\ar[ur]|-{1\otimes p^*}\ar@{=}[rr]\rrtwocell<\omit>{<-3>}&&
A^2
\enddiagram
$$
$$
\parallel
$$
$$
\diagramcompileto{centre2gamma}
&A^3\ar[r]^-{1\otimes p}\ar[dr]|-{p\otimes 1}
\drrtwocell<\omit>{'\cong}&
A^2\ar[dr]^p\ar@{=}[rr]\rrtwocell<\omit>{<3>}&&
A^2\\
A^2\ar[ur]^{p^*\otimes 1}\ar@{=}[rr]
\rrtwocell<\omit>{<-3>} &&
A^2\ar[r]_-p&
A\ar[ur]_{p^*}
\enddiagram
$$
\end{proof}

\begin{cor}\label{lc=c.cor}
Any autonomous map pseudomonoid in a braided monoidal bicategory with
Eilenberg-Moore objects has both a centre and a lax centre, and the two 
coincide.
\end{cor}


\section{\V-modules}\label{VMod.sec}

In this section we interpret the results 
of the previous section in the particular context of the 
bicategory of \V-modules. 

\subsection{Review of the bicategory of \V-modules}\label{ss:reviewVmod}

Throughout this section \V\ will be a complete and cocomplete closed
symmetric monoidal category. There is a 
bicategory $\V\text{-}\mathbf{Mod}$  whose objects are the small \V-categories
and hom-categories $\V\text{-}\mathbf{Mod}(\mathscr A,\mathscr
B)=[\mathscr A^{\mathrm{op}}\otimes \mathscr B,\V]_0$, the category of
\V-functors  from the tensor product of
the \V-categories $\mathscr A^{\mathrm{op}}$ and $\mathscr B$ to \V,
and \V-natural transformations between them. Objects of this category
are called \V-modules and arrows morphisms of \V-modules. The composition
of two \V-modules $M:\mathscr A\to\mathscr B$ and $N:\mathscr
B\to\mathscr C$ is given by $(NM)(a,c)=\int^x N(x,c)\otimes M(a,x)$.
This coend exists because \V\ is cocomplete. The identity module
$1_{\mathscr A}$ is given by $1_{\mathscr A}(a,a')=\mathscr
A(a,a')$. A \V-module $M:\A\to\B$ can also be thought as a
$(\mathrm{ob}\A\times\mathrm{ob}\B)$-indexed family of objects
$M(a,b)$ of \V\ with compatible actions of $\A$ on the right and of
$\B$ on the left; this is, actions $\B(b,b')\otimes M(a,b)\to M(a,b')$
and $M(a,b)\otimes\A(a',a)\to M(a',b)$ subject to compatibility
conditions.

Our convention is that a 
\V-module from \A\ to \B\ as a \V-functor
$\A^{\mathrm{op}}\otimes\B\to\V$, but 
some authors prefer to use \V-functors
$\A\otimes\B^{\mathrm{op}}\to\V$. A different approach was taken in
\cite{variation}, where the bicategory of \V-matrices was used to
define \V-modules for a {\em bicategory}\/ \V. 

There is a pseudofunctor
$(-)_*:\V\text{-}\mathbf{Cat}^{\mathrm{co}}\to\V\text{-}\mathbf{Mod}$
which is the identity on objects and on hom-categories $[\mathscr
A,\mathscr B]_0^{\mathrm{op}}\to[\mathscr A^{\mathrm{op}}\otimes
\mathscr B,\V]_0$  sends a \V-functor $F$ to the \V-functor
$F_*(a,b)=\mathscr B(F(a),b)$. Moreover, the \V-module $F_*$ has right
adjoint $F^*$ given by $F^*(b,a)=\mathscr B(b,F(a))$.
The pseudofunctor $(-)_*$ is easily shown to
be strong monoidal and symmetry-preserving.

The tensor product of \V-categories induces a structure of a
monoidal bicategory on $\V\text{-}\mathbf{Mod}$, which on
hom-categories $[\mathscr A^{\mathrm{op}}\otimes\mathscr B,\V]_0\otimes
[{\mathscr A'}^{\mathrm{op}}\otimes\mathscr B',\V]_0\to[{\mathscr
A'}\otimes\mathscr A^{\mathrm{op}}\otimes\mathscr B\otimes \mathscr B',\V]_0$
is given by point-wise tensor product in \V. Moreover, the usual
symmtery of $\V\text-\mathbf{Cat}$ together with the symmetry of \V\
induce a structure of symmetric monoidal bicategory on
$\V\text{-}\mathbf{Mod}$, or rather, induce a symmetry in the sense of
\cite{monbicat} in any Gray
monoid monoidally equivalent to $\V\text{-}\mathbf{Mod}$.
The natural isomorphisms
$\V\text{-}\mathbf{Mod}(\mathscr B,\mathscr
A^{\mathrm{op}}\otimes\mathscr C)\cong[\B^{\mathrm{op}}\otimes
\A^{\mathrm{op}}\otimes\mathscr 
C,\V]_0\cong\V\text{-}\mathbf{Mod}(\mathscr A\otimes\mathscr
B,\mathscr C)$ show that our monoidal bicategory is also autonomous
with (right and left) biduals given by the opposite \V-category. The
coevaluation $\n:I\to\mathscr A^{\mathrm{op}}\otimes\mathscr A$ and
coevaluation $\e:\mathscr A\otimes\mathscr A^{\mathrm{op}}\to I$
modules are given respectively by $\n(a,b)=\mathscr A(a,b)$ and
$\e(a,b)=\mathscr A(b,a)$, and the bidual of a \V-module $M:\mathscr
A\to\mathscr B$ can be taken as $M^{\circ}(b,a)=M(a,b)$. (Note that
$\e$ and $\n$ do not induce the isomorphisms above, but pseudonatural
equivalences isomorphic to these). 

One of the many pleasant properties of $\V\text{-}\mathbf{Mod}$ is
that it has right liftings. If $M:\B\to \C$ and $N:\A\to \C$ are
\V-modules, a right lifting of $N$ through $M$ is given by the formula
${}^M\! N(a,b)=\int_{c\in\C}[M(b,c),N(a,c)]$. 
As explained in Section \ref{s:thelaxcentre}, the existence of right
liftings endows each hom-category $\V\text{-}\mathbf{Mod}(I,\A)$ with
a canonical structure of a $\V\text{-}\mathbf{Mod}(I,I)$-category,
where $I$ is the trivial \V-category. Henceforth, each
$\V\text{-}\mathbf{Mod}(I,\A)$ is canonically a \V-category  via the
monoidal isomorphism $\V\text{-}\mathbf{Mod}(I,I)\cong\V$. The
\V-enriched hom $\V\text{-}\mathbf{Mod}(I,\A)(M,N)$ is given by
${}^M \! N$,
or explicitly by the object $\int_{a\in\A}[M(a),N(a)]$. 
This is exactly the usual \V-category structure of $[\A,\V]$. 
In fact, each hom-category $\V\text{-}\mathbf{Mod}(\A,\B)$ is
canonically a \V-category, in a way such that the equivalence
$\V\text{-}\mathbf{Mod}(\A,\B)\simeq
\V\text{-}\mathbf{Mod}(I,\A^{\mathrm{op}}\otimes\B)$ is a \V-functor. 

Another feature of $\V\text- \mathbf{Mod}$ we will need is the
existence of Kleisli and Eilenberg-Moore constructions for monads. The
existence of the former was shown in
\cite{Street:CauchyChar}. If $(M,\eta,\mu)$ is a monad in
$\V\text{-}\mathbf{Mod}$ on \A, $\mathsf{Kl}(M)$ has the same objects
as \A\ and homs $\mathsf{Kl}(M)(a,b)=M(a,b)$. Composition is given by
$$
M(b,c)\otimes M(a,b)\to\int^{b\in\A}M(b,c)\otimes
M(a,b)\xrightarrow{\mu_{a,c}}M(a,c) 
$$
and the units by
$I\xrightarrow{\mathrm{id}}\A(a,a)\xrightarrow{\eta_{a,a}}M(a,a)$.
One can verify that the \V-module $K_*$ induced by the 
\V-functor $K:\A\to\mathsf{Kl}(M)$ given by the
identity on objects and by $\eta_{a,b}:\A(a,b)\to M(a,b)$ on homs has
the universal property of the Kleisli construction. With regard to
Eilenberg-Moore constructions, $K^*$ induces, for each $\mathscr X$,  a
functor  $(K^*\circ-):\V\text-
\mathbf{Mod}(\X,\mathsf{Kl}(M))\to\V\text- \mathbf{Mod}(\X,\A)$. This
functor is isomorphic to the one 
sending $M:\mathscr X^{\mathrm{op}}\otimes \mathsf{Kl}(M)\to\V$ to the
$M(\mathscr X^{\mathrm{op}}\otimes K):\mathscr
X^{\mathrm{op}}\otimes \A\to \V$. Therefore, $(K^*\circ-)$ is
conservative because $K$ is the identity on objects: if
$\alpha:M\Rightarrow N$ is a \V-module morphism such that
$\alpha(\mathscr X^{\mathrm{op}}\otimes K)$ is an isomorphism, we have
$\alpha_{x,a}=\alpha_{x,K(a)}$ is an isomorphism for all $x\in
\mathscr X$ and $a\in\A$. It follows that $(K^*\circ- )$ is monadic,
as it is clearly cocontinuous and has a left adjoint. This being true
for all $\mathscr X$, we deduce $K^*$ is monadic in
$\V\text-\mathbf{Mod}$, and hence it is an Eilenberg-Moore
construction for $M$.

\subsection{Left autonomous pseudomonoids and left autonomous \V-categories}

A pseudomonoid in $\V\text{-}\mathbf{Mod}$ is a promonodial
\V-category \cite{closedcatfun}. The pseudomonoid
structure amounts to a multiplication and a unit
\V-functors $P:\mathscr A^{\mathrm{op}}\otimes\mathscr
A^{\mathrm{op}}\otimes\mathscr A\to\V$ and $J:\mathscr A\to \V$
together with associativity and unit \V-natural constraints satisfying
axioms. Any monoidal $\V$-category can be thought of
as a promonoidal \V-category, in fact a map pseudomonoid, 
by using the monoidal pseudofunctor
$(-)_*:\V\text{-}\mathbf{Cat}^{\mathrm{co}}\to\V\text{-}\mathbf{Mod}$;
explicitly, if $\mathscr A$ is a monoidal \V-category, then the
induced promonoidal structure is given by $P(a,b;c)=\mathscr
A(b\otimes a,c)$ and $J(a)=\mathscr A(I,a)$.

Next we show how the results on Hopf modules specialise to the
bicategory of \V-modules, and give explicit descriptions of the main
constructions. Although these descriptions carry over to arbitrary left
autonomous map pseudomonoids, here we will concentrate on the simpler
case of the left
autonomous monoidal \V-categories \A.
The opmonoidal 
monad $T:\mathscr A^{\mathrm{op}}\otimes
\mathscr A\to\mathscr A^{\mathrm{op}}\otimes
\mathscr A$ defined in Section \ref{inthopf.s} is given as a \V-module
by
$
T(a,b;c,d)=\int^x\mathscr A(b\otimes x,d)\otimes\mathscr A(c,a\otimes x).
$
The multiplication is has components 
\begin{multline*}
T^2(a,b;c,d)=\int^{u,v}\int^x\mathscr A(v\otimes x,d)\otimes \mathscr
A(c,u\otimes x)\otimes\int^y\mathscr A(b\otimes y, v)\otimes\mathscr
A(u,a\otimes y)\\
\cong\int^{x,y}\mathscr A((b\otimes y)\otimes x,d)\otimes\mathscr
A(c,(a\otimes y)\otimes x)\\
\cong
\int^{x,y}\mathscr A(b\otimes (y\otimes x),d)\otimes\mathscr
A(c,a\otimes (y\otimes x))\longrightarrow T(a,b;c,d)\\
\end{multline*}
where the last arrow is induced by the obvious arrows
$\mathscr A(b\otimes (y\otimes x),d)\otimes\mathscr
A(c,a\otimes (y\otimes x))\to\int^x\mathscr A(b\otimes
x,d)\otimes\mathscr A(c,a\otimes x)$. 
The unit has components
$$
(\mathscr A^{\mathrm{op}}\otimes \mathscr A)(a,b;c,d)=\mathscr
A(b,d)\otimes \mathscr
A(c,a)\xrightarrow{\eta}\int^x\mathscr A(b\otimes
x,d)\otimes\mathscr A(c,a\otimes x),
$$
the component corresponding to $I\in\mathrm{ob}\mathscr A$. 


The existence of Eilenberg-Moore constructions in $ \V\text- \mathbf{Mod}$ 
implies the following. 

\begin{prop}\label{p:hopfmodinvmod}
Any map pseudomonoid in $\V\textrm{-}\mathbf{Mod}$ has a Hopf module
construction. 
\end{prop}

Following the remaks at the end of the previous section, one can give
an explicit description of the Hopf module construction for a map
pseudomonoid \A. 
The \V-category $(\mathscr A^{\mathrm{op}}\otimes\mathscr
A)^T=(\mathscr A^{\mathrm{op}}\otimes\mathscr A)_T$ has the same
objects as $\mathscr A^{\mathrm{op}}\otimes\mathscr A$, homs
$(\mathscr A^{\mathrm{op}}\otimes\mathscr A)(a,b;c,d)=T(a,b;c,d)$ and
composition and identities induced by the multiplication and unit of
$T$. The unit of the monad $T$  defines a \V-functor
$\eta:\mathscr A^{\mathrm{op}}\otimes \mathscr A\to(\mathscr
A^{\mathrm{op}}\otimes \mathscr
A)^T$; the Kleisli construction for $T$ is just the module  $\eta_*$
and the Eilenberg-Moore construction is $\eta^*$. The module
$L:\mathscr A\to(\mathscr
A^{\mathrm{op}}\otimes \mathscr
A)^T$ in \eqref{lamrep.eq} is then
\begin{equation}\label{Lexplicit.eq}
L=\big(\mathscr A \xrightarrow{(J^*)^\mathrm{\circ}\otimes \mathscr
A}\mathscr A^{\mathrm{op}}\otimes\mathscr
A\xrightarrow{\eta_*}(\A^{\mathrm{op}}\otimes\mathscr A)^T \big)
\end{equation}

When the promonoidal structure is induced by a monoidal structure on
\A, {\em i.e.}, $P(a,b;c)=\A(b\otimes a,c)$ and $J(a)=\A(I,a)$, we can
compute $L$ explicitly. Firstly note that for any \V-functor
$F:\mathscr B\to\mathscr C$ there exists a canonical isomorphism of \V-modules
$(F^*)^\mathrm{\circ}\cong (F^\mathrm{op})_*:\mathscr
B^\mathrm{op}\to\mathscr C^\mathrm{op}$, where
$F^\mathrm{op}:\mathscr B^\mathrm{op}\to\mathscr C^\mathrm{op}$ is the 
usual opposite functor.
Then 
$$
L\cong \eta_*((J^\mathrm{op})_*\otimes \mathscr
A)\cong(\eta(J^\mathrm{op}\otimes \mathscr A))_* .
$$ 
In components, 
$$
L(a;b,c)\cong(\mathscr A^\mathrm{op}\otimes\mathscr
A)^T(\eta(I,a),(b,c)) =T(I,a;b,c)\cong\mathscr A(a\otimes b,c)
$$
with right $\mathscr A$-action and left $(\mathscr
A^\mathrm{op}\otimes \mathscr A)^T$-action. 
The latter is given by the composition of
$(\A^{\mathrm{op}}\otimes\A)^T$, while the \A-action can be shown to
be given as 
$$
\A(a\otimes b,c)\otimes \A(a',a)\xrightarrow{1\otimes (-\otimes
  b)}\A(a\otimes b,c)\otimes\A(a'\otimes b,a\otimes
b)\xrightarrow{\mathrm{comp}}\A(a'\otimes b,c).
$$
The fact that $L$ is a
faithful \V-module (Proposition \ref{lffstrongmon.prop})  means exactly
that the \V-functor $\eta(J^{\mathrm{op}}\otimes\mathscr A)$ is fully
faithful. This can be also verified directly, for the effect of this
\V-functor on homs is 
$$
\mathscr A(b,d)\xrightarrow{1\otimes 1_I}\mathscr
A(b,d)\otimes\mathscr A(I,I)\xrightarrow{\eta}\int^x\mathscr
A(b\otimes x,d)\otimes \mathscr A(I,I\otimes x)\cong \mathscr A(b,d)
$$
sending an arrow $f$ to $\big(b\xrightarrow{\cong}b\otimes
I\xrightarrow{f\otimes 1_I}d\otimes I\xrightarrow{\cong}d\big)$.

Consider for a moment a general 
promonoidal \V-category $\mathscr A$. It is a left autonomous
pseudomonoid in $\V\text{-}\mathbf{Mod}$ when there exists a \V-module
$D:\mathscr A^\mathrm{op}\to\mathscr A$ such that $(P\otimes\mathscr
A)(\mathscr A\otimes D\otimes \mathscr A)(\mathscr A\otimes N)$ is
right adjoint to the multiplication $P:\mathscr A\otimes \mathscr A\to
\mathscr A$. The former \V-module is given in components by 
$$
(P\otimes\mathscr
A)(\mathscr A\otimes D\otimes \mathscr A)(\mathscr A\otimes
N)(a;b,c)\cong
P(D\otimes\mathscr A)(a,c;d)\\
\cong\int^xP(a,x;b)\otimes D(c,x).
$$
When $\mathscr A$ is a monoidal \V-category and the $D$ is induced by
a \V-functor, denoted by $(-)^\vee:\mathscr A^\mathrm{op}\to\mathscr A$,
then the expression above reduces to  
$$
(P\otimes\mathscr
A)(\mathscr A\otimes D\otimes \mathscr A)(\mathscr A\otimes
N)(a;b,c)\cong
\mathscr A(c^\vee\otimes a,b)
$$
and we obtain a \V-natural isomorphism $\mathscr A(c^\vee\otimes
a,b)\cong P^*(a;b,c)=\mathscr A(a,c\otimes b)$. We see, thus, that a
monoidal \V-category is a left autonomous pseudomonoid in
$\V\text{-}\mathbf{Mod}$ with representable dualization if and only if
it is left autonomous in the classical sense. This was first shown in
\cite{dualizations}. 

By the Theorem \ref{maintheorem.th} and the definition of $L$ in
\eqref{lamrep.eq} we have
\begin{prop}
Every promonoidal \V-category $\mathscr A$ for which $P$ and $J$ have
right adjoints has a structure (and a fortiori
unique up to isomorphism) of left
autonomous pseudomonoid in $\V\text{-}\mathbf{Mod}$ if and only if the
\V-module $L$ \eqref{Lexplicit.eq} is an equivalence. In pariticular,
this is true for any monoidal \V-category. 
\end{prop}

\subsection{Lax centres in $\V\textrm{-}\mathbf{Mod}$}

In this section we study the centre and lax centre of pseudomonoid
in the monoidal bicategory of \V-modules by means of the theory developed
in Section \ref{centre.s}. In the way, we compare our work with
\cite{Day:LaxCentre,Day:CentresMonCatsFun}. 

First we consider lax centres of arbitrary pseudomonoids. We shall
show that the results in Section \ref{s:thelaxcentre} apply to
$\V\text{-}\mathbf{Mod}$. To this aim, we have to verify all the
hypothesis required in that section. 

We already saw in Section \ref{ss:reviewVmod} that liftings of arrows out of $I$
through arrows out of $I$ exist.
In order to show $\V\text{-}\mathbf{Mod}$ satisfies the other two
hypothesis required in Section \ref{s:thelaxcentre}, it is enough to
prove that the arrow 
\eqref{eq:M(IA)faith} is an isomorphism for \M\ the bicategory of
\V-modules. In this case \eqref{eq:M(IA)faith} becomes
\begin{equation}\label{eq:M(IA)faithVmod}
[\A^{\mathrm{op}}\otimes \B,\V](M,N)
\to
[[\A,\V],[\B,\V]]((M\circ-),(N \circ-)), 
\end{equation}
where $(M\circ -)$ is the \V-functor given by composition with $M$. To
show that \eqref{eq:M(IA)faithVmod} is an isomorphism, recall that the
\V-functor 
  \begin{equation}\label{bcect4.51}
  [\A^{\mathrm{op}}\otimes \B,\V]\cong[\B,[\A^{\mathrm{op}},\V]]
    \to\mathrm{Cocts}[[\B^{\mathrm{op}},\V],[\A^{\mathrm{op}},\V]]
  \end{equation}
  into the sub-\V-category of
  cocontinuous \V-functors is an equivalence by \cite[Theorem
  4.51]{Kelly:BCECT}. This \V-functor sends 
  $R:\C^{\mathrm{op}}\otimes\C\to\V$ to the left extension of
  the corresponding $R':\C\to[\C^{\mathrm{op}},\V]$ along the Yoneda
  embedding $\mathsf{y}:\C\to[\C^{\mathrm{op}},\V]$,
  $\operatorname{Lan}_{\mathsf{y}}R'$, which is exactly $(R\circ -)$. 


Theorem \ref{t:MZsimeqZM} gives:

\begin{cor}\label{c:ZlinVmod}
Suppose the lax centre of the promonoidal \V-category \A\ exists. 
Then there exits an equivalence of
\V-categories $[Z_\ell \A,\V]\simeq Z_\ell[\A,\V]$, where on the left
hand side appears the lax centre in $\V\text{-}\mathbf{Mod}$ and on
the right hand side the lax centre in $\V\text{-}\mathbf{Cat}$. 
The composition of this equivalence with the forgetful \V-functor
$Z_\ell[\A,\V]\to[\A,\V]$ is canonically isomorphic to the \V-functor
given by composing with the universal \V-module $Z_\ell\A\to \A$. 
If the centre of \A, rather than the lax centre, exists, then the
above holds substituting lax centres by centres throughout.
\end{cor}

Now we turn our attention to autonomous pseudomonoids. 

The existence of Eilenberg-Moore constructions in
$\V\text{-}\mathbf{Mod}$ together with Corollary
\ref{laxcentreEMconstr.cor} and Theorem \ref{lax=centre.th} implies

\begin{prop}\label{p:laut->Zl}
Any left autonomous map pseudomonoid in $\V\text{-}\mathbf{Mod}$ has a lax
centre. Moreover, if the pseudomonoid is also right autonomous then
the lax centre is the centre.
\end{prop}

We shall describe the lax centre explicitly. In order to simplify the
description, we will suppose $\mathscr A$ is a left autonomous monoidal
\V-category, and not merely a promonoidal one. 
However, all the following description carries over to the map
pseudomonoid case with little modification.

By Corollary \ref{laxcentreEMconstr.cor}, the lax centre of
$\mathscr A$ in $\V\text{-}\mathbf{Mod}$ is the Eilenberg-Moore
construction for the monad $\tilde S$ given by
\begin{equation}\label{eq:tildeSmonoidalcat}
\mathscr A\xrightarrow{J\otimes1}\mathscr A\otimes\mathscr
A\xrightarrow{P^*\otimes1}\mathscr A\otimes\mathscr A\otimes\mathscr
A\xrightarrow{1\otimes \mathrm{sw}_*}\mathscr A\otimes\mathscr A\otimes\mathscr
A \xrightarrow{P\otimes1}\mathscr A\otimes\mathscr
A\xrightarrow{P}\mathscr A
\end{equation}
where $\mathrm{sw}$ denotes the usual symmetry in $\V\text{-}\mathbf{Cat}$ that
switches the two factors. 
Explicitly, 
\begin{equation*}
\tilde S(a;b)\cong\int^{x,y}\mathscr A(y\otimes(a\otimes
x),b)\otimes\mathscr A(I,y\otimes
x)\cong\int^y\mathscr A (y\otimes(a\otimes y^\vee),b), 
\end{equation*} 
where $y^\vee$ denotes the left dual of $y$ in $\mathscr A$.
The multiplication of this monad is given by
\begin{multline*}
\tilde S^2(a;b)\cong\int^{u,y,z}\mathscr A\big(y\otimes(u\otimes
y^\vee),b\big)\otimes\mathscr A\big(z\otimes(a\otimes z^\vee),u\big)\\
\cong\int^{y,z}\mathscr A\big(y\otimes(z\otimes(a\otimes z^\vee))\otimes
y^\vee,b\big) \cong\int^{y,z}\mathscr A\big((y\otimes
z)\otimes(a\otimes(y\otimes z)^\vee),b\big)\to\\
\longrightarrow\int^x\mathscr A\big(x\otimes(a\otimes x^\vee),b)\cong
\tilde S(a;b)
\end{multline*}
where the last arrow is induced by the components 
$\zeta_{y\otimes z}^{a,b}:\mathscr A\big((y\otimes
z)\otimes(a\otimes(y\otimes z)^\vee),b\big)\to\int^x\mathscr
A\big(x\otimes(a\otimes x^\vee),b\big)$ of the universal dinatural
transformation defining the latter coend. The unit of $\tilde S$ is
given by components 
$$
\mathscr A(a,b) \xrightarrow{\zeta_I^{a,b}}\int^{x}\mathscr
A\big(x\otimes(a\otimes x^\vee),b)
$$
of the same dinatural transformation corresponding to $x=I$. Now we
have all the ingredients to describe the lax centre $Z_\ell(\mathscr
A)$, that is, a Kleisli construction for $\tilde S$.
It has the same objects as $\mathscr A$, enriched homs
$Z_\ell(\mathscr A)(a,b)=\tilde S(a,b)$, composition
given by the multiplication and unit given by 
$$
I\to\mathscr
A(a,a)\xrightarrow{\zeta_{I}^{a,a}}\tilde S(a,a),
$$ 
where the first arrow
denotes the identity of $a$ in $\mathscr A$. The arrows
$\zeta_I^{a,b}:\mathscr A(a,b)\to \tilde S(a,b)$ define a \V-functor,
which we also call $\zeta$,
and the universal
$Z_\ell(\mathscr A)\to\mathscr A$ is none but $\zeta^*$.

\begin{obs}
  The monad $\tilde S$ is closely related to the monad $\check{M}$ in
  \cite[Section 5]{Day:CentresMonCatsFun}. There the authors show that for a general small
  promonoidal \V-category \C\ there exists a monad $\check{M}$ on \C\
  in $\V\text-\mathbf{Mod}$ with the following property. 
  Whenever $[\C,\V]$ has a small dense
  sub-\V-category of objects with left duals (it is {\em right-dual
    controlled}, in the terminology of
  \cite{Day:CentresMonCatsFun}), 
  the forgetful \V-functor $Z_\ell[\C,\V]\to[\C,\V]$ is a
  (bicategorical) Eilenberg-Moore construction for the monad $M$ on
  $[\C,\V]$ in $\V\text-\mathbf{Cat}$ given 
  by composition with $\check M$. 
  The monad $\check{M}$
  is given by 
  $$\check{M}(a,b)=\int^{x,y} P(P\otimes \C)(y,a,x,b)\otimes
  x^{\wedge}(y),$$ where $x^\wedge$ is the internal hom $\llbracket
  \C(x,-),J\rrbracket\in[\C,\V]$ ($J:I\to \C$ is the  unit of the
  promonoidal structure).
  
  When $\C$ is equipped with a left dualization
  $D:\C^{\mathrm{op}}\to\C$, each \V-module $I\to \C$ with right
  adjoint in $\V\text-\mathbf{Mod}$ has a left dual in the monoidal
  \V-category $\V\text-\mathbf{Mod}(I,\C)=[\C,\V]$. This was first
  shown in \cite{dualizations}. Explicitly, a left dual for $M:I\to \C$ is
  given by $D(M^*)^\circ\cong (\C\otimes M^*)P^*J$. In particular,
  $\C(x,-)$, which is the 
  \V-module induced by the \V-functor $I\to\C$ constant on $x$,  has
  left dual $(\C\otimes \C(-,x))P^*J\cong
  \llbracket\C(x,-),J\rrbracket$. 
  It follows that $[\C,\V]$ has a small dense sub-\V-category with
  left duals, and the results of \cite{Day:CentresMonCatsFun}
  mentioned above apply.
 
  In this situation, if we assume $J$ is a map,
  so that $\tilde S$ exists, 
  we claim that the monads $\check M$ and $\tilde S$ are
  isomorphic, or more precisely, that both are isomorphic as monoids in
  the monoidal \V-category
  $\V\text-\mathbf{Mod}(\C,\C)=[\C^{\mathrm{op}}\otimes\C,\V]$. To
  show this, it is enough to prove that the monads $(\check M\circ-)$
  and $(\tilde S\circ -)$ on $\V\text-\mathbf{Mod}(I,\C)=[\C,\V]$
  given by composition with $\check M$ and $\tilde S$ respectively are
  isomorphic. For, the \V-functor 
  $$
  [\C^{\mathrm{op}}\otimes \C,\V]\cong[\C,[\C^{\mathrm{op}},\V]]
    \to\mathrm{Cocts}[[\C^{\mathrm{op}},\V],[\C^{\mathrm{op}},\V]]
  $$
  into the sub-\V-category of
  cocontinuous \V-functors in \eqref{bcect4.51} is an equivalence by \cite[Theorem
  4.51]{Kelly:BCECT}. This \V-functor is monoidal and sends 
  $R:\C^{\mathrm{op}}\otimes\C\to\V$ to $(R\circ -)$.  
  
  Now, the monad $(\tilde S\circ-)$ is $\V\text-\mathbf{Mod}(I,\tilde
  S)$, and then it has the forgetful \V-functor
  $Z_\ell[\C,\V]\to[\C,\V]$ as an (bicategorical) Eilenberg-Moore
  construction by Corollary \ref{c:ZlinVmod} and Proposition
  \ref{p:laut->Zl}. Then, $(\tilde S\circ-)$ and $M=(\check M\circ-)$
  have the same Eilenberg-Moore construction in $\V\text-\mathbf{Cat}$
  and it follows that both monads are isomorphic as required. 

  More explicitly, by the description of the left dual of a \V-module
  $I\to \C$, we have
  $$
  \llbracket\C(x,-),J\rrbracket(y)\cong \int^{u,v}
  \C(u,y)\otimes \C(v,x)\otimes (P^*J)(u,v)\cong (P^*J)(y,x)
  $$
  and then $\check{M}(a,b)\cong \int^{x,y}P(P\otimes \C)(y,a,x,b)\otimes
  P^*J(y,x)$. In other words,  $\check{M}(a,b)\cong \tilde S(a,b)$;
  see \eqref{eq:tildeSmonoidalcat}.
  
  In conclusion, for a left autonomous map pseudomonoid in
  $\V\text-\mathbf{Mod}$, the 
  monads $\check{M}$ and $\tilde S$ are isomorphic. 
\end{obs}
\begin{exmp}
Let $\mathscr G$ be a groupoid. Write $\Delta:\mathscr G\to\mathscr
G\times\mathscr G$ for the diagonal functor and $E:\mathscr G\to 1$
the only possible functor. These give $\mathscr G$ a structure of
comonoid in $\mathbf{Cat}$ and thus $P=\Delta^*$ and $J=E^*$ is a
promonoidal structure on $\mathscr G$. Explicitly, $P(a,b;c)=\mathscr
G(a,c)\times\mathscr G(b,c)$ and $J(a)=1$; the monoidal structure
induced in $[\mathscr G,\mathbf{Set}]$ is given by the point-wise
cartesian product. Define a functor $D:\mathscr
G^{\mathrm{op}}\to\mathscr G$ as the identity on objects an
$D(f)=f^{-1}$ on arrows. In \cite[Example 10]{monbicat} was
essentially shown that $D$ is a left and right dualization for
$\mathscr G$. Then, by Corollary \ref{lc=c.cor}, $\mathscr G$ has
centre and lax centre in $\mathbf{Set}\text{-}\mathbf{Mod}$ and both
coincide. On the other hand, $[Z(\mathscr G),\mathbf{Set}]\simeq
Z([\mathscr G,\mathbf{Set}])$ by Theorem \ref{t:MZsimeqZM}, which
together with \cite[Theorem 4.5]{Day:LaxCentre} shows that the centre
of $\mathscr G$ in $\mathbf{Set}\text{-}\mathbf{Mod}$ is equivalent to
the category called (lax)centre of $\mathscr G$ in the latter article.
\end{exmp}

\section{Comodules}\label{comodules.sec}

This section deals with the case of the monoidal bicategories of
comodules 
$\mathbf{Comod}(\V)$. In general, \V\ will be a braided monoidal
category with a certain completeness condition. However, when we consider
the lax centre of pseudomonoids the braiding will be a symmetry. 
Our aim is to show how the general theory developed in previous
sections specialises to some of the most basic results of the theory of 
Hopf algebras. 

Throughout the section we will use {\em string diagrams}\/ to denote
arrows in \V. Our convention is that arrows go downwards: 
the domain of the arrow is the top of the string while
the codomain is the bottom string. Arrows are depicted as nodes with
the exception of the comultiplication of a comonoid, which is 
pictured as the bifurcation of one string into two. For background on
string diagrams see \cite{Joyal:GeoTensorCalculus}.

Given a monoidal category \V,   
there is a monoidal 2-category
$\mathbf{Comon}(\V)$ called the {\em 2-category of comonoids.} Its
objects are comonoids in \V, its 1-cells comonoid morphisms and
2-cells $\sigma:f\Rightarrow g:C\to D$ are arrows $\sigma:C\to I$ in
\V\ such that
$$
\xy
(-7,-7)*{\sigma};
(-7,-7)*\xycircle(2.65,2.65){-};
(0,0)**\dir{-};
(7,-7)*{f};
(7,-7)*\xycircle(2.65,2.65){-}="a";
(0,0)**\dir{-};
(0,0)*{};(0,7)*{}**\dir{-}?(.7)+(3,0)*{C};
"a";(13,-13)**\dir{-}?(.7)+(3,0)*{D};
\endxy
=
\xy
(-7,-7)*{g};
(-7,-7)*\xycircle(2.65,2.65){-}="a";
(0,0)**\dir{-};
(7,-7)*{\sigma};
(7,-7)*\xycircle(2.65,2.65){-};
(0,0)**\dir{-};
(0,0)*{};(0,7)*{}**\dir{-}?(.7)+(3,0)*{C};
"a";(-13,-13)**\dir{-}?(.7)-(3,0)*{D};
\endxy
$$
Vertical composition of 2-cells is the usual convolution product:
$\sigma * \sigma'= (\sigma\otimes\sigma')\Delta$, where $\Delta$
denotes the comultiplication. The horizontal compositions 
$$
\diagramcompileto{comod3}
A\rtwocell^f_g{\sigma}&B\ar[r]^h&C
\enddiagram
\qquad\text{and}\qquad
\diagramcompileto{comod4}
D\ar[r]^k&A\rtwocell^f_g{\sigma}&B
\enddiagram
$$
are $A\xrightarrow{\sigma}I$ and
$D\xrightarrow{k}A\xrightarrow{\sigma}I$ respectively.

\begin{obs}
$\mathbf{Comon}(\V)$ is the full sub 2-category of
    ${\V}^{\mathrm{op}}\text{-}\mathbf{Cat}^{\mathrm{op}}$ consisting of
    those $\V^{\mathrm{op}}$-categories with just one object. In
    particular, it is triequivalent to a strict 3-category.
\end{obs}

A pseudomonoid $(C,j,p)$ in $\mathbf{Comon}(\V)$ amounts to a comonoid
$C$ with two comonoid morphisms $j:I\to C$ and $p:C\otimes C\to C$ and 
the invertible 2-cells $\phi:p(p\otimes C)\Rightarrow p(C\otimes p)$,
$\lambda:p(j\otimes C)\Rightarrow 1$ and $\rho:p(C\otimes
j)\Rightarrow1$ satisfying axioms. These 2-cells are
convolution-invertible 
arrows $\phi:C\otimes C\otimes C\to I$ and $\lambda,\rho:C\to I$. 

\begin{exmp}
Normal pseudomonoids, that is,
pseudomonoids whose unit constraints $\lambda, \rho$ are identities,
in the monoidal bicategory  $\mathbf{Comod}(\mathbf{Vect})$ are {\em
coquasibialgebras.} The dual of this algebraic structure, called {\em
quasibialgebra,} was first defined in
\cite{Drinfeld:QuasiHopfAlg} where also were defined the quasi-Hopf
algebras. Then, a coquasibialgebra amounts to a coalgebra
$(C,\epsilon,\Delta)$ with a multiplication $p:C\otimes C\to C$,
denoted by $p(x\otimes y)=x\cdot y$, a
unit $j\in C$, where $k$ is the field, and an additional functional
$\phi:C\otimes C\otimes C\to k$ satisfying
$p(j\otimes x)=x=p(x\otimes j)$, 
$$
\sum \phi(x_1\otimes y_1\otimes z_1)(x_2\cdot y_2)\cdot
z_2=\sum x_1\cdot(y_1\cdot z_1)\phi(x_2\otimes y_2\otimes z_2)
$$
$$
(\phi\otimes \epsilon)*\phi(1\otimes p \otimes1)*(\epsilon\otimes
\phi)= \phi(p\otimes1\otimes 1)*\phi(1\otimes1\otimes p)
$$
where $*$ denotes the convolution product in the dual of $C\otimes
C\otimes C\otimes C$. We used Sweedler's notation $\Delta(x)=\sum
x_1\otimes x_2\in C\otimes C$, as is usual in the theory of Hopf algebras.
\end{exmp}

Now suppose further that $\V$ has equalizers of reflexive pairs and
each functor $X\otimes -$ preserves them. Then we can construct the
{\em bicategory of comodules}\/ over \V, denoted by
$\mathbf{Comod}(\V)$. It has comonoids in \V\ as objects and homs
$\mathbf{Comod}(\V)(C,D)$ the category of
$C$-$D$-bicomodules; this is the category of Eilenberg-Moore algebras
for the comonad $C\otimes-\otimes D$ on \V. The composition of two
comodules $M:C\to D$ and $N:D\to E$ is given by the equalizer of the
following reflexive pair
$$
\xymatrixcolsep{1.5cm}
\diagramcompileto{comod5}
{M\square_{D}N}\ar@{->}[r]&M\otimes N\ar@<1ex>[r]^-{\chi_{r}^M\otimes
  N}\ar@<-1ex>[r]_-{M\otimes\chi_\ell^N}& 
M\otimes D\otimes N
\enddiagram
$$
where the various $\chi$ denote the obvious coactions, and 
with $C$-$E$-comodule structure induced by the structures of $M$ and
$N$. The comodule $M\square_DN$ is sometimes called the {\em cotensor
product}\/ of $M$ and $N$ over $D$.
The identity 1-cell corresponding to a comonoid $C$ is the {\em
regular comodule} $C$, {\em i.e.}\/ it is $C$ with coaction
$(\Delta\otimes1)\Delta:C\to C\otimes C\otimes C$.

There is a pseudofunctor
$(-)_*:\mathbf{Comod}(\V)\to\mathbf{Comod}(\V)$ acting
as the identity on objects, sending a comonoid morphism $f:C\to D$ to
the comodule, denoted by $f_*:C\to D$, with underlying object $C$ and
coaction 
$$
\xy
(0,0)*{};(0,-13)*{}**\dir{-};
(-13,-13)*{};(0,0)**\dir{-};
(7,-7)*{f};
(7,-7)*\xycircle(2.65,2.65){-}="a";
(0,0)**\dir{-};
(0,0)*{};(0,7)*{}**\dir{-}?(.7)+(3,0)*{C};
"a";(13,-13)**\dir{-}?(.7)+(3,0)*{D};
\endxy
$$
and sending a 2-cell $\sigma:f\Rightarrow g$ to the comodule morphism
$\sigma_*:f_*\Rightarrow g_*$  given by 
$$
\xy
(-7,-7)*{};(0,0)**\dir{-};
(7,-7)*{\sigma};
(7,-7)*\xycircle(2.65,2.65){-};
(0,0)**\dir{-};
(0,0)*{};(0,7)*{}**\dir{-}?(.7)+(3,0)*{C};
\endxy
$$
The axioms of coaction and of comodule morphism follow from the ones
of comodule morphism and 2-cell in $\mathbf{Comon}(\V)$ respectively.
It is easy to show that the pseudofunctor  $(-)_*$ is locally fully
faithful (in fact, locally  it can be viewed as a
$\V^{\mathrm{op}}$-enriched Yoneda embedding).

An important property of $(-)_*$ is that it sends any 1-cell in
$\mathbf{Comon}(\V)$ to a map in
$\mathbf{Comod}(\V)$. For, if $f:C\to D$ is a comonoid morphism, then
$f_*$ has a right adjoint, denoted by  $f^*$, with underlying object
$C$ and coaction 
$$
\xy
(0,0)*{};(0,-13)*{}**\dir{-};
(13,-13)*{};(0,0)**\dir{-};
(-7,-7)*{f};
(-7,-7)*\xycircle(2.65,2.65){-}="a";
(0,0)**\dir{-};
(0,0)*{};(0,7)*{}**\dir{-}?(.7)+(3,0)*{C};
"a";(-13,-13)**\dir{-}?(.7)+(-3,0)*{D};
\endxy
$$
The composition $f_*f^*$ is the comodule with object $C$ and coaction $$
\xy
(7,-7)*{f};
(7,-7)*\xycircle(2.65,2.65){-}="b";
(0,0)*{}**\dir{-};
"b";(13,-13)*{}**\dir{-}?(.7)+(3,0)*{D};
(-7,-7)*{f};
(-7,-7)*\xycircle(2.65,2.65){-}="a";
(0,0)**\dir{-};
(0,0)*{};(0,7)*{}**\dir{-}?(.7)+(3,0)*{C};
"a";(-13,-13)**\dir{-}?(.7)+(-3,0)*{D};
(0,0)*{};(0,-13)*{}**\dir{-}?(.9)+(2,0)*{C}
\endxy
$$
and the counit of the adjunction is just the arrow $f:C\to D$, which turns
out to be a comodule morphism; the
unit is the unique map such that 
$$
\diagramcompileto{comod6}
f^*f_*=f_*\square_D f^*\ar[r]&C\otimes C\\
C\ar[ur]_\Delta\ar@{..>}[u]^\eta
\enddiagram
$$
where the horizontal arrow is the defining equalizer of $f^*f_*$.

\begin{obs}\label{comodhasEMcomonads.o}
  The bicategory $\mathbf{Comod}(\V)$ has Eilenberg-Moore objects for
  comonads. If $G$ is a comonad on the comonoid $C$ with
  comultiplication $\delta:G\to G\square_C G$ and counit
  $\epsilon:G\to C$, its Eilenberg-Moore object admits the following
  description (which is dual to the description of Kleisli objects for
  monads in $\V\text{-}\mathbf{Cat}$ in \cite{Street:CauchyChar}). As a
  comonoid, it is $G$ equipped with comultiplication and counit
$$
G\xrightarrow{\delta}G\square_CG\rightarrowtail G\otimes
G\qquad\text{and}\qquad G\xrightarrow{\epsilon}C\xrightarrow{\varepsilon}I.
$$
Note that the arrow $\epsilon :G\to C$ in \V\ becomes a morphism of
comonoids. The universal 1-cell is just the comodule
$\epsilon_*:G\to C$.
\end{obs}

\begin{obs}
The bicategory $\mathbf{Comod}(\V)$ can be viewed as the full sub bicategory of
  $\V^{\mathrm{op}}\text{-}\mathbf{Mod}^{\mathrm{coop}}$ determined by
  the $\V^{\mathrm{op}}$-categories with one object. However, for 
  $\V^{\mathrm{op}}\text{-}\mathbf{Mod}$ to exist
  further completeness assumptions on \V\ have to be made.

\end{obs}

When \V\ is braided, $\mathbf{Comon}(\V)$ and $\mathbf{Comod}(\V)$
have the structure of monoidal
2-categories with tensor product given by the
tensor product of $\V$; note that the braiding is used in defining the
comultiplication and coactions on the tensor product of comonoids and
comodules.
The pseudofunctor $(-)_*$ is strong monoidal, so
that through $(-)_*$ we can think of $\mathbf{Comon}(\V)$ as a
monoidal sub bicategory of $\mathbf{Comod}(\V)$. 
Since its  tensor product is a 2-functor, by
\cite{tricategories}, $\mathbf{Comod}(\V)$ 
is triequivalent to a strict 3-category.

\begin{figure}\label{fig:coevandev}
\begin{center}
\includegraphics[scale=1.2]{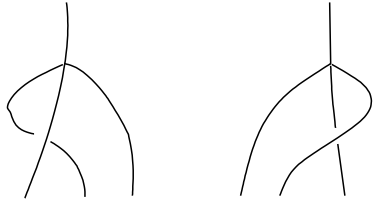}
\end{center}
\caption{Comodule structures of the coevaluation and evaluation.}
\end{figure}

The bicategory $\mathbf{Comod}(\V)$ is not just monoidal but it is
also left and right autonomous. The right bidual of a
comonoid $C$ is the opposite comonoid $C^\circ$. The braiding provides
pseudonatural equivalences 
$$
\mathbf{Comod}(\V)(C\otimes
D,E)\simeq\mathbf{Comod}(\V)(D,C^\circ\otimes E).
$$
The coevaluation
 $\n:I\to C^{\circ}\otimes C$ and evaluation $\e:C\otimes C^\circ\to
 I$ comodules are  the object $C$ with coactions depicted in Figure
 \ref{fig:coevandev}. 
The left bidual is defined by using the inverse of the braiding.

\begin{exmp}
As shown in \cite{dualizations}, Coquasi-Hopf algebras are exactly
the left autonomous normal pseudomonoids in
$\mathbf{Comod}(\mathbf{Vect})$ whose unit, multiplication and dualization are
representable by coalgebra morphisms.
\end{exmp}


\subsection{Hopf modules}

From now on, \V\ will not only have equalizers of reflexive pairs, but
all equalizers.  The reason for this is that we want the following
proposition to hold.
Equalizers are necessary as the proof uses the Adjoint Triangle
Theorem \cite{Dubuc:AdjTri}.

\begin{prop}[\cite{dualizations}]\label{p:MmapiffeMmap}
  A comodule $M:C\to D$ has a right adjoint if and only if its
  composition with $\varepsilon_*:D\to I$ has a right adjoint. 
\end{prop}

We shall describe the monad $\theta$ for a map pseudomonoid $(C,j,p)$,
which for
simplicity we will suppose arising from a pseudomonoid in
$\mathbf{Comon}(\V)$. 

\begin{figure}\label{fig:cotimesmV2}
\begin{center}
\includegraphics{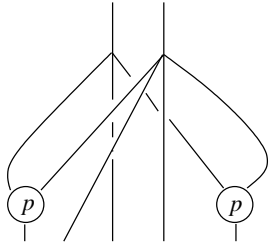}
\caption{Comodule structure of $\theta_D(M)$.}
\end{center}
\end{figure}

Recall from Definition \ref{d:theta} 
that the monad $\theta_D$ on $\mathbf{Comod}(C\otimes D,C)$ is
just the monad $\mathbf{Kl}(C\otimes -)(D,p)$ on $\mathbf{Kl}(C\otimes
-)(D,C)$. In terms of comodules, $\theta_D(M)$ has underlying object
$C\otimes M$ and coaction the arrow depicted in Figure \ref{fig:cotimesmV2}.
The multiplication \eqref{mult.eq} and the unit
\eqref{unit.eq} become in this case 
$$
C\otimes C\otimes M\to C^{\otimes 3}\otimes D \otimes
  C\otimes M \otimes C^{\otimes 3}\xrightarrow{\phi^{-1}\otimes
    \varepsilon_D\otimes 1\otimes 1\otimes \phi}C\otimes C\otimes
  M\xrightarrow{p\otimes 1}C\otimes M
$$
and
$$
M\to C\otimes D\otimes M\otimes
C\xrightarrow{\lambda\otimes\varepsilon_D\otimes 1 \otimes\lambda^{-1}}
M\xrightarrow{j\otimes 1}C\otimes M
$$
where $\phi^{-1}$ and $\lambda^{-1}$ are the convolution
inverses of $\phi:C\otimes C\otimes C\to I$ and $\lambda:C\to I$.

\begin{figure}\label{fig:thetaN}
\begin{center}
\includegraphics{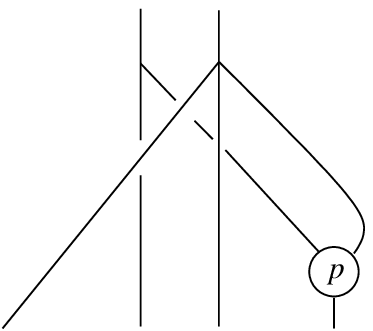}
\end{center}
\caption{}
\end{figure}

In view of Observation \ref{upsilonlambda.obs},
the functor $\upsilon_D\lambda_D$ (see Definition \ref{lambda.def})
is isomorphic to the one sending a comodule
$M:D\to C$  to $C\otimes M$ with coaction the arrow in Figure
\ref{fig:thetaN}. 
The theorem of Hopf modules holds for $C$ exactly when every
$\theta_D$-algebra is isomorphic to one of this form.

\begin{figure}\label{fig:tV2}
\begin{center}
\includegraphics{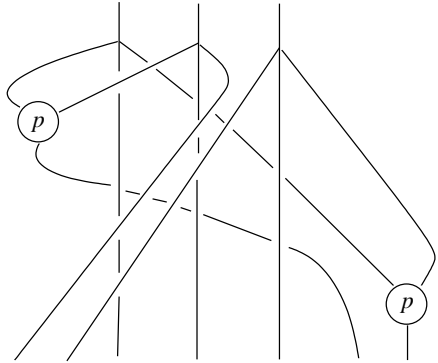}
\end{center}
\caption{Comodule structure of the monad $t$.}
\end{figure}

Now we shall describe for a pseudomonoid $C$ in $\mathbf{Comon}(\V)$
the underlying comodule of the monad $t$ on
$C^\circ\otimes C$ representing $\theta$. Recall from
\eqref{eq:twhenbiduals} that
$$
t\cong (C^\circ\otimes p_*)(C^\circ\otimes C\otimes\e\otimes C)(C^\circ\otimes
p^*\otimes C^*\otimes C)(\n\otimes C^*\otimes C)$$
and so it has underlying object $C\otimes C\otimes C$ with coaction
depicted in Figure \ref{fig:tV2}.

The Hopf module construction for a map pseudomonoid in
$\mathbf{Comod}(\V)$ may not exist,
as this bicategory does
not have an Eilenberg-Moore objects for monads. However, it does have
Eilenberg-Moore construction for comonads (Observation
\ref{comodhasEMcomonads.o}).

\begin{prop}\label{CdualHMC.p}
Given a map pseudomonoid $C$ in $\mathbf{Comod}(\V)$, if the monad
$t:C^\circ\otimes C\to C^\circ\otimes C$ has right adjoint, then $C$ has a
Hopf module construction. In particular, this holds  if
$C\in\mathrm{ob}\V$ has a dual.
\end{prop}
\begin{proof}
The 1-cell $t^*$ has a canonical structure of a right adjoint comonad
to the monad $t$. It well-known that the Eilenberg-Moore construction
for the comonad $t^*$ is an Eilenberg-Moore construction for the monad
$t$. 
To finish, we show that if $C$ has a dual in \V\ then
$t\cong((p^*)^\circ\otimes p)(C^\circ\otimes\n\otimes C)$ has a right
adjoint, and for that it suffices to prove that $\n$ does. But by
Proposition \ref{p:MmapiffeMmap}, $\n$ is a map if and only if $C$ has a
dual. 
\end{proof}

When $\V$ is the category of vector spaces and 
$C$ is a coquasi-bialgebra, the assertion that the functor
$\lambda_I$ from $\mathbf{Comod}(\V)(I,C)$ to the category of Hopf
modules is an equivalence is what Schauenburg
\cite{Schauenburg:TwoCharacterizations} calls the theorem of Hopf
modules. We shall show that when $C$ has a Hopf module construction
both notions are equivalent. 

  Let $\mathscr W$ be a braided monoidal replete full subcategory of
  \V\ closed under equalizers of reflexive pairs. There is an inclusion
  monoidal pseudofunctor $\mathbf{Comod}(\mathscr
  W)\to\mathbf{Comod}(\V)$. This inclusion, being monoidal, preserves
  biduals. 

\begin{cor}\label{comodWcomodV.p}
Let $\mathscr W$ and \V\ be as above. Suppose $C$ is a map
pseudomonoid in $\mathbf{Comod}(\mathscr W)$ such that $C$ has a dual
in $\mathscr W$. Then, the theorem of Hopf modules holds for $C$ in
$\mathbf{Comod}(\mathscr W)$ if and only if it holds for $C$ in
$\mathbf{Comod}(\V)$. 
\end{cor}
\begin{proof}
We begin by observing that since $C$ has dual in $\mathscr W$, and
hence in \V, by Proposition \ref{CdualHMC.p}, $C$ has a Hopf module
construction both
in $\mathbf{Comod}(\mathscr W)$ and in $\mathbf{Comod}(\V)$. Moreover,
the two coincide. To see this, observe that the monad $t$ is given by
\eqref{eq:twhenbiduals} and each of the 1-cells in the composition
lies in $\mathbf{Comod}(\mathscr W)$. Since $C$ has a dual, $t$ has a
right adjoint comonad, whose Eilenberg-Moore construction, described
in Observation \ref{comodhasEMcomonads.o}, is the Hopf module
construction for $C$. By the description of this Eilenberg-Moore
construction, one sees that it lies in $\mathbf{Comod}(\mathscr W)$.  

Hence, we have to prove that the 
1-cell $\ell:C\to(C^\circ\otimes
C)^t$ (see Proposition \ref{lffstrongmon.prop}) is an equivalence in
$\mathbf{Comod}(\mathscr W)$ if and only if it is one in
$\mathbf{Comod}(\V)$.
One direction is trivial, so we shall suppose $\ell$ is an equivalence in
$\mathbf{Comod}(\V)$. We have, then, an adjoint equivalence
$\ell\dashv\ell^*$; as $\ell$ is always a map (by Proposition
\ref{lffstrongmon.prop}), this adjoint equivalence lifts to
$\mathbf{Comod}(\mathscr W)$.
\end{proof}

In the particular case when \V\ is the category of vector spaces and
$\mathscr W$ is the subcategory of finite-dimensional vector spaces, we
have:
\begin{cor}
For any finite-dimensional coquasi-bialgebra $C$ there exists a map
pseudomonoid $D$ in $\mathbf{Comod}(\mathbf{Vect})$ such that the
category of Hopf modules for $C$ (as defined in
\cite{Schauenburg:TwoCharacterizations}) is monoidally equivalent to  the
category of right $D$-comodules $\mathbf{Comod}(\mathbf{Vect})(I,D)$. Moreover,
$D$ can be taken to be the Hopf module construction for $C$, and in
particular, finite-dimensional.
\end{cor}
Note that, in general, the forgetful functor
$\mathbf{Comod}(\mathbf{Vect})(I,D)\to\mathbf{Vect}$ is not monoidal.

By Observation \ref{comodhasEMcomonads.o}, the Hopf module
construction $(C^\circ\otimes C)^t\to C^\circ\otimes C$ can be taken to be 
of the form $\epsilon_*$, where $\epsilon:(C^\circ\otimes C)^t\to
C^\circ\otimes C$ is a comonoid morphism. 

\begin{cor}\label{c:thmeqlamdaIeq}
  Suppose that $C$ is a map
  pseudomonoid  in $\mathbf{Comod}(\mathbf{Vect})$. If $C$ is
  finite-dimensional, the theorem of Hopf
  modules holds for $C$ if and only if the functor 
$$
\lambda_I:\mathbf{Comod}(\mathbf{Vect})(I,C)\to\mathbf{Comod}(\mathbf{Vect})(C,C)^{\theta_I} 
$$
(see Definition \ref{lambda.def}) is an equivalence.
\end{cor}
\begin{proof}
Only the converse is non trivial.
Write $\V$ for $\mathbf{Vect}$.
By Proposition \ref{comodWcomodV.p}, it is enough to show that the
theorem of Hopf modules holds for $C$ in $\mathbf{Comod}(\V_f)$.

The functor $\lambda _I$ is represented by the 1-cell $\ell:C\to
(C^\circ\otimes C)^t$. We have that the functor
$\mathbf{Comod}(\V_f)(I,\ell)$ is an equivalence, and the result follows
from the following lemma.  
\end{proof}

\begin{lem}\label{comodI.l}
\begin{enumerate}
\item
The functor $\mathbf{Comod}(\mathbf{Vect}_f)(I,-)$ reflects equivalences.
\item
Any finitely continuous functor from 
the category 
$\mathbf{Comod}(\mathbf{Vect}_f)(I,D)$ to $\mathbf{Comod}(\mathbf{Vect}_f)(I,E)$ is isomorphic to a functor 
given by composition with a comodule $M:D\to E$.
\end{enumerate}
\end{lem}
\begin{proof}
(1)  Assume that $\mathbf{Comod}(\mathbf{Vect}_f)(I,M)$ is an
  equivalence, for $M:D\to 
  E$. Taking duals, we see that the functor from
  $D^*\text{-}\mathrm{Mod}_f$ to $E^*\text{-}\mathrm{Mod}_f$ given by
  tensoring with $M^*$ is an equivalence; this implies that $M^*$ is
  an equivalence in $\mathbf{Mod}(\mathbf{Vect}_f)$ (a Morita
  equivalence), and hence, taking duals, $M$ is an equivalence. 

(2) Taking duals, the result follows from the fact that any finitely
cocontinuous functor between categories of finite-dimensional modules
over a finite-dimensional algebra, is isomorphic to a functor induced
by tensoring with a bimodule.
\end{proof}

We obtain the following  generalisation of 
\cite[Thm. 3.1]{Schauenburg:TwoCharacterizations}.
\begin{cor}
Let $C$ be a map pseudomonoid in $\mathbf{Comod}(\mathbf{Vect})$ whose
underlying 
space is finite-dimensional. Then $C$  has a left dualization if and only
if the functor $\lambda_I:\mathbf{Comod}(\mathbf{Vect})(I,C)\to
\mathbf{Comod}(\mathbf{Vect})(C,C)^{\theta_I}$ is an equivalence. 
\end{cor}
\begin{proof}
By the corollary above, the theorem of Hopf modules holds for $C$;
hence, $C$ has a left dualization by Theorem \ref{maintheorem.th}.
\end{proof}

\begin{figure}\label{fig:F}
\begin{center}
\includegraphics{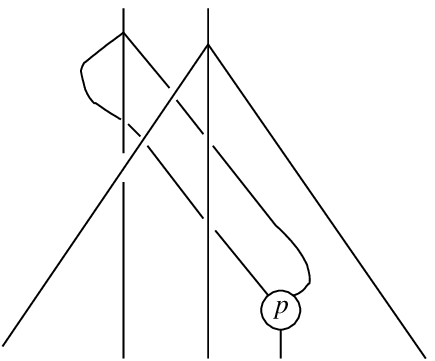}
\end{center}
\caption{}
\end{figure}

Now suppose that $C$ is a left autonomous map pseudomonoid in
$\mathbf{Comod}(\V)$. The existence of a left dualization forces
the multiplication to be a map \cite[Prop. 1.2]{dualizations}. On the
other hand, 
the unit of $C$ is a map because its underlying object $I\in \V$ has
(right) dual by Proposition \ref{p:MmapiffeMmap}, it follows
that any left autonomous pseudomonoid in $\mathbf{Comod}(\V)$ 
is a map pseudomonoid. 
A Hopf module construction for $C$ is provided by 
$(C^\circ\otimes p)(\n\otimes C)\cong
 (p(d\otimes
C))^*:C\to C^\circ\otimes C$. 
In the case when $C$ is a coquasibialgebra, 
the comodule $(C^\circ \otimes p_*)(\n\otimes C)$ is $C\otimes C$ with
coaction depicted in Figure \ref{fig:F}.

\subsection{Centre and Drinfel'd double}

We now consider the results of Section \ref{centre.s} on the lax
centre in the context of comodules. 
We suppose the underlying monoidal category \V\ is symmetric, and thus
$\mathbf{Comon}(\V)$ is a symmetric monoidal $\mathbf{Cat}$-enriched
category. Via the monoidal pseudofunctor $(-)_*$ we obtain comodules 
$c_{M,N}:M\otimes N\to N\otimes M$ making the usual diagrams commute
up to canonical isomorphisms. 

\begin{prop}\label{p:lauthaslaxcentrComod}
Any left autonomous pseudomonoid in $\mathbf{Comod}(\V)$ whose
underlying object in \V\ 
has dual has a lax centre. If the pseudomonoid 
is also right autonomous then the lax
centre equals the centre.
\end{prop}
\begin{proof}
We have already mention that any left autonomous pseudomonoid $C$ in
$\mathbf{Comod}(\V)$ is a map pseudomonoid.  
By Corollary \ref{laxcentreEMconstr.cor} we have to show that the
monad $\tilde s:A\to A$ has an Eilenberg-Moore construction, and for
that it is enough to show that it has a right adjoint, since
$\mathbf{Comod}(\V)$ has Eilenberg-Moore objects for comonads. 
Again by Corollary \ref{laxcentreEMconstr.cor}, we have $\tilde s\cong
p(p\otimes C)(C\otimes c_{C,C})(p^*\otimes C)(j\otimes C)$ and
therefore $\tilde s$ has a right adjoint if $p^*j:I\to C\otimes C$ has one; 
but $C$ being left autonomous, this 1-cell is
isomorphic to $(d\otimes C)\n$ which is a composition of maps: $d$ by
\cite[Prop. 1.2]{dualizations} and $\n$ by Proposition \ref{p:MmapiffeMmap}.
\end{proof}
\begin{obs}\label{centreinVf.o}
  In the proposition above, suppose that the full subcategory $\V_f$ of objects
  with a dual in \V\ is closed under equalizers of reflexive
  pairs. Then the lax centre $Z_\ell(C)\to C$ lies in
  $\mathbf{Comod}(\V_f)$, and it is a lax centre in it. 

  To prove this observe that $t:C^\circ\otimes C\to C^\circ\otimes C$
  and its Eilenberg-Moore construction $C\to C^\circ\otimes C$ lie in
  $\mathbf{Comod}(\V_f)$, and the monad $s$ and the distributive law
  between $t$ and $s$ do so too. It follows that the induced monad
  $\tilde s$ on $C$ lies in $\mathbf{Comod}(\V_f)$, and it has right
  adjoint in this bicategory, as shown in the proof above, and it is
  necessarily the same as in $\mathbf{Comod}(\V)$.  It follows from
  Observation \ref{comodhasEMcomonads.o} that $\tilde s^*$ has an
  Eilenberg-Moore construction in $\mathbf{Comod}(\V_f)$ and coincides
  with the respective construction in $\mathbf{Comod}(\V)$. Moreover,
  this construction is given by $\epsilon_*:C\to C^{\tilde s^*}$,
  where $\epsilon$ is the counit of the comonad $\tilde s^*$.
  Therefore, the lax centre of $C$ in both bicategories coincide.
\end{obs}

The {\em Drinfel'd double}\/ or {\em quantum double}\/ of a
finite-dimensional Hopf
algebra is a finite-dimensional 
braided (also called quasitriangular) Hopf algebra $D(H)$
with underlying vector space $H^*\otimes H$ (one can also take $H\otimes H$)
and suitably defined structure. It is a classical  result
that the category of left $D(H)$-modules is monoidally equivalent to
the category of (two-sided) $H$-Hopf modules and to the centre of the
category of $H$-modules. 
The Drinfel'd double of
a finite-dimensional
{\em quasi-Hopf algebra}\/ was defined in 
\cite{Majid:QDQHA}
using a reconstruction theorem, and explicit constructions were given
in \cite{Hausser:DoublesQuasiQG,Schauenburg:HopfMods}. This last
paper shows that the category of $D(H)$-modules is monoidally equivalent to the
centre of the category of $H$-modules, via a generalisation of
the Yetter-Drinfel'd modules. 
The quantum double of a coquasi-Hopf algebra was described in
\cite{Bulacu:DrinfeldDoubleDiagCrossCopr}. Alternatively, it can be
described by dualising the explicit constructions for the quasi-Hopf
case. Then the Drinfel'd or quantum double
$D(H)$ of a finite-dimensional coquasi-Hopf $H$ algebra is
finite-dimensional and has the property that the category of
$D(H)$-comodules $\mathrm{Comod}(D(H))$ is monoidally equivalent to the
centre of $\mathrm{Comod}(H)$, and the equivalence commutes with the
forgetful functors. 


Given a finite-dimensional coquasi-Hopf algebra $H$, we would like to
study the relationship between the centre $Z(H)$ in
$\mathbf{Comod}(\mathbf{Vect})$ and the Drinfel'd double $D(H)$. To
this aim we will need some of the machinary of {\em Tannakian
  reconstruction}, of which we give the most basic aspects following
\cite{McCrudden:Maschkean}. 

Let \V\ be a monoidal category and $\V_f$ the full sub-monoidal
category with objects with left duals.  
We denote by $\V_f\text-\mathrm{Act}$ the
2-category of pseudoalgebras for the pseudomonoad $(\V_f\times-)$ on
$\mathbf{Cat}$. 
Objects of this 2-category are pseudoactions of $\V_f$ and
1-cells are pseudomorphisms of pseudoactions. Observe that $\V_f$ has
a canonical $\V_f$-pseudoaction given by the tensor product. 
We form the 2-category
$\V_f\text-\mathrm{Alg}/\V_f$ with objects 1-cells $\sigma:\A\to\V_f$ in
$\V_f\text-\mathrm{Act}$. The 1-cells are pairs
$(F,\phi):\sigma\to\sigma'$ where $F:\A\to \A'$ is a 1-cell in
$\V_f\text-\mathrm{Act}$ and $\phi:\sigma'F\cong\sigma$ is a 2-cell in
$\V_f\text-\mathrm{Act}$. 2- cells $(F,\phi)\Rightarrow(F',\phi')$ are
just 2-cells $F\Rightarrow F'$ in $\V_f\text-\mathrm{Act}$. There is a
2-functor
$\mathrm{Comod}_f:\mathbf{Comon}(\V)\to\V_f\text-\mathrm{Act}/\V_f$
sending a comonoid $C$ to the forgetful functor
$\omega_C:\mathrm{Comod}_f(C)\to\V_f$; here $\mathrm{Comod}_f(C)$ is
the category of right coactions of $C$ with underlying object in
$\V_f$. This category has a canonical $\V_f$-pseudoaction such that
$\omega$ is an object of $\V_f\text-\mathrm{Act}/\V_f$. The definition
of $\mathrm{Comod}_f$ on 1-cells and 2-cells should be more or less
obvious; see \cite{McCrudden:Maschkean}. 

Under certain hypothesis on \V, the 2-functor $\mathrm{Comod}_f$ is
bi-fully faithful. Here is the case we will need. 

\begin{prop}[\cite{McCrudden:Maschkean}]\label{p:Comodfbiff}
The 2-functor
$$
\mathrm{Comod}_f:
\mathbf{Comon}(\mathbf{Vect})\to\mathbf{Vect}_f\text-\mathrm{Act}/\mathbf{Vect}_f
$$ 
is bi-fully faithful. 
\end{prop}

\begin{thm}
For any finite-dimensional coquasi-Hopf alebra $H$, 
the coalgebras $H^{\tilde s^*}$ and $D(H)$ are equivalent coquasibialebras. 
Moreover, they are isomorphic as coalgebras. 
\end{thm}
\begin{proof}
By Observatoin \ref{centreinVf.o}, $H^{\tilde s^*}$ is a centre for
the pseudomonoid $H$ in $\mathbf{Comod}(\mathbf{Vect}_f)$. Hence
we have an equivalence in
$\mathbf{Vect}_f\text-\mathrm{Act}/\mathbf{Vect}_f$ from the forgetful
functor $\mathrm{Comod}_f(H^{\tilde s^*})\to\mathbf{Vect}_f$ to the
forgetful functor $Z(\mathrm{Comod}_f(H))\to\mathbf{Vect}_f$. On the
other hand, there is an equivalence from the latter to
$\mathrm{Comod}_f(D(H))\to\mathbf{Vect}_f$. In this way we get an
equivalence from $\mathrm{Comod}_f(H^{\tilde s^*})$ to
$\mathrm{Comod}_f(D(H))$ in
$\mathbf{Vect}_f\text-\mathrm{Act}/\mathbf{Vect}_f$. By Propositon
\ref{p:Comodfbiff} we have an equivalence $f:H^{\tilde s^*}\to D(H)$ in
$\mathbf{Comod}(\mathbf{Vect})$. That is, both coquasibialgebras are
equivalent. 

Note that the equivalence $\mathrm{Comod}_f(f)$, given by
corestriction along $f$, is just the functor
$\mathbf{Comod}(\mathbf{Vect}_f)(I,f_*)$, the composition with the
comodule $f_*:H^{\tilde s^*}\to D(H)$. Then $f_*$ is an equivalence in
$\mathbf{Comod}(\mathbf{Vect}_f)$ by Lemma \ref{comodI.l}. 
But this implies that $f$ is an
isomorphism in $\mathbf{Vect}_f$, as the counit of the adjunction
$f_*\dashv f^*$ is given by $f$ itself. Hence, $f$ is an isomorphism
of coalgebras. 
\end{proof}


\begin{thebibliography}{10}
\expandafter\ifx\csname url\endcsname\relax
  \def\url#1{\texttt{#1}}\fi
\expandafter\ifx\csname urlprefix\endcsname\relax\def\urlprefix{URL }\fi

\bibitem{variation}
R.~Betti, A.~Carboni, R.~Street, R.~Walters, Variation through enrichment, J.
  Pure Appl. Algebra 29 (1983) 109--127.

\bibitem{Bulacu:DrinfeldDoubleDiagCrossCopr}
D.~Bulacu, B.~Chiri{{t}}{\u{a}}, Dual {D}rinfeld double by diagonal crossed
  coproduct, Rev. Roumaine Math. Pures Appl. 47~(3) (2002) 271--294 (2003).

\bibitem{closedcatfun}
B.~Day, On closed categories of functors, in: Reports of the Midwest Category
  Seminar, IV, Lecture Notes in Mathematics, Vol. 137, Springer, Berlin, 1970,
  pp. 1--38.

\bibitem{dualizations}
B.~Day, P.~McCrudden, R.~Street, {Dualizations and antipodes}, Appl. Categ.
  Struct. 11~(3) (2003) 229--260.

\bibitem{Day:LaxCentre}
B.~Day, E.~Panchadcharam, R.~Street, Lax braidings and the lax centre, in: Hopf
  albebras and generalizations, vol. 441 of Contemp. Math., Amer. Math. Soc.,
  Providence, RI, 2007.

\bibitem{monbicat}
B.~Day, R.~Street, Monoidal bicategories and {H}opf algebroids, Adv. Math.
  129~(1) (1997) 99--157.


\bibitem{Day-Street:quantumcat}
B.~Day, R.~Street, Quantum categories, star autonomy, and quantum groupoids,
  in: Galois theory, Hopf algebras, and semiabelian categories, vol.~43 of
  Fields Inst. Commun., Amer. Math. Soc., Providence, RI, 2004, pp. 187--225.

\bibitem{Day:CentresMonCatsFun}
B.~Day, R.~Street, Centres of moniodal categories of functors, in: Categories
  in Algebra, Geometry and Mathematical Physics, No. 431 in Contep. Math.,
  Amer. Math. Soc., 2007, pp. 187--202.

\bibitem{Drinfeld:QuasiHopfAlg}
V.~G. Drinfel{\cprime}d, Quasi-{H}opf algebras, Algebra i Analiz 1~(6) (1989)
  114--148.

\bibitem{Dubuc:AdjTri}
E.~Dubuc, Adjoint triangles, in: Reports of the Midwest Category Seminar, II,
  Springer, Berlin, 1968, pp. 69--91.

\bibitem{tricategories}
R.~Gordon, A.~J. Power, R.~Street, Coherence for tricategories, Mem. Amer.
  Math. Soc. 117~(558) (1995) vi+81.

\bibitem{math.QA/9904164}
F.~Hausser, F.~Nill, {Integral Theory for Quasi-Hopf Algebras}, (preprint)
  arXiv:math.QA/9904164.

\bibitem{Hausser:DoublesQuasiQG}
F.~Hausser, F.~Nill, Doubles of quasi-quantum groups, Comm. Math. Phys. 199~(3)
  (1999) 547--589.

\bibitem{Hermida:FromCoherent}
C.~Hermida, From coherent structures to universal properties, J. Pure Appl.
  Algebra 165~(1) (2001) 7--61.

\bibitem{Joyal:GeoTensorCalculus}
A.~Joyal, R.~Street, The geometry of tensor calculus. {I}, Adv. Math. 88~(1)
  (1991) 55--112.

\bibitem{braidedtencat}
A.~Joyal, R.~Street, {Braided tensor categories}, Adv. Math. 102~(1) (1993)
  20--78.

\bibitem{Kelly:BCECT}
G.~M. Kelly, Basic concepts of enriched category theory, vol.~64 of London
  Mathematical Society Lecture Note Series, Cambridge University Press,
  Cambridge, 1982.

\bibitem{Kelly-Lack:MonFunctors}
G.~M. Kelly, S.~Lack, Monoidal functors generated by adjunctions, with
  applications to transport of structure, in: Galois theory, Hopf algebras, and
  semiabelian categories, vol.~43 of Fields Inst. Commun., Amer. Math. Soc.,
  Providence, RI, 2004, pp. 319--340.

\bibitem{Lack:Acoherentapproach}
S.~Lack, {A coherent approach to pseudomonads}, Adv. Math. 152~(2) (2000)
  179--202.

\bibitem{Lack-Street:ftm2}
S.~Lack, R.~Street, The formal theory of monads. {II}, J. Pure Appl. Algebra
  175~(1-3) (2002) 243--265, special volume celebrating the 70th birthday of
  Professor Max Kelly.

\bibitem{Larson-Sweedler:AnAssocOrtogonal}
R.~G. Larson, M.~E. Sweedler, An associative orthogonal bilinear form for
  {H}opf algebras, Amer. J. Math. 91 (1969) 75--94.

\bibitem{Majid:QDQHA}
S.~Majid, Quantum double for quasi-{H}opf algebras, Lett. Math. Phys. 45~(1)
  (1998) 1--9.

\bibitem{Marmolejo:Doctrines}
F.~Marmolejo, Doctrines whose structure forms a fully faithful adjoint string,
  Theory Appl. Categ. 3 (1997) No.\ 2, 24--44 (electronic).

\bibitem{McCrudden:opmonoidalmonads}
P.~McCrudden, {Opmonoidal monads}, Theory Appl. Categ. 10 (2002) 469--485.

\bibitem{McCrudden:Maschkean}
P.~McCrudden, Tannaka duality for {M}aschkean categories, J. Pure Appl. Algebra
  168~(2-3) (2002) 265--307, category theory 1999 (Coimbra).

\bibitem{Street:DoublesMonCatsPreprint}
C.~Pastro, R.~Street, Doubles for monoidal categories, preprint (2007).

\bibitem{Schauenburg:HopfMods}
P.~Schauenburg, Hopf modules and the double of a quasi-{H}opf algebra, Trans.
  Amer. Math. Soc. 354~(8) (2002) 3349--3378 (electronic).

\bibitem{Schauenburg:TwoCharacterizations}
P.~Schauenburg, Two characterizations of finite quasi-{H}opf algebras, J.
  Algebra 273~(2) (2004) 538--550.

\bibitem{street.ftm}
R.~Street, {The formal theory of monads}, J. Pure Appl. Algebra 2 (1972)
  149--168.

\bibitem{Street:CauchyChar}
R.~Street, Cauchy characterization of enriched categories, Rend. Sem. Mat. Fis.
  Milano 51 (1981) 217--233 (1983).

\bibitem{Street:monoidalcentre}
R.~Street, The monoidal centre as a limit, Theory Appl. Categ. 13 (2004) No.
  13, 184--190 (electronic).

\bibitem{yonedastructures}
R.~Street, R.~Walters, {Yoneda structures on 2-categories}, J. Algebra 50
  (1978) 350--379.

\end{thebibliography}

\def\cprime{$'$}

\end{document}